\documentclass[12pt]{article}
\usepackage{graphicx}
\usepackage{amsmath}
\usepackage{geometry}
\usepackage{amsfonts}
\usepackage{amssymb}
\newtheorem{theorem}{Theorem}

\newtheorem{condition}[theorem]{Condition}

\newtheorem{corollary}[theorem]{Corollary}

\newtheorem{lemma}[theorem]{Lemma}

\newtheorem{proposition}[theorem]{Proposition}
\newtheorem{remark}[theorem]{Remark}

\newenvironment{proof}[1][Proof]{\textbf{#1.} }{\ \rule{0.5em}{0.5em}}

\begin{document}

\begin{center}
A Shooting Approach to Layers and Chaos in a Forced Duffing Equation, I

\bigskip

Shangbing Ai and Stuart P. Hastings\newline 

Department of Mathematics, University Pittsburgh, Pittsburgh, PA 15260
\end{center}

\section{abstract}

We study equilibrium solutions for the problem
\begin{align*}
u_{t}  &  =\varepsilon^{2}u_{xx}-u^{3}+\lambda u+\cos x\\
u_{x}\left(  0,t\right)   &  =u_{x}\left(  1,t\right)  =0.
\end{align*}
Using a shooting method we find solutions for all non-zero $\varepsilon.$
\ For small $\varepsilon$ \ we add to the solutions found by previous authors,
especially Angennent, Mallet-Paret and Peletier, and by Hale and Sakamoto, and
also give new elementary ode proofs of their results. \ \ Among the new
results is the existence of internal layer-type solutions. \ \ Considering the
ode satisfied by equilibria, but on an infinite interval, we obtain chaos
results for $\lambda\geq\lambda_{0}=\frac{3}{2^{2/3}}$ and $0<\varepsilon
\leq\frac{1}{4}.$ \ We also consider the problem of bifurcation of solutions
as $\lambda$ increases from $0.$

\section{\label{sec1}Introduction}

This is the first of a series of papers studying the existence of bounded
solutions of the equation%

\begin{equation}
\varepsilon^{2}u^{\prime\prime}=u^{3}-\lambda u+\cos t, \label{1.1}%
\end{equation}
where $\lambda$ \ and $\varepsilon$ \ are positive parameters. \ This equation
is a standard model in the theory of nonlinear oscillations, one of several
which have been called a forced Duffing equation in the literature. \ With the
signs shown it is often referred to as the equation of a ``soft spring.'' \ A
comprehensive reference to the early theory of $\left(  \ref{1.1}\right)  $ is
\cite{nm}, which gives a detailed account of the results obtained by classical
perturbation methods, such as averaging or multi-scale techniques. \ More
recent efforts have used dynamical systems concepts to establish results about
chaotic behavior of one sort or another. \ An important reference is by
Angenent, Mallet-Paret, and Peletier \cite{ampp}, who studied stable steady
states for a reaction-diffusion equation
\begin{equation}
u_{\tau}=\varepsilon^{2}u_{xx}+f\left(  x,u\right)  \label{1.2}%
\end{equation}
with boundary conditions
\begin{equation}
u_{x}\left(  0,\tau\right)  =u_{x}\left(  L,\tau\right)  =0 \label{1.3}%
\end{equation}
for a class of functions $f$ which were cubic in the state variable $u.$
\ Equation $\left(  \ref{1.1}\right)  $ is obtained from $\left(
\ref{1.2}\right)  $ (with $f\left(  x,u\right)  =u^{3}-\lambda u+\cos x)$ by
setting $x=t$ and assuming that $u$ \ is independent of $\tau.$ \ While the
specific form $\left(  \ref{1.1}\right)  $ is not mentioned in \cite{ampp},
\ the methods there are easily seen to apply and to prove important results
about the existence of stable steady-states of the problem $\left(
\ref{1.2}\right)  -\left(  \ref{1.3}\right)  $ with this function $f.$
\ Compactness requirements, however, seem to make it more difficult to study
the problem on an infinite interval, and thereby get ``chaos'', using infinite
dimensional methods. \ \ 

\bigskip

In this paper we introduce a shooting technique, which we have not seen in
this form elsewhere, which is the basis of our approach. \ We use this
technique to obtain new periodic solutions for $\left(  \ref{1.1}\right)  $
over a range of $\varepsilon$ which is not ``small''. \ We show that one can
give rigorous results about a weak form of ``chaos'' over this larger range of
$\varepsilon.$ \ We also introduce the problem of determining how these
solutions arise as $\lambda$ varies, since for $\lambda\leq0$, $\left(
1.1\right)  $ will be shown to have a unique bounded solution. \ Then, letting
$\varepsilon$ return to its traditional role as a small parameter, we
reproduce and extend the results in \cite{ampp} as they apply to this
equation, using elementary ode methods. \ \ We study the bifurcation problem
in $\lambda$ in more detail, \ and we begin our study of a further class of
bounded solutions which we can find for small $\varepsilon.$ \ 

\bigskip

However, a detailed study of this new class of solutions has proven rather
lengthy, so we will leave its elaboration to a later paper. \ The goal of that
paper will be to obtain a general chaos result which includes all of the
solutions discussed in each paper, in some sense. \ It is anticipated that
there will be still a third paper devoted to this topic, in which we use
techniques from \cite{hm2} to study the stability of the new solutions, as
equilibria of $\left(  \ref{1.2}\right)  -\left(  \ref{1.3}\right)  .$
\ \ This will lead to a consideration of the Morse index of these solutions,
and their role as part of the global attractor for $\left(  \ref{1.2}\right)
-\left(  \ref{1.3}\right)  .$ \ \ Future work will include extension of the
results to more general equations. \ 

\bigskip

In most of this paper we will consider the ode $\left(  \ref{1.1}\right)  $
for its own sake. \ The literature is large, and it is surprising that, as far
as we have been able to determine, many of the solutions we will find have not
been described before. \ Some of these solutions are periodic and some are
``chaotic'' in a sense to be defined below. \ Those which are periodic are
unstable as steady-states of $\left(  \ref{1.2}\right)  -\left(
\ref{1.3}\right)  $ and therefore apparently not found by the pde techniques
of \cite{ampp}.\ However recent work of K. Nakashima on a similar problem uses
continuation methods which may yield results similar to some of ours for this
problem, such as the existence of these new periodic solutions \cite{nak}.
There have been many other papers on the problem $\left(  \ref{1.2}\right)  -
\left(  \ref{1.3}\right)  $ with cubic nonlinearities, and some of these also
consider unstable solutions. Two recent ones, which cite other related work,
are by Hale and Salazar \cite{hs1},\cite{hs2}. An earlier paper by Kurland
gives many oscillatory solutions which are similar to some of those we find in
section 4.2, but for a different class of nonlinear functions $f(x,u)$
\cite{kur}. The use of formal asymptotic analysis on similar problems has been
studied by Ockendon, Ockendon and Johnson \cite{ooj}, Norbury and Yeh
\cite{ny} and Mays and Norbury \cite{mn}. \bigskip

All of these papers use methods which appear to quite different from ours, and
study, for the most part, equations different from $\left(  \ref{1.1}\right)
$. We have not investigated whether these methods, many using infinite
dimensional functional analysis or sophisticated topology, apply to $\left(
\ref{1.1}\right)  $. Most do not mention equations of the general form of
ours, (as in $\left(  \ref{ns1}\right)  $ below), \cite{ampp} being a notable
exception. Most do not study problems on an infinite interval, so that chaos
is not considered. We have seen no other work on the bifurcation problem
considered in section 4.5.

\bigskip The proof of a weak form of chaos which we give seems to us to be
very simple. \ No analysis is required beyond a simple phase plane argument
and the continuity of solutions with respect to initial conditions. \ See
Theorem \ref{thm4} \ and also the first part of section \ref{sec3.4}. \ This
is for small $\varepsilon.$ \ For larger $\varepsilon$ a few estimates are
needed to verify the hypotheses of Theorem \ref{thm4}. \ \ These also appear
in section \ref{sec3.4}. \ \ 

\bigskip

To relate our results to those of traditional nonlinear oscillation theory, we
recall that in \cite{nm} the undamped form of the equation is written as
\begin{equation}
\ddot{u}+\omega_{0}^{2}u=\varepsilon\left(  u^{3}+k\cos\left(  \left(
w_{0}+\varepsilon\sigma\right)  t\right)  \right)  . \label{1.4}%
\end{equation}
This form is chosen to study ``near-resonance'' phenomena, which were the main
interest of much previous work. \ Resonance, or near-resonance, will play no
role in our approach, since we are not studying a small perturbation of a
linear problem.\ \ A rescaling of $\left(  \ref{1.4}\right)  $ to put the
equation in the form $\left(  \ref{1.1}\right)  $ shows that $\varepsilon$
small in $\left(  \ref{1.4}\right)  $ corresponds to $\lambda$ \ large in
$\left(  \ref{1.1}\right)  ,$ though the analogy is not exact because this
rescaling also results in a small amplitude high frequency forcing. \ 

\bigskip

The paper is organized partly according to the restrictions place on
$\varepsilon.$ \ Propositions \ref{thm0}, \ \ref{thm01}, \ref{thm02}, and
Theorem \ref{thm02a}\ in section \ref{sec2}, are for any $\varepsilon>0.$
\ The results in Section 3 are for a range of $\varepsilon$ \ which can be
stated explicitly, while the results in section 4 \ are for ``sufficiently
small'' $\varepsilon.$ \ A more detailed outline of the paper can be found in
section 5. \ 

\section{\label{sec2}\bigskip Results for all positive $\varepsilon$}

At this initial stage it is easy, and we believe of some interest, to study a
more general class of equations. \ So to start with, consider the equation
\begin{equation}
\varepsilon^{2}u^{\prime\prime}=u^{3}-\lambda u+g\left(  t\right)  \label{ns1}%
\end{equation}
where $g\ $is any continuous bounded function on $[0,\infty),$ with initial conditions

\bigskip%
\begin{equation}
u\left(  0\right)  =\alpha,\,\,u^{\prime}\left(  0\right)  =0. \label{1.5}%
\end{equation}
We denote the unique solution by $u_{\alpha}.$ Our initial interest for a
general $g$ is to find solutions which are bounded on $[0,\infty).$ \ Fixing
$\varepsilon>0,$ this can be viewed as a bifurcation problem in $\lambda.$
\ There is a ``main branch'' \ of bounded solutions which exist for all
$\lambda,$ as given in the following result.

\bigskip

\begin{proposition}
\label{thm0} \ For any $\varepsilon>0$ \ and any $\lambda,$ $\left(
\ref{ns1}\right)  $ has a solution satisfying $\left(  \ref{1.5}\right)  $
which is bounded on $[0,\infty).$
\end{proposition}

\begin{proof}
\ Choose $b>0$ \ so large that if $\left|  u\right|  \geq b$ \ then
$u^{3}-\lambda u+g\left(  t\right)  \neq0,$ \ for any $t.$ Then,
\ $u^{\prime\prime}>0$ \ when $u\geq b$ \ and $u^{\prime\prime}<0$ \ when
$u\leq-b.$ \ We define two subsets of the $\alpha$ axis:%
\begin{align*}
A  &  =\left\{  \alpha\,|\,u_{\alpha}\left(  x\right)  >b\text{ for some
}x\geq0\right\} \\
B  &  =\left\{  \alpha\,|\,u_{\alpha}\left(  x\right)  <-b\text{ for some
}x\geq0\right\}  .
\end{align*}
These sets are clearly non-empty (e.g. $b+1\in A$) and open. \ They are
disjoint because $u_{\alpha}$ cannot have a local maximum in the region $u>b$
or a local minimum in the region $u<-b.$ \ Hence there is an $\alpha^{\ast
}\notin A\cup B,$ and this corresponds to a solution which is bounded on
$[0,\infty).$ This proves Proposition \ref{thm0}.
\end{proof}

\bigskip

It should be remarked that for any even function $g\left(  \cdot\right)  ,$
such as cosine, the solutions $u_{\alpha}$ are even and so $u_{\alpha^{\ast}}$
is bounded on the whole real line. \ 

\begin{proposition}
\label{thm01}\ \ For any $\lambda\leq0$ there is only one solution of $\left(
\ref{ns1}\right)  -\left(  \ref{1.5}\right)  $ which is bounded on $[0,\infty).$
\end{proposition}

\begin{proof}
\ \ Suppose that there are two bounded solutions, say $u_{1}$ and $u_{2},$ and
let $v=u_{1}-u_{2}.$ Then
\[
\varepsilon^{2}v^{\prime\prime}=\left(  u_{1}^{2}+u_{1}u_{2}+u_{2}^{2}%
-\lambda\right)  v,\,\,v^{\prime}\left(  0\right)  =0,
\]
and since $u_{1}^{2}+u_{1}u_{2}+u_{2}^{2}$ is positive definite and
$\lambda\leq0,$ \ $v$ \ cannot be bounded on $[0,\infty)$. \ This proves the result.
\end{proof}

\begin{remark}
A similar argument shows that for a general continuous bounded function $g,$
when $\lambda\leq0$ there is exactly one solution of $\left(  \ref{ns1}%
\right)  $ which is bounded on $\left(  -\infty,\infty\right)  .$ However the
existence part of this result is a little longer, requiring a two-parameter
``shooting'' argument and some more complicated topology. \ 
\end{remark}

We now consider what happens as $\lambda$ increases from zero. \ Consider the
function $p_{\lambda}\left(  u\right)  =u^{3}-\lambda u.$ \ As $\lambda$
crosses zero, $p$ develops, by a ``pitchfork bifurcation'', two new roots
(besides $u=0),$ a local maximum at $u=-\sqrt{\frac{\lambda}{3}}$ \ and a
local minimum at $u=\sqrt{\frac{\lambda}{3}}.$ \ As $\lambda$ increases
further, \ the value of $p_{\lambda}$ at these maxima and minima will
eventually exceed the maximum of $\left|  g\left(  t\right)  \right|  .$ Let
\[
\lambda_{0}=\sup\left\{  \lambda\,\,|\,\,\,\,\left|  g\left(  t\right)
\right|  \geq p_{\lambda}\left(  -\sqrt{\frac{\lambda}{3}}\right)  \text{
\ for some }t\right\}  .
\]

\begin{proposition}
\label{thm02}If $\lambda>\lambda_{0},$ \ then for any $\varepsilon>0$ \ there
is a unique bounded solution with $u>0,$ \ $3u^{2}-\lambda>0$ \ , \ and
another bounded solution, also unique, with $u<0,$ \ $3u^{2}-\lambda>0.$
\end{proposition}

\bigskip

\begin{proof}
\ Existence of these solutions is a simple application of the shooting
technique from Proposition \ref{thm0}. \ \ It is easy to show that if
$\lambda>\lambda_{0,}$ then at $u=\sqrt{\frac{\lambda}{3}},$ where
$u^{3}-\lambda u$ \ has a local minimum, $u^{\prime\prime}<0.$ \ We consider
solutions $u_{\alpha}$ with $\sqrt{\frac{\lambda}{3}}<\alpha<b.$ For $\alpha$
in this region and close to $\sqrt{\frac{\lambda}{3}},$ $u_{\alpha}$ decreases
and crosses $\sqrt{\frac{\lambda}{3}}$ transversely. \ For $\alpha$ in this
region and close to $b,$ $u_{\alpha}$ increases and crosses $u=b$, again
transversely. \ It easily follows that there is a solution $u_{\alpha}$ such
that $\sqrt{\frac{\lambda}{3}}<u_{\alpha}\left(  t\right)  <b$ \ for all
$t\geq0.$ \ Similarly, a second solution is found with $-b<u_{\alpha}%
<-\sqrt{\frac{\lambda}{3}}.$ Uniqueness of both of these solutions follows as
in Proposition~\ref{thm01}.
\end{proof}

\begin{remark}
We believe that there should be a third bounded solution as well, but we don't
have a proof for a general continuous bounded function $g.$ \ There appears to
be a significant difference between boundary value problems on finite
intervals, such as those in \cite{ampp}, and the semi-infinite interval
problem discussed here. \ It is easy to show, for a bounded $g$ and any given
$L>0,$ \ that for large enough $\lambda$ there are at least three solutions of
$\left(  \ref{ns1}\right)  -\left(  \ref{1.5}\right)  $ such that $u^{\prime
}\left(  L\right)  =0.$ \ But it does not appear to us to be so easy to find
three bounded solutions on $[0,\infty).$ \ ``Generically'' there should be
three, perhaps by some sort of degree theory argument, but the non-compactness
of $[0,\infty)$ seems to make the use of such arguments more challenging. \ 
\end{remark}

\textbf{\bigskip}

So now we specialize to the equation of principle interest, namely $\left(
\ref{1.1}\right)  .$ It is useful to consider the curve in the $\left(
t,u\right)  $ plane defined implicitly by the equation $u^{\prime\prime}=0;$
that is , by $f\left(  t,u\right)  \equiv u^{3}-\lambda u+\cos t=0.$ \ This
curve has one component for $\lambda\leq\lambda_{0}=\frac{3}{2^{\frac{2}{3}}}$
\ and three components for $\lambda>\lambda_{0}$. \ Here arethe graphs for
$0<\lambda<\lambda_{0,}$ $\lambda=\lambda_{0},$ \ and $\lambda>\lambda_{0}$: \ 

\bigskip%

{\includegraphics[
height=1.9043in,
width=2.0583in
]%
{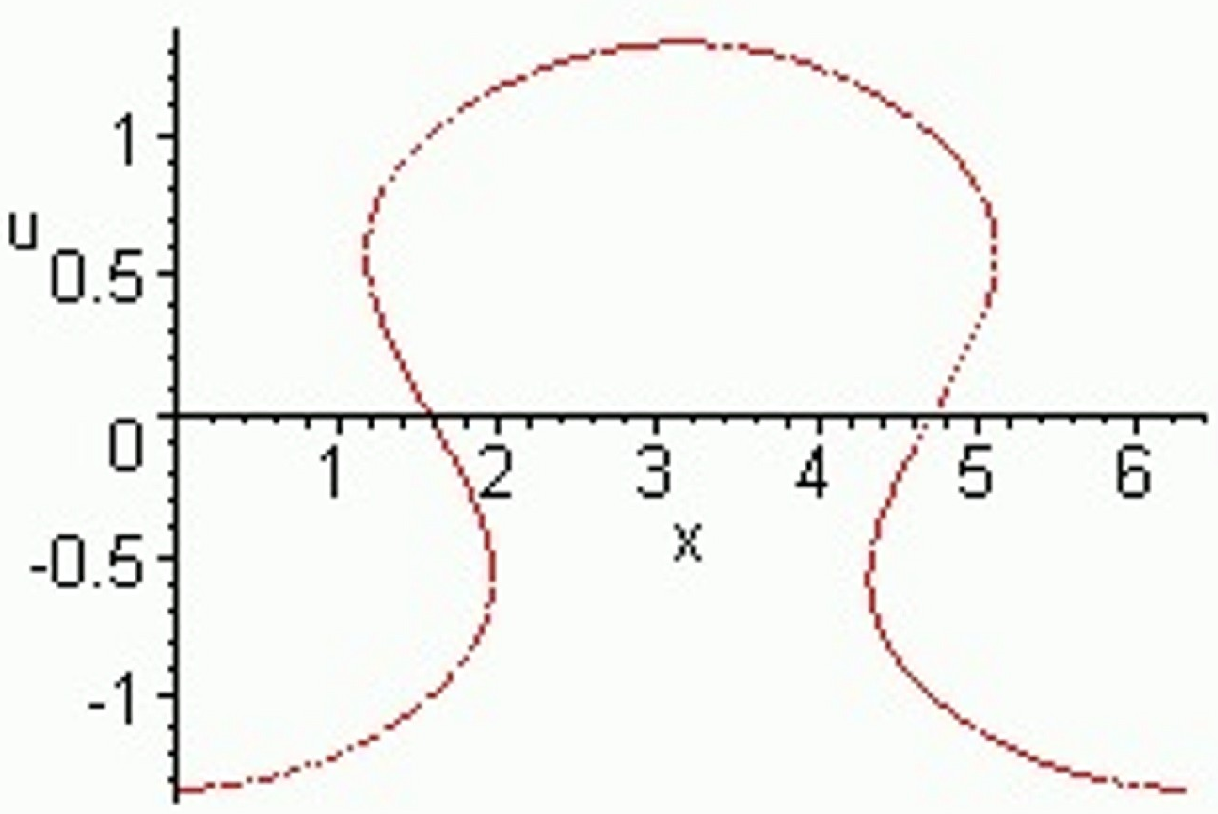}%
}%
{\includegraphics[
height=1.9069in,
width=2.1828in
]%
{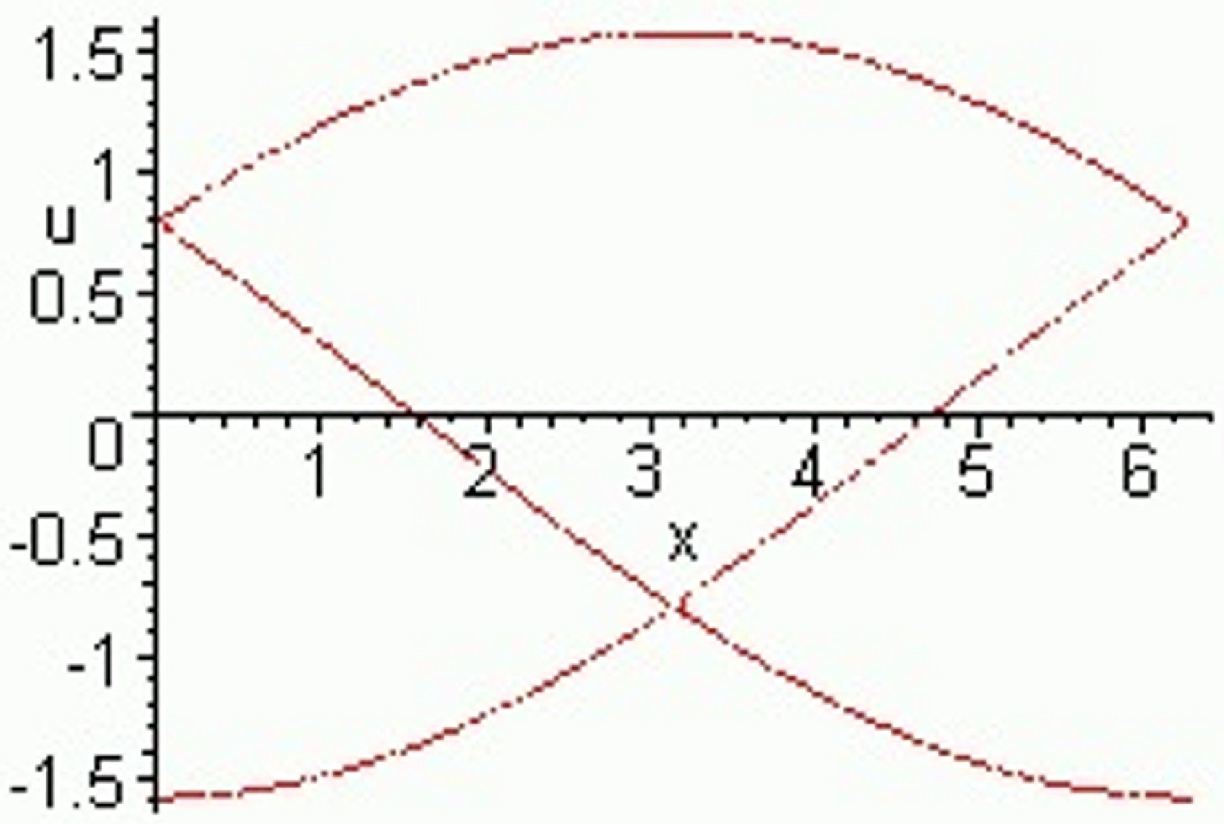}%
}%
{\includegraphics[
height=1.9052in,
width=1.868in
]%
{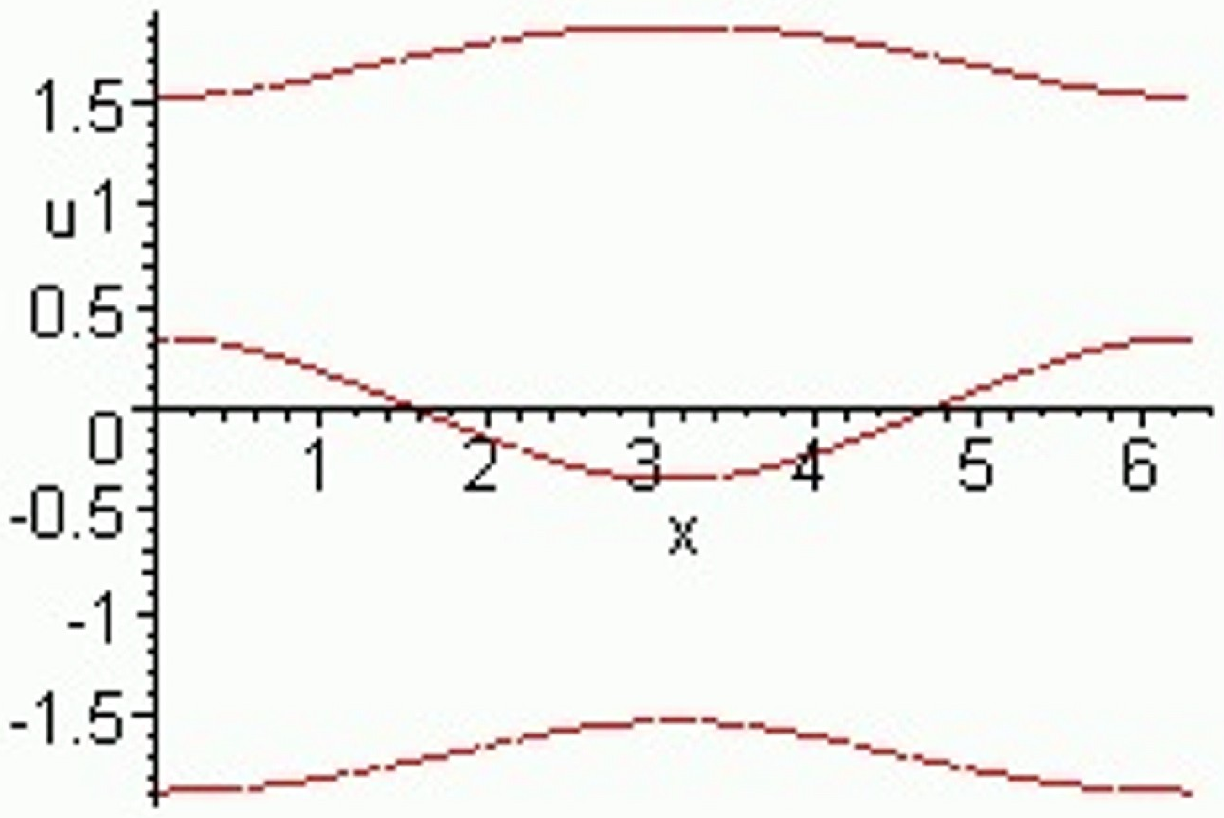}%
}%

\begin{center}
Figure 1.
\end{center}

In all three cases, the curve crosses $u=0$ at odd multiples of $\frac{\pi}%
{2}.$ For $\lambda>\lambda_{0},$ \ let $\bar{U}\left(  t\right)  >U_{0}\left(
t\right)  >$\underline{$U$}$\left(  t\right)  $ be the three solutions of
$u^{3}-\lambda u+\cos\left(  t\right)  =0.$

\bigskip

\begin{theorem}
\label{thm02a} When $g\left(  t\right)  =\cos t,$ let $u_{1}$ and $u_{5}$ be
the solutions found in Proposition \ref{thm02}, with $u_{1}<0$ \ and
$u_{5}>0.$ Then these solutions\ are $2\pi$-periodic, $u_{1}^{\prime}>0$ in
$\left(  0,\pi\right)  $ and $\underline{U}\left(  0\right)  <u_{1}%
<\underline{U}(\pi)$ in $\left[  0,\pi\right]  .$ Also, $u_{5}\left(
t\right)  =-u_{1}\left(  t+\pi\right)  $ and $u_{1}$ and $u_{5}$ are the
minimal and maximal bounded solutions of $\left(  \ref{1.1}\right)  $ on
$\left(  -\infty,\infty\right)  ,$ in the sense that every other bounded
solution lies between these two. \ Also there is a third $2\pi$ periodic
solution, $u_{p}\left(  =u_{\alpha_{p}}\right)  ,$ which satisfies $u\left(
\frac{\pi}{2}\right)  =0$ and $u^{\prime}>0$ on $(0,\frac{\pi}{2}).$ \ A
unique solution with these properties exists for all $\lambda$ (and so for
$\lambda\leq0$ \ this is the solution found in Proposition \ref{thm0}).
\end{theorem}

\begin{proof}
In proving that the solutions of Proposition \ref{thm02} \ are periodic, we
give a second proof of their existence, for the case $g\left(  t\right)  =\cos
t$. Let $\alpha_{0}=\underbar{U}(0)$. First observe that $u_{\alpha_{0}%
}^{\prime\prime}(0)=u_{\alpha_{0}}^{\prime\prime\prime}(0)=0$ while
$\varepsilon^{2}u_{\alpha_{0}}^{(4)}(0)=-1$. Therefore $u_{\alpha_{0}}(t)$
decreases monotonically to $-\infty$. If $\alpha_{0}<\alpha<U_{0}(0)$, then
$u_{\alpha}^{\prime\prime}(0)>0$ and so $u_{\alpha}$ initially increases.
However, for sufficiently small positive values of $\alpha-\alpha_{0}$,
$u_{\alpha}^{\prime}$ has a first zero, which we denote by $t_{1}(\alpha)$.
\ This means, in turn, that $u_{\alpha}^{\prime\prime}$ must have a first zero
at $\tau_{1}\left(  \alpha\right)  \in$ $\left(  0,t_{1}\left(  \alpha\right)
\right)  $. \ For small $\alpha-\alpha_{0},$ $t_{1}\left(  \alpha\right)
<\frac{\pi}{2}.$ \ Further, $\varepsilon^{2}u_{\alpha}^{\prime\prime\prime
}\left(  \tau_{1}\left(  \alpha\right)  \right)  \leq0.$ From
\[
\varepsilon^{2}u^{\prime\prime\prime\prime}\left(  t\right)  =\left(
3u^{2}-\lambda\right)  u^{\prime\prime}+6u\,u^{\prime2}-\cos t
\]
we conclude that for small $\alpha-\alpha_{0}>0,$ $u_{\alpha}^{\prime\prime
}<0$ on $(\tau_{1}\left(  \alpha\right)  ,t_{1}\left(  \alpha\right)  ].$
Hence, $u\left(  t_{1}\left(  \alpha\right)  \right)  <\underline{U}\left(
t_{1}\left(  \alpha\right)  \right)  .$ \ 

\bigskip

Clearly, $t_{1}(\alpha)$ approaches $0$ as $\alpha$ tends to $\alpha_{0}$ from
above. Extend the function $t_{1}(\cdot)$ to larger $\alpha$ by continuity, as
long as possible. That is, $t_{1}(\cdot)$ is the continuous function such that
$u_{\alpha}^{\prime}(t_{1}(\alpha))=0$ and, for $\alpha$ sufficiently close to
$\alpha_{0}$, $t_{1}(\alpha)$ is the first positive zero of $u_{\alpha
}^{\prime}$. Then $t_{1}\left(  \alpha\right)  $ remains the first zero of
$u_{\alpha}$ until either (i) $u_{\alpha}^{\prime\prime}(t_{1}(\alpha))=0$
(since we can then no longer use the implicit function theorem to solve
$u_{\alpha}^{\prime}(t)=0$ for $t$) or (ii) there exists a first $t_{0}%
\in(0,t_{1}\left(  \alpha\right)  )$ such that $u_{\alpha}^{\prime}(t_{0})=0$
and $u_{\alpha}^{\prime\prime}(t_{0})=0$ (since then $t_{1}\left(
\alpha\right)  $ ceases to be the first positive zero of $u_{\alpha}^{\prime}%
$). But if (i) occurs at some first $\alpha_{1}>\alpha_{0}$, and $t_{1}\left(
\alpha_{1}\right)  $ is the first zero of $u_{\alpha_{1}}^{\prime},$ with
$t_{1}\left(  \alpha_{1}\right)  <\pi,$ then $u_{\alpha_{1}}^{\prime
\prime\prime}\left(  t_{1}\left(  \alpha_{1}\right)  \right)  =-\sin
t_{1}\left(  \alpha_{1}\right)  <0.$ This implies that $u_{\alpha_{1}}%
^{\prime}=0$ \ to the left of $t_{1}\left(  \alpha_{1}\right)  $, a
contradiction. If (ii) occurs we get the same contradiction at $t_{0}\left(
\alpha\right)  .$ \ If (i) \ occurs at $\alpha_{1}$ and $t_{1}\left(
\alpha_{1}\right)  =\pi,$ then $u_{\alpha_{1}}^{\prime\prime\prime}\left(
t_{1}\left(  \alpha_{1}\right)  \right)  =0$ \ and $\varepsilon^{2}%
u_{\alpha_{1}}^{\prime\prime\prime\prime}\left(  t_{1}\left(  \alpha
_{1}\right)  \right)  =1.$ This means that $u_{\alpha_{1}}^{\prime\prime}>0$
on either side of $t_{1}\left(  \alpha_{1}\right)  =\pi,$ \ and $t_{1}\left(
\alpha_{1}\right)  $ could not be the first zero of $u_{\alpha_{1}}^{\prime}.$
Hence $t_{1}\left(  \alpha\right)  $ is continuous and remains the first zero
of $u_{\alpha}^{\prime}$ \ as long as $t_{1}\left(  \alpha\right)  \leq\pi.$

\bigskip

On the other hand, if $\alpha=\underline{U}(\pi)$, then the first maximum of
$u_{\alpha}$ is larger than $U_{0}(\pi)$. Therefore, $t_{1}(\alpha)$ must vary
continuously as $\alpha$ increases from $\alpha_{0}$ until it eventually takes
on the value $\pi$ for some $\alpha\in(\alpha_{0},\underline{U}(\pi))$. For
any such $\alpha$, $u_{\alpha}$ gives a periodic solution of (\ref{1.1}) with
the properties that $u_{\alpha}^{\prime}(t)>0$ for $t\in(0,\pi)$ and
$\underline{U}(0)<u_{\alpha}(t)<\underline{U}(\pi)$ for $t\in\lbrack0,\pi]$.
From Proposition~\ref{thm02}, it follows that such an $\alpha$ is unique and
the corresponding solution $u_{\alpha}$ coincides with the unique bounded
solution with $u<0$ and $3u^{2}-\lambda<0$ in Proposition~\ref{thm02}. The
reflection $-u_{\alpha}(\pi+t)$ of this periodic solution gives a second
periodic solution, which corresponds to the unique bounded solution with $u>0$
and $3u^{2}-\lambda>0$ in Proposition~\ref{thm02}. \ 

\bigskip

The proof that $u_{1}$ and $u_{5}$ are the minimal and maximal bounded
solutions of $\left(  \ref{1.1}\right)  $ on $\left(  -\infty,\infty\right)  $
follows from that of Proposition \ref{thm02}. \ \bigskip

To obtain the third periodic solution, $u_{p},$ again let $\alpha
_{0}=\underline{U}\left(  0\right)  ,$ the smallest root of $u^{3}-\lambda
u+1=0.$ We saw that if $\alpha<\alpha_{0}$ then $u^{3}-\lambda u+\cos t<0$
\ for all $u\leq\alpha,$ and all $t,$\ and $u_{\alpha}$\ decreases
monotonically to $-\infty$, (in finite time). If $\alpha=\alpha_{0},$ so that
$u^{\prime\prime}\left(  0\right)  =0,$ then $u^{\prime\prime\prime}\left(
0\right)  =0,\,\,u^{\prime\prime\prime\prime}\left(  0\right)  <0$ \ and again
$u$ \ decreases below $\alpha_{0}$ and tends to $-\infty.$ \ On the other
hand, \ if $\alpha<0$ \ is sufficiently small in magnitude then $\varepsilon
^{2}u^{\prime\prime}\left(  0\right)  >\frac{1}{2}$ \ and it is easy to show
that $u_{{}}$ increases to cross zero at some first $t_{0}=t_{0}\left(
\alpha\right)  ,$ with $u^{\prime}>0$ \ on $(0,t_{0}\left(  \alpha\right)  ].$
\ \ Lowering $\alpha$ \ from $0,$ the function $t_{0}\left(  \alpha\right)  $
is continuous as long as $u_{\alpha}^{\prime}\left(  t_{0}\left(
\alpha\right)  \right)  >0$. \ \ If $u_{\alpha}^{\prime}\left(  t_{0}\left(
\alpha\right)  \right)  =0$ \ then $u^{\prime\prime}\left(  t_{0}\left(
\alpha\right)  \right)  \leq0$ since $u$ cannot have a local minimum at
$t_{0}\left(  \alpha\right)  .$ This implies that $t_{0}\left(  \alpha\right)
\geq\frac{\pi}{2}.$ \ Since $t_{0}\left(  \alpha\right)  $ \ is not defined
for $\alpha\leq\alpha_{0}$, it must increase beyond $\frac{\pi}{2}$ as
$\alpha$ decreases from zero. Thus there is some first (largest)
$\alpha=\alpha_{p}<0$ with $t_{0}\left(  \alpha_{p}\right)  =\frac{\pi}{2}.$
Let $u_{p}=u_{\alpha_{p}}.$ An easy symmetry argument shows that $u_{p}\left(
\frac{\pi}{2}+t\right)  =-u_{p}\left(  \frac{\pi}{2}-t\right)  $ for all $t,$
showing that $u_{p}$ is $2\pi-$periodic, with a maximum at $\pi$. \ 

\bigskip

Suppose, therefore, that $u=u_{p}$ satisfies $u^{\prime}\left(  \tau\right)
=0$ \ for some first $\tau\in\left(  0,\pi\right)  .$ \ By the anti-symmetry
of $u_{p}$ around $\frac{\pi}{2}$ we can assume that $\tau\leq\frac{\pi}{2}.$
Then $u^{\prime\prime}\left(  \tau\right)  \leq0.$ \ If $u^{\prime\prime
}\left(  \tau\right)  =0$ \ then $\varepsilon^{2}u^{\prime\prime\prime}\left(
\tau\right)  =-\sin\tau<0$ so $u^{\prime\prime}$ \ becomes negative. \ In
either case $u^{\prime}<0$ on some interval to the right of $\tau$. \ 

\bigskip

From the graphs of $u^{\prime\prime}=0$ it is clear that if $u\left(
0\right)  <0,\,\,u^{\prime\prime}(0)>0,$ \ $u<0$ on $\left(  0,\frac{\pi}%
{2}\right)  $ and $u^{\prime}=0$ \ before $t=\frac{\pi}{2},$ \ then after that
point, $u^{\prime\prime}<0$ at least until $t=\pi.$ \ This implies that
$u_{p}^{\prime}>0$ \ on $(0,\pi).$

\bigskip

The uniqueness of $u_{p}$ follows from Lemma~\ref{alem5.7} below.

\bigskip

\begin{lemma}
\label{alem5.7} Assume that $\alpha_{1}<\alpha_{2}<0$. If $0<T\leq\pi/2$ and
$u_{\alpha_{2}}\leq0$ in $[0,T]$, then $u_{\alpha_{1}}<u_{\alpha_{2}}$ in
$[0,T]$.
\end{lemma}

\begin{proof}
First, using the Sturm-Liouville technique we get
\begin{equation}
\varepsilon^{2}(u_{\alpha_{1}}^{\prime}u_{\alpha_{2}}-u_{\alpha_{1}}%
u_{\alpha_{2}}^{\prime})(t)=\int_{0}^{t}u_{\alpha_{1}}u_{\alpha_{2}}%
(u_{\alpha_{1}}^{2}-u_{\alpha_{2}}^{2})\,ds + \int_{0}^{t}(u_{\alpha_{2}%
}-u_{\alpha_{1}})\cos t\,ds. \label{a5.20}%
\end{equation}

Assume that the lemma is false. Then, there is a first $t=\bar{t}\in(0,\pi/2]$
such that $u_{\alpha_{2}}\leq0$ in $[0,\bar{t}]$ and $u_{\alpha_{1}}%
=u_{\alpha_{2}}$ at $\bar{t}$. Then, $u_{\alpha_{1}}^{\prime}(\bar{t})\geq
u_{\alpha_{2}}^{\prime}(\bar{t})$. \ Evaluating (\ref{a5.20}) at $t=\bar{t}$
shows that the left side is nonpositive and the right side is positive. This
contradiction proves the lemma and completes the proof of Theorem \ref{thm02a}.
\end{proof}
\end{proof}

\bigskip

Below we will want a slight variant of Lemma \ref{alem5.7}, with a similar
proof which we omit.

\begin{lemma}
\label{alem5.8} Let $u$ and $U$ be two solutions of (\ref{1.1}) with
$u^{\prime}(\pi)=U^{\prime}(\pi)=0$. If $\pi>T\geq\pi/2$ and $u\geq0$ in
$[T,\pi]$, then $u<U$ in $[T,\pi]$.
\end{lemma}

\begin{corollary}
\label{acor5.4} Assume that $u$ is a periodic solution of (\ref{1.1}) with
$u^{\prime}(0)=u^{\prime}(\pi)=0$. If $u(0)<u_{p}(0)$, then $u(t)<u_{p}(t)$
for all $t\in(-\infty,\infty)$.
\end{corollary}

\begin{proof}
From Lemma~\ref{alem5.7} it suffices to show that $u(t)<u_{p}(t)$ for
$t\in(\pi/2,\pi]$. If $u(\pi)>u_{p}(\pi)$, then Lemma~\ref{alem5.8} implies
$u(t)>u_{p}(t)$ for all $t\in\lbrack\frac{\pi}{2},\pi]$ so that $u$ has a jump
at $\frac{\pi}{2}$, a contradiction. Hence we have $u(\pi)<u_{p}(\pi)$. Assume
now that the corollary is false. Since $u_{p}^{\prime}>0$ in $(0,\pi)$, there
is a $\bar{t}\in(\frac{\pi}{2},\pi)$ such that $u^{\prime}(\bar{t})=0$ and
$u(\bar{t})>u_{p}(\bar{t})$. The Sturm-Liouville technique shows that
\[
\varepsilon^{2}[u^{\prime}(t)u_{p}(t)-u(t)u_{p}^{\prime}(t)]=-u(\bar{t}%
)u_{p}^{\prime}(\bar{t})+\int_{\bar{t}}^{t}[uu_{p}(u^{2}-u_{p}^{2}%
)+(u_{p}-u)\cos s]\,ds.
\]
If $u\left(  t\right)  =u_{p}\left(  t\right)  $ for some largest $t\in\left[
\frac{\pi}{2},\bar{t}\right]  $ then $(\frac{u}{u_{p}})^{\prime}<0$ in
$(t,\bar{t})$ a contradiction. We again get that $u$ has a jump at
$t=\frac{\pi}{2}$ and this proves the corollary.
\end{proof}

\bigskip

\bigskip\ 

We have now shown that as $\lambda$ increases from $0$, \ new periodic
solutions appear. \ We can prove one thing about the initial bifurcation of
new solutions for any positive $\varepsilon.$ \ Let
\[
\lambda_{b}=\sup\left\{  \lambda\,|\,\text{There is only one solution, }%
u_{p},\text{with }u^{\prime}\left(  0\right)  =0,\ u^{\prime}\left(
\pi\right)  =0\right\}  .
\]

\bigskip

\begin{theorem}
\label{thm04}For any $\varepsilon>0,$ \ $\lambda_{b}<\lambda_{0}$ . \ 
\end{theorem}

\begin{proof}
\ Suppose that $\varepsilon>0$ \ and let $\lambda=\lambda_{0}.$ \ It will be
helpful to consider the graph of the set of solutions of $u^{\prime\prime}=0,$
\ shown in Figure 1-b. \ Starting from $\alpha=-\sqrt{\frac{\lambda}{3}}$ we
consider the solution $u_{\alpha}$ as $\alpha$ \ is lowered. \ Let
$t_{2}\left(  \alpha\right)  $ denote the first intersection of this solution
with $u=-\sqrt{\frac{\lambda}{3}}.$ \ Then near $\alpha=-\sqrt{\frac{\lambda
}{3}}$, \ $t_{2}\left(  \alpha\right)  $ is defined and continuous. \ Further,
$u_{\alpha}^{\prime}>0$ on $(0,t_{2}\left(  \alpha\right)  ]$ as long as
$t_{2}\left(  \alpha\right)  <\pi.$ The solution cannot be tangent to the line
$u=-\sqrt{\frac{\lambda}{3}}$ \ until we reach a value of $\alpha,$ say
$\alpha_{1},$ \ where $t_{2}\left(  \alpha\right)  =\pi.$ \ Suppose that
$u_{\alpha_{1}}^{\prime}\left(  \pi\right)  =0.$ \ \ Since $\lambda
=\lambda_{0}$ it follows that $u_{\alpha_{1}}^{\prime\prime}\left(
\pi\right)  =0$ and we easily calculate that $u_{\alpha_{1}}^{\prime
\prime\prime}(\pi)=0$ while $u_{\alpha_{1}}^{^{\prime\prime\prime\prime}%
}\left(  \pi\right)  >0.$ \ But this implies that $u_{\alpha_{1}}^{{}}%
=-\sqrt{\frac{\lambda}{3}}$ \ at some earlier point, a contradiction.
\ Therefore $u_{\alpha_{1}}^{\prime}\left(  \pi\right)  >0$. A similar
argument shows that there is an $\alpha_{2}>\sqrt{\frac{\lambda}{3}}$ such
that $u_{\alpha_{2}}(\pi)=\sqrt{\frac{\lambda}{3}}$ and $u_{\alpha_{2}%
}^{\prime}(\pi)<0$. If $\alpha>b$ then $u_{\alpha}^{\prime}(\pi)>0$ (or
$u_{\alpha}(\pi)$ is not defined, which is easily handled), while if
$\alpha<-b$ then $u_{\alpha}^{\prime}(\pi)<0$. Therefore a simple continuity
argument shows that there are three values of $\alpha$ such that $u^{\prime
}(\pi)=0$, and this give three $2\pi$-periodic solutions of \ref{1.1}.

\bigskip

Further, \ this argument extends by continuity to values of $\lambda$ close to
$\lambda_{0}.$ \ That is, keeping $\alpha_{1}$ and $\alpha_{2}$ fixed, we
still have $u_{\alpha_{1}}^{\prime}\left(  \pi\right)  >0,$ $u_{\alpha_{2}%
}^{\prime}\left(  \pi\right)  <0,$ $u_{b}^{\prime}\left(  \pi\right)  >0$ or
undefined, \ and $u_{-b}^{\prime}\left(  \pi\right)  <0$ or undefined. \ This
proves Theorem \ref{thm04}.
\end{proof}

\section{\ \label{sec3}Results for ``large'' $\varepsilon$}

\subsection{\label{sec3.1}A simple proof of chaos}

In this section we show how the shooting method can give a global result about
chaos for this equation. \ Our aim is for a concise statement implying the
existence of uncountably many bounded solutions and a natural map between a
family of symbol sequences and a set of bounded solutions, where moreover, we
can get specific estimates on ranges of $\lambda$ \ and $\varepsilon$ \ which
support this chaotic behavior. \ This will mean that we do not prove
uniqueness of the solution corresponding to a particular sequence, and hence
we can only prove a weaker sensitivity to initial conditions then one obtains
from some asymptotic methods. \ In a later section, in which $\varepsilon$ is
taken as sufficiently small, more precise results will be obtained.

\bigskip

Our results in the next few sections have one main hypothesis, which we state now.

\begin{condition}
\label{A}There is an $\bar{\alpha}\in\left(  \alpha_{0},\text{ }-\sqrt
{\frac{\lambda}{3}}\right)  $ such that $u_{\bar{\alpha}}\left(  t\right)  $
increases monotonically as $t$ increases from $0,$ and crosses $u=b$ \ before
$t=\frac{\pi}{2}.$ \ 
\end{condition}

Recall that $b$ is positive and satisfies $b^{3}-\lambda b-1>0$. An easy phase
plane argument (given in section \ref{sec3.4}) shows that if $\lambda
\geq\lambda_{0}$, \ then Condition $\ref{A}$ is satisfied for sufficiently
small $\varepsilon$. \ Later we will show, with a bit more work, that
Condition $\ref{A}$ is satisfied in a specific range of $\varepsilon$ and
$\lambda$.

\bigskip

We use $\bar{\alpha}$ as given in Condition \ref{A} to define a family
$\left\{  \hat{w}_{k}\right\}  $ of ``special'' solutions as follows:
\begin{equation}
\hat{w}_{k}\left(  t\right)  =\left(  -1\right)  ^{k}u_{\bar{\alpha}}\left(
t-k\pi\right)  . \label{eqnA}%
\end{equation}

Thus, each $\hat{w}_{k}$ is a translation, \ and for odd $k$ a reflection,
\ of the solution $u_{\bar{\alpha}}$. \ It is easily seen that $\left(
-1\right)  ^{k}\hat{w}_{k}$ has its global minimum at $k\pi$ and on each side
of this point increases monotonically to $u=b$ before $\left|  t-k\pi\right|
=\frac{\pi}{2}.$ \ \ 

\bigskip

Let $w_{k}$ denote the restriction of $\hat{w}_{k}$ to the interval $\left[
s_{k},S_{k}\right]  $ in which $\left|  \hat{w}_{k}\right|  \leq b.$ In the
following figure we plot several of the functions $w_{k}.$%

\begin{center}
\includegraphics[
height=1.8862in,
width=2.8945in
]%
{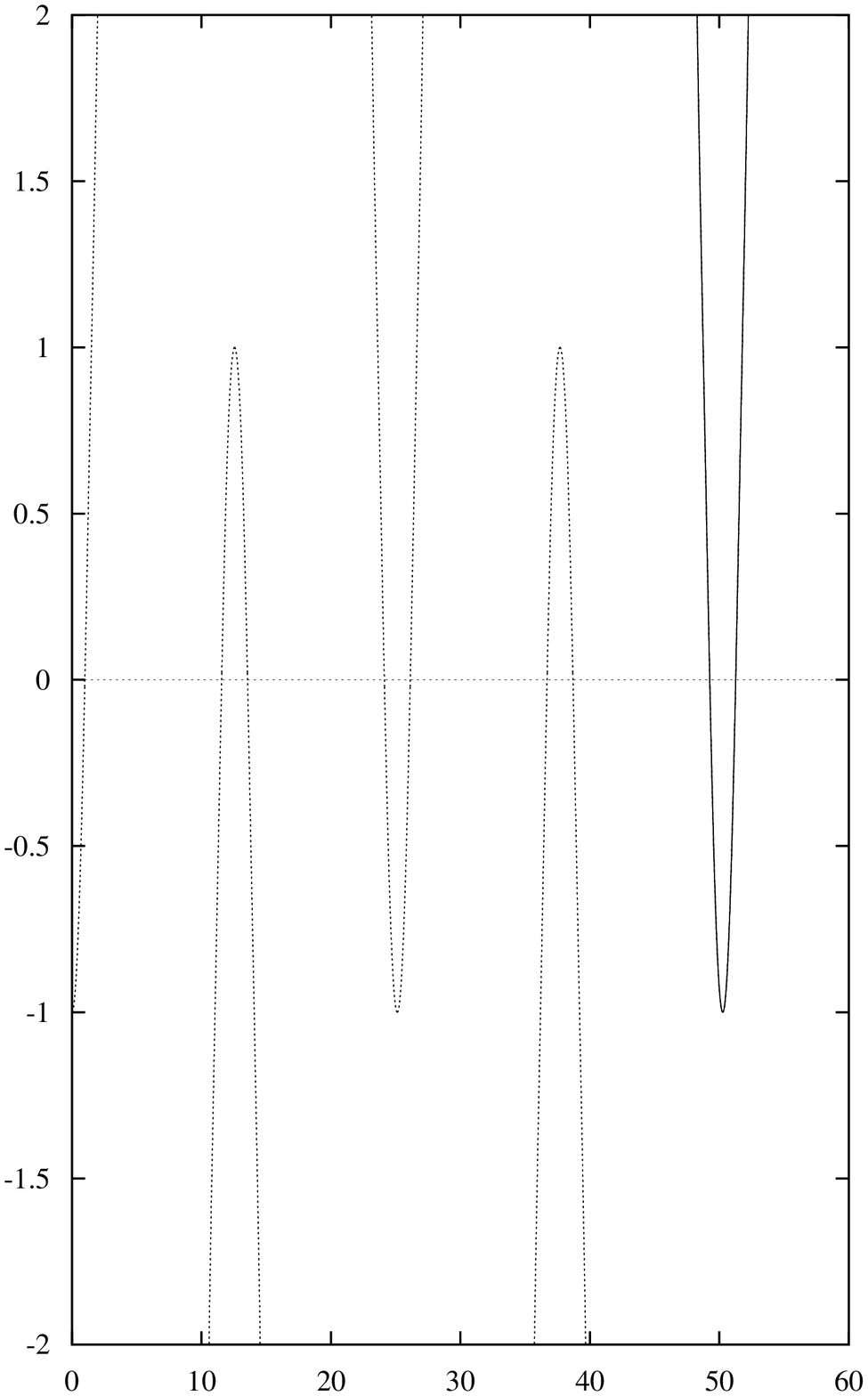}%
\\
Figure 2
\end{center}

\begin{center}
Graphs of $w_{0},w_{1},w_{2},w_{3},$ and $w_{4}$
\end{center}

\bigskip

We can now state our result.

\begin{theorem}
\label{thm4}Assume that $\lambda>\lambda_{0}$ and that Condition \ref{A}
holds. \ Let $\Sigma=\left\{  \sigma_{k}\right\}  _{k=1..\infty}$ be any
sequence of positive integers such that $\ \sigma_{k+1}-\sigma_{k}\geq2$.
\ Then there is an $\alpha_{\Sigma}\in\left(  \alpha_{0},\bar{\alpha}\right)
$ \ \ such that $u_{\alpha_{\Sigma}}$ \ exists on $\left(  -\infty
,\infty\right)  ,$ and on $\left(  0,\infty\right)  $ the graph of
$u_{\alpha_{\Sigma}}$ intersects the graph of each $w_{\sigma_{k}}$ , \ but
does not intersect the graph of any $w_{j}$ with $j\notin\Sigma.$
\end{theorem}

\begin{remark}
\bigskip We allow $\Sigma$ to be finite, or even empty. \ 
\end{remark}

\begin{proof}
We define two functions, \ $f_{+}$ \ and $f_{-},$ \ as follows:%
\begin{align*}
f_{+}\left(  t\right)   &  =\left\{
\begin{array}
[c]{c}%
w_{k}\left(  t\right)  \text{ \ if }k\text{ \ is even and }s_{k}\leq t\leq
S_{k}\\
b\text{ otherwise}%
\end{array}
\right. \\
f_{-}\left(  t\right)   &  =\left\{
\begin{array}
[c]{c}%
w_{k}\left(  t\right)  \text{ if }k\text{ is odd and }s_{k}\leq t\leq S_{k}\\
-b\text{ otherwise}%
\end{array}
\right.
\end{align*}
A key observation is that no solution can intersect the curve $f_{-}$
tangentially from above, \ or the curve $f_{+}$ tangentially from below.
\ This is because $u^{\prime\prime}<0$ when $u=-b,$ $u^{\prime\prime}>0$
\ when $u=b,$ and no two distinct solution graphs can be tangent to each
other. \ Also, no solution can be tangent to either of the lines $u=\pm b$
from between these two lines.

\bigskip

Further, we will need two other functions, $g_{+}$ \ and $g_{-}$, \ where
$g_{-}<0<g_{+}.$ \ These are defined by
\begin{align*}
g_{+}\left(  t\right)  =\max\left\{  \sqrt{\frac{\lambda}{3}},\,\,\,\,f_{-}%
\left(  t\right)  \right\}  , \quad g_{-}\left(  t\right)  =\min\left\{
-\sqrt{\frac{\lambda}{3}},\,\,\ f_{+}\left(  t\right)  \right\}  .
\end{align*}
As in the figure below, the function $f_{+}$ can be described as the function
whose graph is the line $u=b$ \ except at downward bumps when this line meets
a $w_{k}$ with $k$ even. \ Also, \ the graph of $g_{+}$ is the line
$u=\sqrt{\frac{\lambda}{3}}$ except for upward bumps when this line meets the
graph of $w_{k}$ for some odd $k.$ $\ $The functions $\ \ f_{-}$ and $g_{-}$
are reflections and translations of $f_{+}$ and $g_{+}.$ Recall that in
Condition \ref{A} \ we required that $\left|  \bar{\alpha}\right|
>\sqrt{\frac{\lambda}{3}}.$ \ No solution can be tangent to $g_{+}$ from above
or to $g_{-}$ from below. \ \ Also, no solution can be tangent to
$u=\sqrt{\frac{\lambda}{3}}$ from above or to $u=-\sqrt{\frac{\lambda}{3}}$
from below. \ 

\bigskip%
\begin{center}
\includegraphics[
height=1.5in,
width=4.1in
]%
{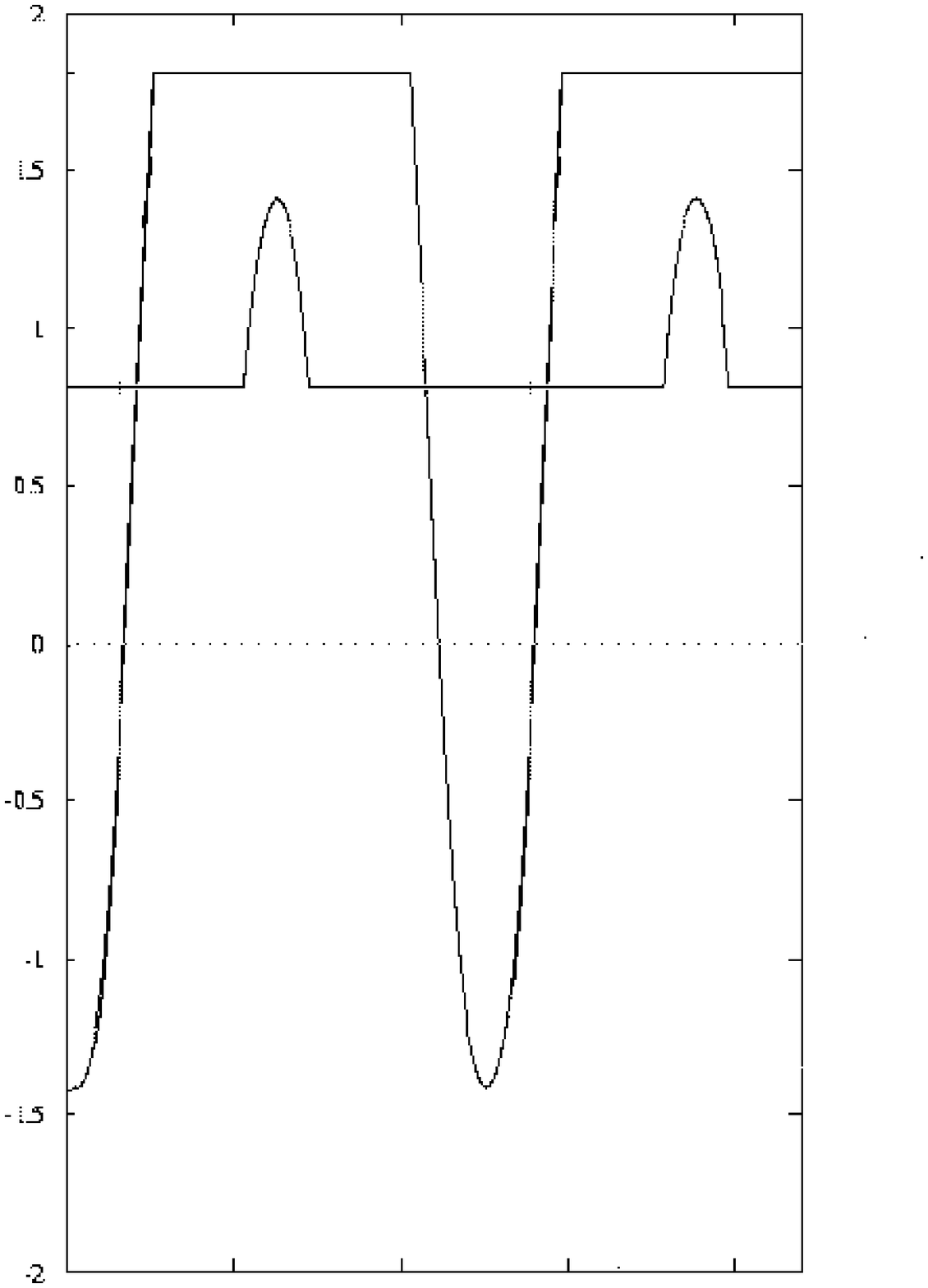}%
\\
Figure 3
\end{center}

Graphs of $f_{+}$ and $g_{+}$

Assuming that the sequence $\Sigma$ is chosen, we shall define a sequence of
closed intervals $\left\{  I_{k}\right\}  $ with the following properties:
(See Figure 4).

\bigskip

$\left(  i\right)  $ \ For each $k,$ $\alpha\in I_{k}$ implies that
$u_{\alpha} $ intersects $w_{\sigma_{1}},...,w_{\sigma_{k}}$ and no other
$w_{j}$ with $j<\sigma_{k}$ . \ 

$\left(  ii\right)  $ Let $I_{k}=\left[  \beta_{k},\alpha_{k}\right]  ,$ and
for $\alpha\in I_{k}$ let $t_{k}\left(  \alpha\right)  $ denote the first $t$
such that $u_{\alpha}\left(  t\right)  =w_{\sigma_{k}}\left(  t\right)  .$
\ Then either $t_{k}\left(  \alpha_{k}\right)  =s_{\sigma_{k}}$ and
$t_{k}\left(  \beta_{k}\right)  =S_{\sigma_{k}},$ \ or $t_{k}\left(  \beta
_{k}\right)  =s_{\sigma_{k}}$ and $t_{k}\left(  \alpha_{k}\right)
=S_{\sigma_{k}}.$

$\left(  iii\right)  $ $I_{k}\subset I_{k-1}.$

\bigskip

Since the intervals will be closed, their intersection will be non-empty, and
a point $\alpha_{\Sigma}\in\cap_{k=1}^{\infty}I_{k}$ will have the properties
stated in the Theorem. \ We do not prove that the intersection is only one
point. \ We will be able to do this for sufficiently small $\varepsilon$.

\bigskip

We will choose $I_{1}$ as a subinterval of $\left(  -b,\bar{\alpha}\right)  .$
\ Suppose, first, that $\sigma_{1}$ is even. \ For $\alpha$ close to
$\bar{\alpha},$ \ $u_{\alpha}$ intersects $f_{+}$ before $t=\frac{\pi}{2}$,
\ and therefore before $u_{\alpha}$ could intersect $w_{1}.$ Let $t_{1}\left(
\alpha\right)  $\ denote the first intersection of $u_{\alpha}$ with $f_{+}.$
Then $t_{1}\left(  \cdot\right)  $ is continuous in some maximal subinterval
$\left(  \beta,\bar{\alpha} \right]  $ of $\left(  -b,\bar{\alpha} \right]  .$
However, $\beta>-b,$ for solutions $u_{\alpha}$ with $\alpha$ \ close to $-b$
\ do not intersect $f_{+}$ at all. \ (Instead they decrease monotonically to
below $u=-b$ and then to $-\infty$.) \ The nontangency of $u_{\alpha}$ with
$f_{+}$ implies that
\[
\lim\sup_{\alpha\rightarrow\beta^{+}}t_{1}\left(  \beta\right)  =\infty.
\]
\ We do not know if $t_{1}\left(  \alpha\right)  $ is monotone increasing.
However, there must be an interval $I_{1}=\left[  \beta_{1},\alpha_{1}\right]
$ such that $t_{1}\left(  \alpha_{1}\right)  =s_{\sigma_{1}},$\thinspace
$t_{1}\left(  \beta_{1}\right)  =S_{\sigma_{1}}$, and $s_{\sigma_{1}}%
<t_{1}\left(  \alpha\right)  <S_{\sigma_{1}}$ in $\left(  \beta_{1},\alpha
_{1}\right)  ,$ $\ t_{1}\left(  \alpha\right)  <S_{\sigma_{1}}$ on $\left(
\beta_{1},\bar{\alpha} \right]  .$ \ As long as $t_{1}\left(  \alpha\right)  $
is continuous, it must remain the first intersection of $u_{\alpha}$ with
$f_{+},$ \ because of the nontangency of $u_{\alpha}$ with $f_{+}.$ \ Further,
while $t_{1}\left(  \alpha\right)  $ is continuous the solution $u_{\alpha}$
cannot intersect $f_{-}$ in $\left(  0,t_{1}\left(  \alpha\right)  \right)  ,$
because of non-tangency with $f_{-}$. \ A similar argument is used to
construct $I_{1}$ if $\sigma_{1}$ is odd. In this case we have $t_{1}%
(\beta_{1})=s_{\sigma_{1}}$ and $t_{1}(\alpha_{1})=S_{\sigma_{1}}$.

\bigskip

\ We now assume that a (decreasing nested) sequence of closed intervals
$I_{1},...,I_{n}$ \ has been constructed with properties $\left(  i\right)
-\left(  iii\right)  $ for $k=1\cdots n.$ \ \ We wish to construct
$I_{n+1}\subset I_{n}.$ \ We will consider the case where $\sigma_{n}$ is even
and $\sigma_{n+1}$ is odd. We will also assume that $t_{n}\left(  \alpha
_{n}\right)  =s_{\sigma_{n}}$ while $t_{n}\left(  \beta_{n}\right)
=S_{\sigma_{n}}$. \ We will show that we can construct $\left[  \alpha
_{n+1},\beta_{n+1}\right]  $ so that $t_{n+1}\left(  \alpha_{n+1}\right)
=S_{\sigma_{n+1}}$ and $t_{n+1}\left(  \beta_{n+1}\right)  =s_{\sigma_{n+1}}$,
or we could make a different choice which would result in $t_{n+1}\left(
\alpha_{n+1}\right)  =s_{\sigma_{n+1}},\,\,t_{n+1}\left(  \beta_{n+1}\right)
=s_{\sigma_{n+1}}.$ \ (Hence, for any given sequence $\Sigma$ \ there will be
many solutions with the properties $\left(  i\right)  -\left(  iii\right)  .$
This does not by itself imply a reduced sensitivity to initial conditions,
because these solutions may be separated from each other. \ This is discussed
in a separate section below. )

\bigskip

\ \ The assumption that $\sigma_{n}$ is even means that the function
$w_{\sigma_{n}}$ is a downward pointing spike and $w_{\sigma_{n}}\left(
s_{\sigma_{n}}\right)  =w_{\sigma_{n}}\left(  S_{\sigma_{n}}\right)  =b,$
while taking $\sigma_{n+1}$ odd \ means that $w_{\sigma_{n+1}}$ is an upward
pointing spike. \ %

\begin{center}
\includegraphics[
height=1.6in,
width=3.5in
]%
{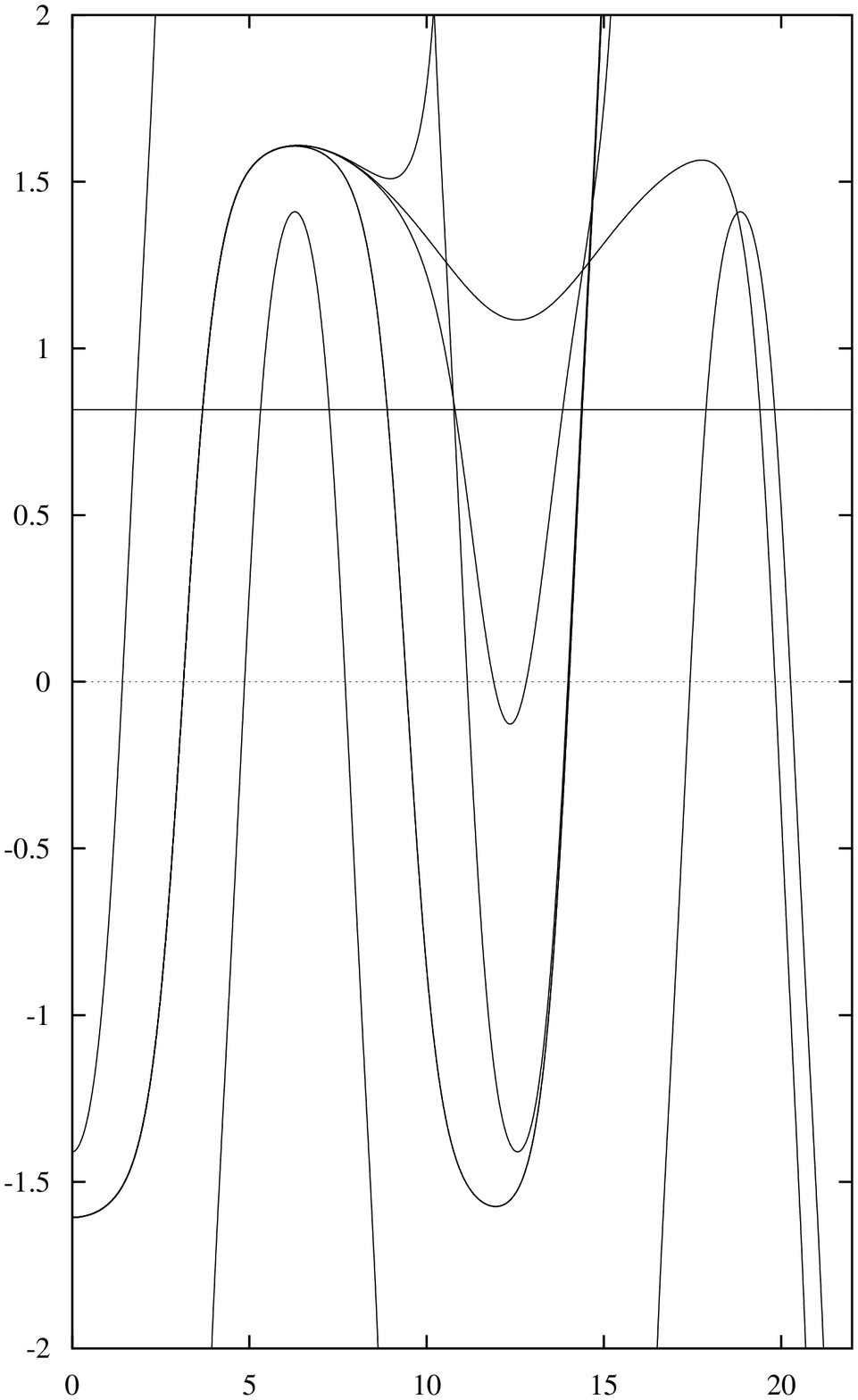}%
\\
Figure 4
\end{center}

Several solutions used in the shooting process, \ with the scale
$t\rightarrow\frac{t}{\varepsilon},$ $\varepsilon=.5$

Solutions shown are $u_{\alpha_{1}},u_{\beta_{1}},$ $u_{\mu_{1}}$and
$u_{\theta_{1}},$ together with $w_{0},w_{1},w_{2},$ and $w_{3}.$

To define $I_{n+1},$ in the case where $\sigma_{n}$ is even, let $\mu_{n}%
=\sup\left\{  \alpha<\alpha_{n}|\,u_{\alpha}\left(  t_{n}\left(
\alpha\right)  \right)  =\sqrt{\frac{\lambda}{3}}\right\}  .$ \ \ Then
$t_{n}\left(  \alpha\right)  \in\left(  s_{\sigma_{n}},\sigma_{n}\pi\right)
,$ $u_{\alpha}\left(  t_{n}\left(  \alpha\right)  \right)  \in\left(
\sqrt{\frac{\lambda}{3}},b\right)  $ for $\mu_{n}<\alpha<\alpha_{n}$ , and
$u_{\mu_{n}}^{\prime}\left(  t_{n}\left(  \mu_{n}\right)  \right)  <0.$ \ The
last inequality is true because $u_{\mu_{n}}$ does not intersect
$w_{\sigma_{n}-1}$, $\ $an upward pointing spike, \ and $\left|  \bar{\alpha
}\right|  >\sqrt{\frac{\lambda}{3}}$. (If we wished to have $t_{n+1}\left(
\alpha_{n+1}\right)  =s_{\sigma_{n+1}}$ we would set $\mu_{n}=\inf\left\{
\beta>\beta_{n}|\,\,u_{\alpha}\left(  t_{n}\left(  \alpha\right)  \right)
=\sqrt{\frac{\lambda}{3}}\right\}  .$ In this case $t_{n}\left(
\alpha\right)  \in\left(  \sigma_{n}\pi,S_{n}\right)  ,$ $\ u_{\alpha}\left(
t_{n}\left(  \alpha\right)  \right)  \in\left(  \sqrt{\frac{\lambda}{3}%
},b\right)  $ for $\beta_{n}<\alpha<\mu_{n}$ , and $u_{\mu_{n}}^{\prime
}\left(  t_{n},\mu_{n}\right)  >0.$ Slight changes are necessary if
$\sigma_{n}$ is odd.) \ \ 

\bigskip

Consider values of $\alpha$ near to $\mu_{n}.$ \ Since $u_{\mu_{n}}\left(
t_{n}\left(  \mu_{n}\right)  \right)  =\sqrt{\frac{\lambda}{3}},$
$\ s_{n}<t_{n}\left(  \mu_{n}\right)  <\sigma_{n}\pi,$\ and $u_{\mu_{n}%
}^{\prime}\left(  t_{n}\left(  \mu_{n}\right)  \right)  <0,$ \ we can define
$\tau_{n}\left(  \alpha\right)  $ continuously by the equations $u_{\alpha
}\left(  \tau_{n}\left(  \alpha\right)  \right)  =g_{+}\left(  \tau_{n}\left(
\alpha\right)  \right)  ,$ $\tau_{n}\left(  \mu_{n}\right)  =t_{n}\left(
\mu_{n}\right)  .$ \ Then $\tau_{n}\left(  \cdot\right)  $ will be continuous
in some maximal interval around $\mu_{n}.$ \ However,\ as $\alpha$ increases
from $\mu_{n},$ $\tau_{n}\left(  \alpha\right)  $ must tend to infinity, since
it is not defined at $\alpha=\alpha_{n},$ where the solution $u$ reaches $b$
at $s_{n}=t_{n}\left(  \alpha_{n}\right)  $ \ and never decreases below $b$
after that.

\bigskip

Let
\[
\theta_{n}=\sup\left\{  \alpha>\mu_{n}|\,\tau_{n}\text{ is continuous on
}\left[  \mu_{n},\alpha\right]  \text{ \ and }\tau_{n}\left(  \alpha\right)
=\left(  \sigma_{n}+1\right)  \pi\right\}  .
\]
The interval $I_{n+1}$ will be a subinterval of $\left[  \theta_{n},\alpha
_{n}\right]  .$ \ \ 

\bigskip

\ Observe that $u_{\theta_{n}}$ crosses $w_{\sigma_{n}+1},$ from above, at its
maximum point, at $t=\left(  \sigma_{n}+1\right)  \pi.$ This is also an
intersection of $u_{\theta_{n}}$ with $f_{-}.$ \ Let the intersection of this
solution with $f_{-}$ \ be denoted by $\rho=\rho_{n}\left(  \theta_{n}\right)
$, \ and extend $\rho_{n}$ as a function of $\alpha$ \ continuously\ for
$\alpha\geq\theta_{n}$, \ $\rho_{n}$ being defined as the intersection of
$u_{\alpha}$ with $f_{-},$ as long as it remains continuous. \ The solution
$u_{\alpha}$ \ may possibly intersect $f_{-}$ at earlier points, but $\rho
_{n}$ \ is defined uniquely by requiring that $\rho_{n}\left(  \theta
_{n}\right)  =\left(  \sigma_{n}+1\right)  \pi.$

\bigskip

\ We see that as $\alpha$ increases from $\theta_{n}$ \ the function $\rho
_{n}$ must eventually be undefined, since it is not defined at $\alpha_{n}$ .
\ So, it must increase (not necessarily monotonically), \ and there must be
some closed subinterval of $\left(  \theta_{n},\alpha_{n}\right)  $ in which
$\rho_{n}\left(  \alpha\right)  $ lies in the interval $\left[  s_{\sigma
_{n+1}},S_{\sigma_{n+1}}\right]  $ and moves from the left end of this
interval to the right end (not necessarily monotonically) \ as $\alpha$
increases. \ This subinterval of $\left(  \theta_{n},\alpha_{n}\right)  $ is
chosen for the interval $I_{n+1}=\left[  \beta_{n+1},\alpha_{n+1}\right]  $.
\ We define it unambiguously by choosing the subinterval with the given
properties which lies nearest to $\theta_{n}.$ \ The possibility of an
infinite set of subintervals with these properties is precluded by bounds on
the variables for a given $\varepsilon$. \ The construction of $I_{n+1}$ shows
that it satisfies conditions (i)-(iii).\ 

\bigskip

A similar construction will give $I_{n+1}$ in the case where $\sigma_{n}$ and
$\sigma_{n+1}$ are both even. \ In this case we can obtain $t_{n+1}\left(
\alpha_{n+1}\right)  =s_{\sigma_{n+1}}$ by decreasing $\alpha$ from
$\alpha_{n}$ and observing that the intersection of $u_{\alpha}$ with
$u=f_{+}$ must tend to $\infty.$ before we reach $\theta_{n}$.\ \ The case
where $\sigma_{n}$ is odd is also handled similarly. At each step we can
obtain $t_{n}(\alpha_{n})=s_{\sigma_{n}}$ if $\sigma_{n}$ is even and
$t_{n}(\alpha_{n})=S_{\sigma_{n}}$ if $\sigma_{n}$ is odd. \ This completes
the induction step and the proof of Theorem \ref{thm4}
\end{proof}

\bigskip

\bigskip

\subsection{\label{sec3.2}``Kneading'' theory}

A one-parameter shooting process as in the proof of Theorem \ref{thm4}
inherently gives an ordering of the initial conditions corresponding to
different sequences. \ \ This sort of ordering is related to ``kneading
theory'' \cite{gh}. \ In earlier work on similar problems formal results about
this were stated, \ for symbol sequences on two symbols \cite{hm1}%
,\cite{hm2},\cite{ht1},\cite{ht2}. \ \ \ Here we will describe the theory only
for the solutions found in Theorem \ref{thm4}. \ \ Further solutions found
below would make the description more complicated and so this part of the
theory will not be discussed in those cases. \ 

\bigskip

Notice that in the induction procedure used in the proof of Theorem
\ref{thm4}, \ repeated use was made of lowering or raising $\alpha$ \ and
letting the intersections with the sets $f_{+}$ \ or $f_{-}$ \ tend to
infinity. \ There was never an assumption that these intersections moved
monotonically. \ However, the intersection will pass through some particular
$w_{k}$ before it first intersects $w_{k+2},$ the next spike pointing in the
same direction as $w_{k}.$ \ The point of intersection might retreat along the
$t$ axis and intersect $w_{k}$ again. \ This would increase the number of
solutions corresponding to some particular sequence, but does not prevent our
giving an order to those chosen as in the proof of theorem \ref{thm4}. \ Our
repeated choice to have $t_{n}\left(  \alpha_{n}\right)  =s_{\sigma_{n}},$
\ rather than $S_{\sigma_{n}},$ \ makes the relation between solutions and
sequences easier to describe. \ (If $\sigma_{1}=1$ \ then we could not make
this choice at step 1, but could thereafter.) \ The ordering of solutions
corresponding to sequences by the algorithm used in the proof is well defined
because of the definitions of $\mu_{n},$ $\theta_{n},$ \ etc. \ \ This is
despite the fact that we don't know that the nested intervals converge to
single points. \ Any choice of points within the limit intervals will give the
same ordering. \ \ We have the following corollary of the proof of Theorem
\ref{thm4}:

\bigskip

\begin{theorem}
\label{thmknead} Suppose that $\Sigma_{1}$ and $\Sigma_{2}$ are two sequences
of positive integers as in the statement of Theorem \ref{thm4}. Suppose that
for some $k\geq0$, $\sigma_{i}(\Sigma_{1})=\sigma_{i}(\Sigma_{2})$ if $1\leq
i\leq k$. If $\sigma_{k+1}(\Sigma_{1})$ is even and $\sigma_{k+1}(\Sigma_{2})$
is odd, then $\alpha(\Sigma_{1})>\alpha(\Sigma_{2})$. If $\sigma_{k+1}%
(\Sigma_{1})>\sigma_{k+1}(\Sigma_{2})$ and both of these integers are even,
then $\alpha(\Sigma_{1})<\alpha(\Sigma_{2})$. If $\sigma_{k+1}(\Sigma
_{1})>\sigma_{k+1}(\Sigma_{2})$ and both integers are odd, then $\alpha
(\Sigma_{1})>\alpha(\Sigma_{2})$.
\end{theorem}

\bigskip

\subsection{\label{sec3.3}Further periodic solutions, and symbolic dynamics
with five symbols.}

\ Earlier we showed that for $\lambda>\lambda_{0},$ \ there are at least three
solutions with period $2\pi.$ \ Numerically, however, we find five periodic
solutions $u_{1},...,u_{5}$ , with $u_{1}<u_{2}<u_{3}<u_{4}<u_{5.}$ \ Here is
a graph of these solutions when $\varepsilon=1,\,\,\lambda=2$:%

\begin{center}
\includegraphics[
height=1.6224in,
width=3.576in
]%
{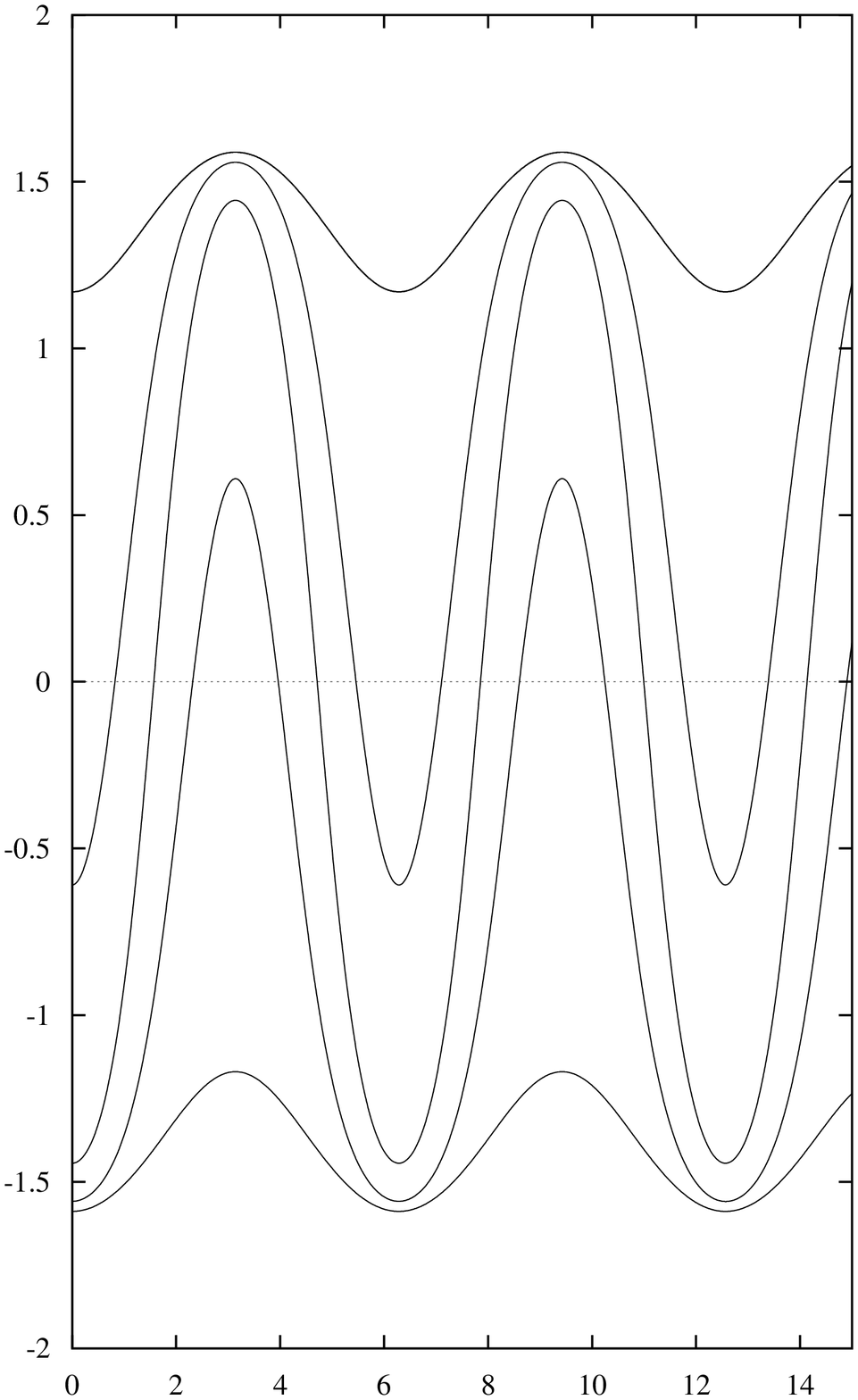}%
\\
Figure 5
\end{center}

\begin{center}
Periodic solutions, \ $u_{1},u_{2},u_{3}\left(  =u_{p}\right)  ,u_{4,}$and
$u_{5}.$ Here $\varepsilon=1,$ $\lambda=2.$
\end{center}

\bigskip

The solution $u_{3}$ is the solution $u_{p}$ which is anti-symmetric around
$\frac{\pi}{2}$, found in Theorem \ref{thm02a}. \ The solutions $u_{1}$ and
$u_{5}$ are the ones found in Proposition \ref{thm02}. \ \ These are the same
ones found in \cite{ampp}, \ at least for small $\varepsilon,$ and correspond
to stable steady states for $\left(  \ref{1.2}\right)  -\left(  \ref{1.3}%
\right)  .$\ \ \ The solutions $u_{2}$ and $u_{4}$ seem to be new (although
their existence is strongly suggested by the existence of three stable
steady-states). \ Other solutions are found which ``follow'' one or another of
these solutions over long intervals, switching in almost arbitrary ways from
one to another. To state a precise theorem, we have to capture what it means
to ``follow'' a solution over an interval in a way that sharply distinguishes
different patterns. \ We need to identify disjoint open sets of solutions
which, over a finite interval, follow the solutions in different orders. \ We
are able to do this because of the various non-tangency conditions stated above.

\bigskip

The following result is stated without reference to $u_{1},...,u_{5,}$ or even
knowledge that $u_{2}$ and $u_{4}$ exist. \ It will be seen later that
existence of $u_{2}$ \ and $u_{4}$ is an easy consequence of the technique
used to prove the result.

\bigskip

\begin{theorem}
\label{thm5} Suppose that the hypotheses of Theorem \ref{thm4} $\ $are
satisfied. \ (That is, $\lambda>\lambda_{0}$ and Condition \ref{A} \ holds.)
\ Let $\left\{  w_{k}\right\}  $ be the family of ``spike'' solutions defined
by equation $\left(  \ref{eqnA}\right)  $. \ Let $\Omega=\left\{  \omega
_{k}\right\}  $ \ be a sequence chosen from the integers $\left\{
1,2,3,4,5\right\}  .$ \ Assume that no $1$ or $2$ is followed immediately by a
$4$ or $5$, \ and no $4$ or $5$ is followed immediately by a $1$ or $2.$
(There must be an intervening $3.)$ \ Then there is a solution $u=u_{\Omega}$
with the following properties: \
\begin{align*}
\left(  i\right)  \text{ \ If }\omega_{k}  &  =1\text{ or }2\text{ \ then
}u\text{ intersects }w_{2k-1}\text{ at two points and does not}\\
&  \text{ intersect }w_{2k}.\\
\left(  ii\right)  \text{ \ If }\omega_{k}  &  =1\text{ \ then }u<g_{-}\text{
in }\left(  \left(  2k-2\right)  \pi,2k\pi\right)  \text{ \ while if }%
\omega_{k}=2,\text{ then }u>g_{-}\text{ somewhere }\\
&  \text{in this interval.}\\
\left(  iii\right)  \text{ \ If }\omega_{k}  &  =3\text{ \ then }u\text{ does
not intersect either }w_{2k-1}\text{ or }w_{2k}.\\
\left(  iv\right)  \text{ If }\omega_{k}  &  =4\text{ or }5\text{ , \ then
}u\text{ does not intersect }w_{2k-1}\text{ \ but does intersect }%
w_{2k},\text{ at two points}.\\
\left(  v\right)  \text{ If }\omega_{k}  &  =4\text{ \ then }u\text{ falls
below }g_{+}\text{ in }\left(  \left(  2k-1\right)  \pi,\left(  2k+1\right)
\pi\right)  ,\text{ while if }\omega_{k}=5,\text{ \ then }\\
u  &  >g_{+}\text{ in this interval.}%
\end{align*}
\end{theorem}

\begin{proof}
\ The proof is a refinement of the proof of Theorem \ref{thm4}. \ The notation
is different, however, for in the statement of Theorem \ref{thm5} all of the
$w_{k}$ are included, not just those picked out by some index set $\left\{
\sigma_{k}\right\}  .$ As before, $s_{k}<S_{k}$ will be the endpoints of the
graphs of $w_{k}$. \ The functions $f_{-}$ and $f_{+},$ $g_{-}$ and $g_{+},$
are the same as before. Also, $t_{k}\left(  \alpha\right)  $ will denote the
first intersection of $u_{\alpha}$ with $w_{k}.$

\bigskip

Rather than carry out a formal induction, we will show how to construct a
solution corresponding to a sequence beginning with $\omega_{1}=4$ or
$5,\omega_{2}=3,\omega_{3}=1$ or $2.$ Hence we want a solution which misses
$w_{1},$ intersects $w_{2},$ misses $w_{3}$ and $w_{4},$ intersects $w_{5}$,
\ and does not intersect $w_{6}.$ It will then be apparent how to deal with
the general case. \ In the notation of Theorem \ref{thm4} we are starting with
$\sigma_{1}=2,\sigma_{2}=5.$

\bigskip

If $\omega_{1}=5,$ then we proceed in the same way as in the proof of Theorem
\ref{thm4}. \ The induction step there (going from $n=1$ to $n=2)$ produced an
interval $I_{2}$ in which the solutions $u_{\alpha}$ do not cross $g_{+}$ in
$\left[  \pi,3\pi\right]  .$ \ Further, the solutions in this interval will
not intersect $w_{3}$ or $w_{4}$ and will intersect the two endpoints of
$w_{5}$ when $\alpha$ is at the endpoints of $I_{2}$. \ The refinement of this
interval according to whether $\omega_{3}=1$ or $2$ \ will be similar to the
way we handle the case $w_{1}=4,$ so we turn to that case now.

\bigskip

When $\omega_{1}=4,$ we follow the proof of Theorem \ref{thm4} up to the
definition of $\mu_{1}.$ However now we set
\[
\mu_{1}=\inf\left\{  \alpha>\beta_{1}|\,\ t_{1}\left(  \alpha\right)
<2\pi\text{ and }u_{\alpha}\left(  t_{1}\left(  \alpha\right)  \right)
=\sqrt{\frac{\lambda}{3}}\right\}  .
\]
Then for $\alpha\in\left[  \beta_{1},\mu_{1}\right]  ,$ $u_{\alpha}$
intersects $w_{2}$ to the right of the first point where $w_{2}=\sqrt
{\frac{\lambda}{3}}.$ We now let $\tau_{1}\left(  \alpha\right)  $ be the
intersection of $u_{\alpha}$ with $g_{+}$ defined by continuity from the point
$u_{\mu_{1}}\left(  t_{1}\left(  \mu_{1}\right)  \right)  .$ \ Set
\[
\theta_{1}=\inf\left\{  \alpha>\mu_{1}|\,\,\tau_{1}\left(  \alpha\right)
=3\pi\right\}  ,
\]
and%
\[
\phi_{1}=\sup\left\{  \alpha<\mu_{1}|\,\,t_{1}\left(  \alpha\right)
=2\pi\right\}  .
\]
These are defined because $u_{\alpha_{1}}$ does not intersect $g_{+}$ after
$\pi,$ and because $t_{1}\left(  \alpha\right)  $ varies between $\tau
_{1}\left(  \mu_{1}\right)  $ and $S_{1}$ as $\alpha$ moves from $\mu_{1}$ to
$\beta_{1}.$ Then $\phi_{1}\in\left(  \beta_{1},\theta_{1}\right)  $ and for
$\alpha\in\lbrack\phi_{1},\theta_{1})$ \ $u_{\alpha}$ intersects $w_{2}$ and
also intersects $g_{+}$ in the interval $\left(  \pi,3\pi\right)  .$

\bigskip

The solution $u_{\phi_{1}}$ intersects $w_{2}$ at its minimum point, at
$t=2\pi,$ and from there increases to cross $b$ before $\frac{5\pi}{2}.$ (This
is shown by a comparison with $w_{2},$ as in the proof of Lemma \ref{alem5.7}%
.) \ The solution $u_{\theta_{1}}$ intersects $w_{2}$ twice, intersects
$w_{3}$ at its maximum point \ and then decreases to cross $-b$ \ before
$\frac{7\pi}{2}.$ \ \ Following the same procedure as in the proof of Theorem
\ref{thm4}, \ we lower $\alpha$ from $\theta_{1}$ \ and find an interval
$I_{2}=\left[  \beta_{2},\alpha_{2}\right]  \subset\left(  \phi_{1},\theta
_{1}\right)  $ such that $u_{\alpha_{2}}\left(  s_{3}\right)  =-b,$
\ $u_{\beta_{2}}\left(  S_{3}\right)  =-b,$ \ and for $\alpha\in\left(
\beta_{2},\alpha_{2}\right)  ,$ $u_{\alpha}$ intersects $g_{+}$ in $\left(
\pi,3\pi\right)  $ and also intersects $w_{3}.$ \ Other cases, and succeeding
steps, being similar, this completes the proof of Theorem \ref{thm5}.
\end{proof}

\bigskip

We can also use this process to prove the existence of $u_{2}$ and $u_{4}$:

\begin{theorem}
\ If $\lambda>\lambda_{0}$ and Condition \ref{A} holds, then there are at
least five $2\pi$-\ periodic solutions of \ref{1.1}, \ $u_{1},u_{2}%
,u_{3}\left(  =u_{p}\right)  ,u_{4}$ and $u_{5},$ \ all with $u^{\prime
}\left(  \pi\right)  =0$. Further, these solutions are ordered by $u_{1}%
<u_{2}<u_{3}<u_{4}<u_{5}.$
\end{theorem}

\begin{proof}
We already know that $u_{1},u_{3},$ \ and $u_{5},$ which exist by Proposition
\ref{thm02}, are increasing in $\left(  0,\pi\right)  .$ We give two proofs of
the existence of $u_{2}$ and $u_{4}$. The first is a little shorter, while the
second shows in addition that $u_{2}$ and $u_{4}$ are increasing on $\left(
0,\pi\right)  .$

\bigskip

Proof 1: \ We will find a solution $u=u_{\alpha}$ with the property that
$\alpha<\alpha_{p},$ $u^{\prime}\left(  \pi\right)  =0$, \ $u\left(
\pi\right)  >g_{-},$ and $u$ \ intersects $w_{1}.$ (None of the three
solutions $u_{1},u_{3},$ \ and $u_{5}$ have all of these properties.) \ We
assume in Theorem \ref{thm5} \ that $\omega_{1}=2.$ \ Then in the construction
used in the proof of Theorem \ref{thm5} \ we find $\mu_{1}$ such that
$u^{\prime}>0$ on $(0,t_{1}\left(  \mu_{1}\right)  ]$ \ and $u_{\mu_{1}%
}\left(  t_{1}\left(  \mu_{1}\right)  \right)  =g_{-}.$ We also find $\phi
_{1}$ such that $t_{1}\left(  \phi_{1}\right)  =\pi,$ and necessarily,
$u_{\phi_{1}}^{\prime}\left(  t_{1}\left(  \phi_{1}\right)  \right)  <0.$ \ If
we decrease $\alpha$ from $\phi_{1}$ we come to a first (largest) $\alpha$,
say $\hat{\alpha},$ where $u_{\hat{\alpha}}$ first intersects $g_{-}$ \ at
$t=\pi.$ At that point, $u_{\hat{\alpha}}^{\prime}\left(  \pi\right)  >0.$
\ \ Somewhere between $\hat{\alpha}$ and $\phi_{1}$ there must be an $\alpha$
with $u_{\alpha}^{\prime}\left(  \pi\right)  =0,$ $u_{\alpha}\left(
\pi\right)  >g_{-}.$ \ This solution must also intersect $w_{1},$ and this
gives the desired fourth periodic solution (beyond those constructed in
Proposition \ref{thm02}.) \ The fifth is obtained by translation and
reflection. \ From Corollary~\ref{acor5.4}, it follows that $u_{1}<u_{2}%
<u_{3}<u_{4}<u_{5}.$

\bigskip

Here is a second proof: \bigskip If $\alpha=u_{1}(0)$, then $v(t)=\frac
{\partial u_{\alpha}}{\partial\alpha}|(t)>0$ and $v^{\prime}(t)>0$ for
$t\in(0,\pi]$, and $u_{1}^{\prime}>0$ in $(0,\pi)$. $\ $Also, $u_{1}$
intersects $w_{1}$ in $(0,\pi)$, and it follows that for $\alpha>u_{1}\left(
0\right)  $ and sufficiently close to $u_{1}\left(  0\right)  ,$ $u_{\alpha
}^{\prime}$ \ is also positive on $\left(  0,\pi\right)  ,$ and again
$u_{\alpha}$ intersects $w_{1}$ at some point in $\left(  0,\pi\right)  .$
\ Let
\begin{align*}
\alpha_{2}=\sup\{\hat{\alpha}  &  \in(u_{1}(0),u_{p}(0))\,|\,\,\text{For
}u_{1}\left(  0\right)  <\alpha\leq\hat{\alpha},\text{ \ }u_{\alpha}^{\prime
}>0\text{ on }\left(  0,\pi\right)  \text{ }\\
&  \text{and }u_{\alpha}\text{ intersects }w_{1}\text{ at some point in
}\left(  0,\pi\right)  \}.
\end{align*}
Then $\alpha_{2}<u_{p}(0)$ since $u_{p}$ does not intersect $w_{1}$ at all.

\bigskip

By continuity, and since no $u_{\alpha}$ is tangent to $w_{1},$ $\ u_{2}%
:=u_{\alpha_{2}}$ intersects $w_{1}$ in $(0,\pi]$ and $u_{2}^{\prime}\geq0$ in
$[0,\pi]$. The first intersection of $u_{2}$ with $w_{1}$ is either at $\pi$
or in $(0,\pi)$. In the first case $u_{2}^{\prime}(\pi)>0,$ which implies
$(u_{2}-w_{1})^{\prime}(\pi)>0$ and so $u_{2}<w_{1}$ just before $\pi$, a
contradiction. Hence $u_{2}$ first intersects $w_{1}$ in $(0,\pi)$. Then by
the definition of $\alpha_{2}$, it follows that $u_{2}^{\prime}=0$ at some
first $\tilde{t}\in(0,\pi]$. We claim that $\tilde{t}=\pi$. If $0<\tilde
{t}<\pi$. then $u_{2}^{\prime\prime}(\tilde{t})=0$ for otherwise,
$u_{2}^{\prime\prime}(\tilde{t})<0$ and this implies that $u_{2}^{\prime}<0$
just after $\tilde{t}$, a contradiction. So $u_{2}^{\prime\prime\prime}%
(\tilde{t})=-\sin\tilde{t}<0$ which gives $u_{2}^{\prime}<0$ just before and
after $\tilde{t}$, again a contradiction. This shows that $u_{2}^{\prime}%
(\pi)=0$ so that $u_{2}$ is a periodic solution which satisfies $u_{1}%
(0)<u_{2}(0)<u_{p}(0)$ and $u_{2}^{\prime}>0$ in $(0,\pi).$ Further, $u_{2}$
intersects $w_{1}$ in $(0,\pi)$, which implies that $u_{2}(\pi)<w_{1}(\pi)$.
Let $u_{4}(t)=-u_{2}(t+\pi)$. From Theorem~\ref{thm02a} and
Corollary~\ref{acor5.4}, it again follows that $u_{1}<u_{2}<u_{3}<u_{4}<u_{5}.$
\end{proof}

\bigskip

\begin{remark}
If the forcing term had less symmetry than cosine, \ but still was symmetric
around $\pi$ \ and $0,$ we could still construct the five periodic solutions
by a similar shooting method. \ 
\end{remark}

\subsection{\bigskip\label{sec3.4}Verification of Condition \ref{A} for a
range of parameters}

\bigskip

\ It is trivial to verify Condition \ref{A} for $\lambda>\lambda_{0}$ and
sufficiently small $\varepsilon.$ \ This is most easily seen by making the
change of variables $\tau=\frac{t}{\varepsilon}$ , letting $v\left(
\tau\right)  =u\left(  t\right)  ,$ \ and considering the resulting equation
\begin{equation}
\ddot{v}=v^{3}-\lambda v+\cos\varepsilon\tau. \label{slow1}%
\end{equation}
Then consider the phase plane for the limiting equation
\[
\ddot{v}=v^{3}-\lambda v+1.
\]
There are three equilibria in the $v,\dot{v}$ plane, $\mathbf{p}_{i}=\left(
\bar{\alpha}_{i},0\right)  ,$ $i=1,2,3,$ where $\bar{\alpha}_{1}<0<\bar
{\alpha}_{2}<\bar{\alpha}_{3}.$ \ The outer two, $\mathbf{p}_{1}$ \ and
$\mathbf{p}_{3}$ \ are saddle points while $\mathbf{p}_{2}$ \ is a center.
\ \ There is a homoclinic orbit at $\mathbf{p}_{3},$ \ so that the left branch
of the unstable manifold at $\mathbf{p}_{3}$ is bounded, \ with $v$ taking its
minimum value at a point $\left(  \alpha^{\ast},0\right)  $ where $\bar
{\alpha}_{1}<\alpha^{\ast}<\bar{\alpha}_{2}$, while the unstable manifold at
$\mathbf{p}_{1}$ is unbounded to the right. \ Suppose that $v\left(  0\right)
=\bar{\alpha}\in\left(  \bar{\alpha}_{1},\alpha^{\ast}\right)  .$ \ Then
$\dot{v}>0$ \ for $t>0$ \ as long as the solution is defined, and
$v\rightarrow\infty$ \ in finite time. \ \ Continuity with respect to
$\varepsilon$ \ implies that for small enough $\varepsilon,$ the solution of
$\left(  \ref{slow1}\right)  $ with $v\left(  0\right)  =\bar{\alpha},$
$\dot{v}\left(  0\right)  =0$ will increase monotonically and cross $v=b$
\ before $\tau=\frac{\pi}{2\varepsilon}.$ \ \ This verifies Condition \ref{A}
for small $\varepsilon.$ \ 

\bigskip

Next we give an estimate on $\varepsilon$ for given $\lambda>\lambda_{0}$ to
ensure that Condition~\ref{A} holds. This estimate of $\varepsilon$ is by no
means best.

\begin{lemma}
\label{alem6.1} For any $\lambda\geq\lambda_{0}$ let $b=b(\lambda
)=\sqrt{\lambda+\frac{1}{2\lambda}}$ and $\varepsilon_{\lambda}=\frac{\pi
}{3T_{\lambda}},$ where $T_{\lambda}=2\sqrt{b+\sqrt{\lambda}}$ for
$\lambda_{0}\leq\lambda<4$ and
\begin{equation}
T_{\lambda}=2\sqrt{2}+\frac{\sqrt{2}\,\ln(\sqrt{\lambda}-1)}{\sqrt{\lambda}%
}+2\Big(\sqrt{b+\sqrt{\lambda}}-\sqrt{2\sqrt{\lambda}-2}\,\Big) \label{a6.3}%
\end{equation}
for $\lambda\geq4.$ Let $v$ be the solution of (\ref{slow1}) with
$v(0)=-\sqrt{\lambda}$ and $\dot{v}(0)=0$. If $0<\varepsilon\leq
\varepsilon_{\lambda}$ in (\ref{slow1}), then there is a $T\in(0,T_{\lambda})$
such that
\begin{equation}
\dot{v}>0\quad\mbox{ in }(0,T]\subset(0,\frac{\pi}{3\varepsilon}],\quad
\mbox{ and}\quad v(T)=b. \label{a6.4}%
\end{equation}
\end{lemma}

\begin{proof}
\ Let
\[
T=\sup\{\tau\in(0,\frac{\pi}{3\varepsilon})|\,\,\dot{v}>0,v<b\quad\mbox
{ in }(0,\tau]\}.
\]
Since $\ddot{v}(0)=1$, it follows that $T$ is well defined. Then for $\tau
\in(0,T)$, $\cos(\varepsilon\tau)\geq1/2$ and then (\ref{slow1}) gives
$\ddot{v}\geq v^{3}-\lambda v+1/2$. Multiply this inequality by $2\dot{v}$ and
integrate over $[0,\tau]$ to give
\begin{equation}
(\dot{v})^{2}\geq\frac{1}{2}v^{4}-\lambda v^{2}+v-(\frac{1}{2}\lambda
^{2}-\lambda^{2}-\sqrt{\lambda})=\frac{1}{2}(v^{2}-\lambda)^{2}+(v+\sqrt
{\lambda})
\end{equation}
and hence, for $\tau\in(0,T]$,
\begin{equation}
\dot{v}\geq\sqrt{\frac{1}{2}(v^{2}-\lambda)^{2}+(v+\sqrt{\lambda})}.
\end{equation}
This implies that
\begin{equation}
T<\int_{-\sqrt{\lambda}}^{b}\frac{dv}{\sqrt{\frac{1}{2}(v^{2}-\lambda
)^{2}+(v+\sqrt{\lambda})}}%
\end{equation}
and $\dot{v}(T)>0$ (recall that $v(0)=-\sqrt{\lambda}$).

\bigskip

If $\lambda_{0}\leq\lambda<4$, then $T<\int_{-\sqrt{\lambda}}^{b}\frac
{1}{\sqrt{v+\sqrt{\lambda}}}dv=2\sqrt{b+\sqrt{\lambda}}=T_{\lambda}$, which
proves (\ref{a6.4}) for $0<\varepsilon<\varepsilon_{\lambda}$. \ \ Assume that
$\lambda\geq4$ and write%

\[
\int_{-\sqrt{\lambda}}^{b}\frac{dv}{\sqrt{\frac{1}{2}(v^{2}-\lambda
)^{2}+(v+\sqrt{\lambda})}}=I_{1}+I_{2}+I_{3},
\]
where
\[
I_{1}:=\int_{-\sqrt{\lambda}}^{-\sqrt{\lambda}+2}\frac{dv}{\sqrt
{v+\sqrt{\lambda}}}=2\sqrt{2},\quad I_{2}:=\int_{-\sqrt{\lambda}+2}%
^{\sqrt{\lambda}-2}\frac{\sqrt{2}dv}{\lambda-v^{2}}=\frac{\sqrt{2}\,\ln
(\sqrt{\lambda}-1)}{\sqrt{\lambda}},
\]
and
\[
I_{3}:=\int_{\sqrt{\lambda}-2}^{b}\frac{dv}{\sqrt{v+\sqrt{\lambda}}}=2\Big
(\sqrt{b+\sqrt{\lambda}}-\sqrt{2\sqrt{\lambda}-2}\,\Big).
\]
It follows from (\ref{a6.3}) that $T<T_{\lambda}$ and therefore (\ref{a6.4})
follows from the definition of $T$ and the assumption that $0<\varepsilon
\leq\varepsilon_{\lambda}$.
\end{proof}

\begin{remark}
Since $\lim_{\lambda\rightarrow\infty}T_{\lambda}=2\sqrt{2}$, it follows from
Lemma~\ref{alem6.1} that for sufficiently large $\lambda>0$, Condition \ref{A}
holds for $0<\varepsilon<\frac{\pi}{6\sqrt{2}}$. \ Also, easy numerical
estimates (assisted by a computer algebra program!) \ show that for any
$\lambda\geq\lambda_{0},$ Condition \ref{A} \ holds if $0<\varepsilon<\frac
{1}{4}.$ \ 
\end{remark}

\bigskip

\section{\label{sec4}Results for small $\varepsilon.$}

\bigskip

\subsection{\label{sec4.1}Asymptotic form of periodic solutions as
$\varepsilon\rightarrow0.$\bigskip}

\ One of the key points made in the work of Angenent, Mallet-Paret, and
Peletier is that as $\varepsilon\rightarrow0,$ \ the solution which we have
called $u_{p},$ or $u_{3},$ \ tends to the lower root, $\ \underline{U}\left(
t\right)  ,$ of $u^{\prime\prime}=0$ \ \ in $[0,\frac{\pi}{2})$ \ and to the
upper root, which we have denoted by $\bar{U}\left(  t\right)  ,$ in
$(\frac{\pi}{2},\pi].$ \ \ (This was for the case $\lambda>\lambda_{0}.)$ We
will need this result here. \ Since one of the aims of this paper is to give
elementary proofs of their results, with no reliance on infinite dimensional
analysis or partial differential equations, we give a new proof. The result is
obtained from the following lemma, where we make no assumption on $\lambda.$
\ For any $\lambda,$ the equation $u^{3}-\lambda u+\cos\left(  t\right)  =0$
\ has a smallest solution $u=\underline{U}\left(  t\right)  $ which is
continuous in any interval $I_{k}=\left[  \frac{\left(  2k-1\right)  \pi}%
{2},\frac{\left(  2k+1\right)  \pi}{2}\right]  $ with $k$ an even integer.
When $k$ is odd, $\bar{U}$ is continuous in $I_{k}$. \ 

\begin{lemma}
\label{lem4} For some even integer $k$ let $J_{k}$ and $M_{k}$ be closed
intervals with $M_{k}\subset int\left(  J_{k}\right)  \subset I_{k}$ (where
$int$ \ denotes interior). \ Then for any $\delta>0,$ there is an
$\varepsilon_{0}$ such that if $0<\varepsilon<\varepsilon_{0},$ \ and if $u$
is a solution of $\left(  1.1\right)  $ with $-b\leq u\leq0$ in $J_{k},$ then
$\left|  u(t)-\underline{U}\left(  t\right)  \right|  <\delta$ in $M_{k}.$ If
$k$ is odd, a similar result holds, stating that solutions which are positive
and bounded by $b$ in $J_{k}$ are close to $\bar{U}$ in $M_{k}$.
\end{lemma}

If we recall that $u_{p}$ is an even function, we see that $u_{p}$ must be
close to the lowest root \underline{$U$}$\left(  t\right)  $ of $u^{\prime
\prime}=0$ on any interval of the form $[0,\left(  1-\delta\right)  \frac{\pi
}{2}]$ \ as $\varepsilon\rightarrow0.$ \ This is true also for any bounded
solution with $u\left(  0\right)  <u_{p}\left(  0\right)  .$

\bigskip

\begin{proof}
Consider the case $k$ even. The branch $\underline{U}\left(  t\right)  $ is
strictly negative in $J_{k}.$ \ Further, given a sufficiently small
$\delta>0,$ depending on $\lambda$ and $J_{k},$ the quantity
\[
f_{\max}=\max_{\substack{\underline{U}\left(  t\right)  -\delta\geq u\\t\in
J_{k}}}f\left(  t,u\right)
\]
is negative and the quantity
\[
f_{\min}=\min_{\substack{\underline{U}\left(  t\right)  +\delta\leq
u\leq0\\t\in J_{k}}}f\left(  t,u\right)
\]
is positive. \ It follows $\varepsilon_{0}$ can be chosen so that if $u\left(
t\right)  \leq\underline{U}\left(  t\right)  -\delta$ for some $t\in M_{k}$
and $u^{\prime}\left(  t\right)  \leq0,$ then $u$ crosses $u=-b$ as $t$
increases in $J_{k}$, while if $u^{\prime}\left(  t\right)  >0,$ then $u$
crosses $-b$ in $J_{k}$ as $t$ decreases. \ Similarly, if $u\left(  t\right)
\geq\underline{U}\left(  t\right)  +\delta,$ then $u$ will cross zero before
$t$ leaves $J_{k},$ either forward or backward depending on the sign of
$u^{\prime}.$ The case $k$ odd is similar. This proves Lemma \ref{lem4}. \ 
\end{proof}

\bigskip

The remaining results in this section give more information about the
asymptotic behavior of solutions as $\varepsilon\rightarrow0.$ \ They will be
needed in the next section. \ \bigskip

\begin{lemma}
\label{lem4'} Assume that $\lambda>\lambda_{0}$. \ Then \bigskip

(i) There are $M_{0}>0$ and $\varepsilon_{0}>0$ such that if $0<\varepsilon
<\varepsilon_{0}$ and $u$ is a solution of (\ref{1.1}) satisfying $u\leq
-\sqrt{\frac{\lambda}{3}}$ over an interval $[c,d],$ with $-\infty\leq
c<d\leq\infty$, then $u(t)<\underline{U}\left(  t\right)  +M_{0}%
\varepsilon^{2}$ in $J_{\varepsilon}:=[c+\frac{3}{K}\varepsilon|\ln
\varepsilon|,d-\frac{3}{K}\varepsilon|\ln\varepsilon|]$. \ If $u$ is a
solution of $\left(  1.1\right)  $ satisfying $u\leq-\sqrt{\frac{\lambda}{3}}$
in $[c,d]$, then $u(t)>\bar{U}\left(  t\right)  -M_{0}\varepsilon^{2}$ in
$J_{\varepsilon}$. Here $K=\sqrt{\bar{U}(0)-\sqrt{\frac{\lambda}{3}}}$,
$[c,d]:=(-\infty,d]$ if $c=-\infty$ and $d<\infty$, and other cases such as
$d=\infty$ and $c>-\infty$ and both $c=-\infty$ and $d=\infty$ are defined similarly.\bigskip

(ii) If $0<\varepsilon\leq\varepsilon_{0}$ and $u$ is a solution of
(\ref{1.1}) which is bounded in $(-\infty,\infty)$ and satisfies
$u<-\sqrt{\frac{\lambda}{3}}$ in $[c,d]$ with $-\infty\leq c<d\leq\infty$,
then
\begin{equation}
|u(t)-\underline{U}(t)|<M_{0}\varepsilon^{2}\quad\mbox{ for }t\in
J_{\varepsilon}. \label{a601}%
\end{equation}
If instead $u$ satisfies $u>\sqrt{\frac{\lambda}{3}}$ over $[c,d]$, then
\begin{equation}
|u(t)-\bar{U}(t)|<M_{0}\varepsilon^{2}\quad\mbox{ for }t\in J_{\varepsilon}.
\label{a602}%
\end{equation}
\end{lemma}

\begin{proof}
It is easy to check that there are $M_{0}>0$ and $\varepsilon_{0}>0$ such that
if $0<\varepsilon\leq\varepsilon_{0}$ and $U_{2}=\underline{U}+(M_{0}%
-1)\varepsilon^{2}$, then
\[
\varepsilon^{2}U_{2}^{\prime\prime}<U_{2}^{3}-\lambda U_{2}+\cos t
\]
for all $t\in(-\infty,\infty)$. Assume for contradiction that there are
$\varepsilon\in(0,\varepsilon_{0}) $ and $\hat{t}\in J_{\varepsilon}$ such
that $u\geq\underline{U}+M_{0}\varepsilon^{2}$ at the point $\hat{t}$. Assume
first that $u^{\prime}(\hat{t})\geq\underline{U}^{\prime}(\hat{t})$. Let
$w=u-U_{2}$. Since $\varepsilon^{2}w^{\prime\prime}>(u^{2}+uU_{2}+U_{2}%
^{2}-\lambda)w$ in $[c,d]$, $w(\hat{t})>0$ and $w^{\prime}(\hat{t})\geq0$, it
follows that $w>0$, $w^{\prime}>0$ and $\varepsilon^{2}w^{\prime\prime}%
>K^{2}w$ for $t\in\lbrack\hat{t},d]$. Then for $\hat{t}\leq t<s\leq d$ we find
that
\[
\frac{d}{ds}\left(  w(s)+\sqrt{w^{2}(s)-w^{2}(t)+\frac{\varepsilon^{2}}{K^{2}%
}(w^{\prime})^{2}(t)}\right)  >\frac{K}{\varepsilon}\left(  w(s)+\sqrt
{w^{2}(s)-w^{2}(t)+\frac{\varepsilon^{2}}{K^{2}}(w^{\prime})^{2}(t)}\right)
\]
and so for $\hat{t}\leq t<T\leq d$ we obtain
\begin{equation}
w(T)+\sqrt{w^{2}(T)-w^{2}(t)+\frac{\varepsilon^{2}}{K^{2}}(w^{\prime})^{2}%
(t)}>[w(t)+\frac{\varepsilon}{K}w^{\prime}(t)]e^{\frac{K}{\varepsilon}(T-t)}.
\label{a600}%
\end{equation}
Assume that $d<\infty$. Evaluating (\ref{a600}) at $T=d$ and $t=\hat{t}$ and
using $d-\hat{t}\geq\frac{3}{K}\varepsilon|\ln\varepsilon|$ gives
\[
w(d)+\sqrt{w^{2}(d)-w^{2}(\hat{t})+\frac{\varepsilon^{2}}{K^{2}}(w^{\prime
})^{2}(\hat{t})}>[w(\hat{t})+\frac{\varepsilon}{K}w^{\prime}(\hat{t})]\frac
{1}{\varepsilon^{3}}%
\]
Since $w(\hat{t})\geq\varepsilon^{2}$ and $w^{\prime}(\hat{t})\geq0$, it
follows that the above inequality does not hold if $\varepsilon$ is small.

\bigskip

Therefore, $\varepsilon_{0}$ can be chosen independent of $\hat{t}$ so that
for $0<\varepsilon\leq\varepsilon_{0}$, $u^{\prime}(\hat{t})<\underline
{U}^{\prime}(\hat{t})$ if $u(\hat{t})\geq\underline{U}+M_{0}\varepsilon^{2}$.
Then let $\hat{u}(t)=u(-t)$ and $\hat{U}_{2}(t)=U_{2}(-t)$. Apply the same
argument for $\hat{u}$ as above to get the same contradiction if
$\varepsilon<\varepsilon_{0}$.\bigskip

If $d=\infty,$ then letting $t\rightarrow\infty$, we see that the right side
of (\ref{a600}) goes to $\infty$, contradicting the boundedness of the left
side of (\ref{a600}). The other assertions of (i) can be proved similarly.\bigskip

(ii) Since $u_{1}$ is the minimal solution of (\ref{1.1}), it follows from (i)
that to prove $\left(  \ref{a601}\right)  $ it suffices to show that
$u_{1}(t)>\underline{U}-M_{0}\varepsilon^{2}$ in $J_{\varepsilon}$. We can
choose $\varepsilon_{0}$ so that for $0<\varepsilon\leq\varepsilon_{0}$,
\[
\varepsilon^{2}U_{1}^{\prime\prime}>U_{1}^{3}-\lambda U_{1}+\cos t,
\]
where $U_{1}=\underline{U}-(M_{0}-1)\varepsilon^{2}$. Assume that there is an
$\varepsilon\leq\varepsilon_{0}$ and $\hat{t}\in J_{\varepsilon}$ such that
$u_{1}(\hat{t})\leq\underline{U}(\hat{t})-M_{0}\varepsilon^{2}$. Then, as
above, by considering $(U_{1}-u_{1})(t)$ for increasing $t$ if $u_{1}^{\prime
}(\hat{t})\leq U_{1}^{\prime}(\hat{t})$ and for decreasing $t$ if
$u_{1}^{\prime}(\hat{t})>U_{1}^{\prime}(\hat{t})$, we obtain that
$(U_{1}-u_{1})(t)\rightarrow\infty$ if $t\rightarrow\pm\infty$ respectively, a
contradiction. This proves (\ref{a601}), and (\ref{a602}) can be proved similarly.
\end{proof}

\bigskip

The next Lemma describes the asymptotic behavior of $u_{1}$, $u_{5}$ and their
first order derivatives as $\varepsilon\rightarrow0$.

\begin{lemma}
\label{alem7.11} Assume that $\lambda>\lambda_{0}$. \ We can choose
$\varepsilon_{0}$ so small that if $0<\varepsilon\leq\varepsilon_{0}$, then
for all $t\in(-\infty,\infty)$,
\begin{equation}
|u_{1}(t)-\underbar{U}(t)|<M_{0}\varepsilon^{2},\quad\mbox{ and }\quad
|u_{5}(t)-\bar{U}(t)|<M_{0}\varepsilon^{2} \label{a7.20}%
\end{equation}
and
\begin{equation}
|u_{1}^{\prime}(t)-\underline{U}^{\prime}(t)|<M_{0}^{\prime}\varepsilon
,\quad\mbox{ and }\quad|u_{5}^{\prime}(t)-\bar{U}^{\prime}(t)|<M_{0}^{\prime
}\varepsilon, \label{a7.21}%
\end{equation}
where $M_{0}^{\prime}$ is a constant independent of $\varepsilon$.
\end{lemma}

\begin{proof}
The inequalities in (\ref{a7.20}) follow from (ii) of the above lemma. To show
(\ref{a7.21}) we let $z=\frac{1}{\varepsilon^{2}}(u_{1}-\underline{U})$, which
satisfies
\begin{equation}
\varepsilon^{2}z^{\prime\prime}=(3\underline{U}^{2}-\lambda)z+3\varepsilon
^{2}\underline{U}z^{2}+\varepsilon^{4}z^{3}-\underline{U}^{\prime\prime}.
\label{a7.24}%
\end{equation}
Further, $z^{\prime}(0)=0$ and $|z(t)|<M_{0}$ for all $t\in(-\infty,\infty)$
and $\varepsilon\in(0,\varepsilon_{0})$. Multiply (\ref{a7.24}) by $z^{\prime
}$ to get
\begin{align}
\frac{1}{2}\varepsilon^{2}(z^{\prime})^{2}(t)  &  =\frac{1}{2}(3\underline
{U}^{2}(t)-\lambda)z^{2}(t)-\frac{1}{2}(3\underline{U}^{2}(0)-\lambda
)z^{2}(0)-3\int_{0}^{t}\underline{U}(s)\underline{U}^{\prime}(s)z^{2}%
(s)\,ds\nonumber\label{a7.25}\\
&  +\varepsilon^{2}\underline{U}(t)z^{3}(t)-\varepsilon^{2}\underline
{U}(0)z^{3}(0)-\int_{0}^{t}z^{3}(s)\underline{U}^{\prime}(s)\,ds+\frac{1}%
{4}\varepsilon^{4}(z^{4}(t)-z^{4}(0))\nonumber\\
&  -\underline{U}^{\prime\prime}(t)z(s)+\underline{U}^{\prime\prime
}(0)z(0)+\int_{0}^{t}z(s)\underline{U}^{\prime\prime\prime}(s)\,ds
\end{align}
Since $z$, $\underline{U}$ and the derivatives of $\underline{U}$ are all
bounded with the bounds independent of $\varepsilon$, it follows that the
right-hand side of (\ref{a7.25}) is bounded by a constant, say, $\frac{1}%
{2}(M_{0}^{\prime})^{2}$ over $[0,\pi]$, that is, $\varepsilon^{2}(z^{\prime
})^{2}(t)\leq(M_{0}^{\prime})^{2}$ for $t\in\lbrack0,\pi]$ and $\varepsilon
\in(0,\varepsilon_{0}]$. Hence (\ref{a7.21}) holds for $t\in\lbrack0,\pi]$.
Since $u_{1}$ and $\underline{U}$ are $2\pi$-periodic and even functions, we
see that (\ref{a7.21}) holds for all $t\in(-\infty,\infty)$.
\end{proof}

\bigskip

\begin{lemma}
\label{alem7.10} Assume that $\lambda>\lambda_{0}$. Then for each small
$\mu>0$ there is an $\varepsilon_{\mu}\in\left(  0,1\right)  $ such that if
$0<\varepsilon<\varepsilon_{\mu}$ and $u$ is a solution with $u>\sqrt
{\frac{\lambda}{3}}$ in $[c,d]$, then
\begin{equation}
|u-u_{5}|+\frac{\varepsilon}{2K}|u^{\prime}-u_{5}^{\prime}|\leq M_{1}%
e^{-\frac{K}{\varepsilon}\mu}\quad\mbox{ in }[c+\mu,d-\mu]. \label{a621}%
\end{equation}
If $u\leq-\sqrt{\frac{\lambda}{3}}$ in $[c,d]$, then
\[
|u(t)-u_{1}(t)|+\frac{\varepsilon}{2K}|u^{\prime}(t)-u_{1}^{\prime}(t)|\leq
M_{1}e^{-\frac{K}{\varepsilon}\mu}\quad\mbox{ in }[c+\mu,d-\mu].
\]
Here $K$ is defined in Lemma~\ref{lem4'} and $M_{1}=2(\bar{U}(\pi)-\sqrt
{\frac{\lambda}{3}})$.
\end{lemma}

\begin{proof}
(i) Let $w=u_{5}-u$. Since $w>0$ and $u>\sqrt{\frac{\lambda}{3}}$ in $[c,d]$,
it follows that $\varepsilon^{2}w^{\prime\prime}>K^{2}w$ in $[c,d]$. From
Lemma~\ref{lem4} we see that $w<M_{0}\varepsilon^{2}$ in $[c+\mu,d-\mu]$. We
first assume that $w^{\prime}\geq0$ at $t_{0}:=c+\mu$. Then (\ref{a600}) holds
for this $w$ with $t_{0}\leq t<T\leq d$, and so for $t_{0}\leq t<T\leq d$,

\bigskip%
\begin{equation}
w(t)+(1-e^{-\frac{K}{\varepsilon}(T-t)})\frac{\varepsilon}{K}|w^{\prime
}(t)|\leq2w(T)e^{-\frac{K}{\varepsilon}(T-t)} \label{a622}%
\end{equation}
Since $w(T)<\bar{U}(\pi)-\sqrt{\frac{\lambda}{3}}$, and $e^{-\frac
{K}{\varepsilon}\mu}<\frac{1}{2}$ if $\varepsilon$ is small, it follows that
for $t\in\lbrack c-\mu,d-\mu]$,
\begin{equation}
w(t)+\frac{\varepsilon}{2K}|w^{\prime}(t)|\leq M_{1}e^{-\frac{K}{\varepsilon
}(T-t)}. \label{a620}%
\end{equation}
Let $T=d$. Since $d-t\geq\mu$ for $t\in\lbrack t_{0},d-\mu]$, (\ref{a621})
follows immediately from (\ref{a620}) for $t\in\lbrack c+\mu,d-\mu]$. \bigskip
Assume that $w^{\prime}(t_{0})<0$. Let $\hat{t}=\sup\{t\in\lbrack t_{0}%
,d-\mu]:w^{\prime}<0\mbox{ in }[t_{0},t)\}$. From what we just proved, we can
assume that $\hat{t}=d-\mu$. Since $w^{\prime}$ is decreasing in $[c,d]$,
$w^{\prime}<0$ in $[c-\mu,\hat{t}]$, and $w>0$ in $[c,d]$, it suffices to show
that (\ref{a621}) holds at $t=t_{0}$. Integrating $w^{\prime\prime}w^{\prime
}<\frac{K^{2}}{\varepsilon^{2}}ww^{\prime}$ over $[c,t_{0}]$ gives
\[
w(t_{0})+\frac{\varepsilon}{K}|w^{\prime}(t_{0})|\leq\Big[w(c)+\sqrt
{w^{2}(c)-w^{2}(t_{0})+\frac{\varepsilon^{2}}{K^{2}}(w^{\prime})^{2}(t_{0}%
)}\,\,\Big]\,e^{-\frac{K}{\varepsilon}\mu}%
\]
which implies that (\ref{a621}) holds at $t=t_{0}$. \bigskip The proof of the
inequality for $u_{1}$ is similar and therefore is omitted.
\end{proof}

\bigskip

From Lemmas~\ref{lem4'}, \ref{alem7.11} and \ref{alem7.10} one immediately
obtains a refinement of Lemma \ref{lem4'}.

\begin{corollary}
\label{acor40} Assume that $\lambda>\lambda_{0}$. For any small $\mu>0$, there
is an $\varepsilon_{\mu}>0$ such that if $0<\varepsilon<\varepsilon_{\mu}$ and
$u$ is a bounded solution of (\ref{1.1}) over $(-\infty,\infty)$ satisfying
$u\leq-\sqrt{\frac{\lambda}{3}}$ over $[c,d]$, then for $t\in\lbrack
c+\mu,d-\mu]$
\begin{equation}
|u(t)-\underline{U}(t)|+\varepsilon|u^{\prime}(t)-\underline{U}^{\prime
}(t)|\leq(M_{0}+M_{0}^{\prime}+1)\varepsilon^{2}. \label{a650}%
\end{equation}
If, instead, $u$ satisfies $u\geq\sqrt{\frac{\lambda}{3}}$ over $[c,d]$, then
for $t\in\lbrack c+\mu,d-\mu],$
\begin{equation}
|u(t)-\bar{U}(t)|+\varepsilon|u^{\prime}(t)-\bar{U}^{\prime}(t)|\leq
(M_{0}+M_{0}^{\prime}+1)\varepsilon^{2}. \label{a651}%
\end{equation}
\end{corollary}

\bigskip

\bigskip The next lemma is a refinement of Lemma~\ref{lem4} if $\lambda
>\lambda_{0}$.

\begin{lemma}
\label{alem41} Assume that $\lambda>\lambda_{0}$. For any small $\mu>0$, there
is an $\varepsilon_{\mu}>0$ such that if $0<\varepsilon<\varepsilon_{\mu}$ and
$u$ is a bounded solution of (\ref{1.1}) which satisfies $u<-\sqrt{\lambda}$
at some point in $J_{k\mu}:=[\frac{(2k-1)\pi}{2}+\frac{\mu}{2},\frac
{(2k+1)\pi}{2}-\frac{\mu}{2}]$ for some even integer $k$, then (\ref{a650})
holds in $M_{k\mu}:=[\frac{(2k-1)\pi}{2}+\mu,\frac{(2k+1)\pi}{2}-\mu]$.
Similarly, (\ref{a651}) holds in $M_{k\mu}$ if $u$ satisfies $u>\sqrt{\lambda
}$ at some point in $J_{k\mu}$ for some odd integer $k$.
\end{lemma}

\begin{proof}
To show the first part of lemma, from Corollary~\ref{acor40} it suffices to
show that for $\varepsilon$ small, $u<-\sqrt{\frac{\lambda}{3}}$ in $J_{k\mu}%
$. Assume that this is false. Then there is a sequence $\varepsilon_{n}$ with
$\lim_{n\rightarrow\infty}\varepsilon_{n}=0$ and a sequence of $t_{n}\in
J_{k\mu}$ with $\lim_{n\rightarrow\infty}t_{n}=\hat{t}$ for some $\hat{t}\in
J_{k\mu}$ such that $u_{n}(t_{n})=-\sqrt{\lambda}$. Then by a phase plane
argument we obtain that $u_{n}$ will reach $b$ by time $\hat{t}+O(\varepsilon
_{n})$ as $n\rightarrow\infty$, contradicting that $u_{n}$ is bounded.
\end{proof}

\bigskip

Since $u_{p}$ is antisymmetric around $\pi/2$, from Lemma~\ref{alem41} we
immediately get

\begin{theorem}
\label{athm100} Assume that $\lambda>\lambda_{0}$. For any small $\mu>0$ there
is an $\varepsilon_{\mu}>0$ such that if $0<\varepsilon<\varepsilon_{\mu}$,
then (\ref{a650}) and (\ref{a651}) hold for $u=u_{p}$ in $[-\frac{\pi}{2}%
+\mu,\frac{\pi}{2}-\mu]$ and $[\frac{\pi}{2}+\mu,\frac{3\pi}{2}-\mu]$ respectively.
\end{theorem}

\bigskip

The next theorem is about the asymptotic behavior of $u_{2}$ as $\epsilon\to0$.

\begin{theorem}
(i) Let $V_{-1}$ be the homoclinic solution of $\ddot{V}=V^{3}-\lambda V-1$
such that $\lim_{\tau\to\pm\infty}V_{-1}(\tau)=\underline{U}(\pi)$ and
$\dot{V}_{-1}(0)=0$. Then $\lim_{\epsilon\to0}u_{2}(\pi)=V_{-1}(0)$. For any
given $T>0$, $u_{2}(\pi+\epsilon\tau)-V_{-1}(\tau)$ approaches zero uniformly
for $\tau\in[-T,0]$ as $\epsilon\to0$. \bigskip(ii) For any small $\mu>0$
there is an $\epsilon_{\mu}>0$ such that if $0<\epsilon<\epsilon_{\mu}$, then
(\ref{a650}) holds in $[0,\pi-\mu]$ for $u=u_{2}$.
\end{theorem}

\begin{proof}
Let $v(\tau)=u_{2}(\pi+\epsilon\tau)$. Then $v$ satisfies
\[
\ddot{v}=v^{3}-\lambda v+cos (\pi+\epsilon\tau), \quad v(0)=u_{2}(\pi),
\quad\dot{v}(0)=0.
\]
Then (i) follows easily by a phase plane argument. To show (ii), we take $T$
so large that $V_{-1}(-T)<-\sqrt{\frac{\lambda}{3}}$. Then from (i) it follows
that if $\epsilon$ is sufficiently small, then $u_{2}(\pi-\epsilon
T)<-\sqrt{\frac{\lambda}{3}}$. Therefore, if $\epsilon$ is small enough to
satisfy $\epsilon T <\frac{\mu}{2}$, then, since $u_{2}^{\prime}>0$ in
$(0,\pi)$, $u_{2}(t)\leq u_{2}(\pi-\frac{\mu}{2})<-\sqrt{\frac{\lambda}{3}}$
for $t\in[-\frac{\mu}{2},\pi-\frac{\mu}{2}]$. Hence (ii) follows from
Corollary~\ref{acor40}.
\end{proof}

\bigskip

\subsection{\label{sec4.2}Further periodic solutions}

\bigskip

As pointed out earlier, the classical analysis of Duffing's equation, with
results such as those in \cite{nm}, \ corresponded in a sense to taking
$\lambda$ large. \ These solutions are bounded as $\lambda\rightarrow\infty.$
\ The solutions $u_{1},u_{2},u_{3}$ found in section \ref{sec2} \ all have
minimum values (which are also their initial values) below $-\sqrt{\lambda},$
while $u_{4}$ and $u_{5}$ have maxima \ above $\sqrt{\lambda}.$ \ In this
section we begin by finding solutions with oscillations, distributed
reasonably evenly in $\left(  0,2\pi\right)  ,$ and with initial values in the
interval $\left(  -\sqrt{\lambda},0\right)  .$ \ This is for a fixed $\lambda$
and small $\varepsilon.$ Then we consider the behavior of those solutions if
we keep the number of oscillations fixed (together with $\lambda)$ \ and let
$\varepsilon\rightarrow0.$ \ It is found that the oscillations collect near
odd multiples of $\frac{\pi}{2},$ forming internal ``layers'' between $\bar
{U}$ and $\underline{U}$. (See Figure 6.) The result is related to one in
\cite{hm2} about a similar nonlinear forced oscillator (derived in \cite{ooj}
), where the nonlinearity was quadratic rather than cubic. \ However the
Duffing equation has a richer collection of solutions than the equation in
\cite{hm2}. \ 

\bigskip%

\begin{center}
\includegraphics[
height=3.0882in,
width=3.9384in
]%
{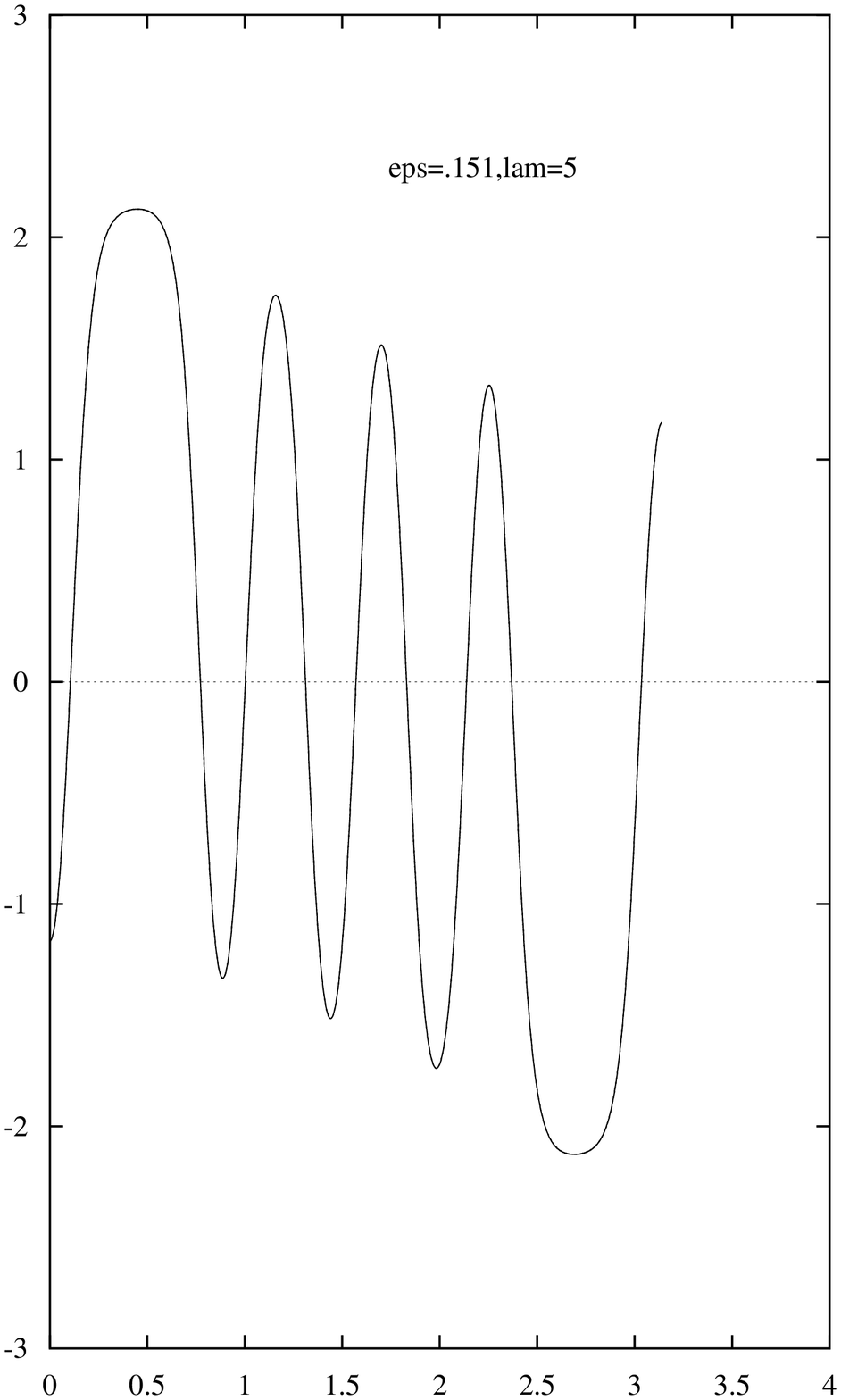}%
\\
Figure 6
\end{center}

\begin{center}
\end{center}

\bigskip

A question of interest is to what extent the results in this paper depend on
the symmetry of the cosine function. \ This is addressed further in section
\ref{sec5}, \ but the techniques in the Theorems \ref{ath7.1} and \ref{ath7.5}
of the current section do not use the symmetry of cosine around $\frac{\pi}%
{2}.$ \ These should apply to more general forcing functions, as we will
discuss in future work. \ On the other hand, the proof of Theorem \ref{ath10}
below shows how use of the symmetry around $\frac{\pi}{2}$ \ can greatly
simplify some proofs. \ The key fact is that if $u$ is a solution of $\left(
\ref{1.1}\right)  $ such that $u\left(  \frac{\pi}{2}\right)  =0,$ then
$u\left(  \frac{\pi}{2}-t\right)  = - u\left(  \frac{\pi}{2}+t\right)  $ for
all $t.$ \ \ 

\bigskip

In some of the following results we will switch back and forth between the
original scaling and that in the slow equation $\left(  \ref{slow1}\right)  ,$
with $t=\varepsilon\tau,$ $u\left(  t\right)  =v\left(  \tau\right)  .$
\ \ For any $\kappa\in(0,1]$ , \ the equation
\begin{equation}
\ddot{V}=V^{3}-\lambda V+\kappa\label{new1}%
\end{equation}
has a unique solution, forming the homoclinic orbit, with $V\left(  \pm
\infty\right)  =\bar{U}\left(  \arccos\kappa\right)  $ and $\dot{V}\left(
0\right)  =0.$ We denote this solution by $V_{\kappa}.$ Let
\[
\Lambda=\inf\{\lambda>\lambda_{0}:\,V_{1}(0)<U_{0}(\pi)\}
\]
Thus, \ for $\lambda>\Lambda,$ the homoclinic orbit at $\tau=0$ ($\kappa=1)$
has a minimum below $U_{0}\left(  \pi\right)  .$ It is easy to show that
$\Lambda<3.$

\begin{theorem}
\label{ath7.1} Suppose that $\lambda>\Lambda.$ \ Then for any $N>0$ \ there is
an $\varepsilon_{0}$ such that if $0<\varepsilon<\varepsilon_{0},$ and if $M$
is any positive integer with $M\leq N,$ then $\left(  \ref{1.1}\right)  $ has
a periodic solution with exactly $M$ local maxima in $\left(  0,\pi\right)  ,$
the last being at $\pi.$ Also, these periodic solutions have a local minimum
at $t=0$ and satisfy $-\sqrt{\lambda}<u(0)<U_{0}(\pi)$ and $U_{0}%
(\pi)<u\left(  \pi\right)  <\bar{U}(\pi)$.
\end{theorem}

\begin{proof}
We will need the following lemma, which is similar to Theorem 2 in \cite{hm2}.

\begin{lemma}
\label{alem7.2} Let $u$ be a solution of (\ref{1.1}) with $u^{\prime}(0)=0$.
If $t_{1}$ and $t_{2}$ are two successive minima (maxima) of $u$ with
$0<t_{1}<t_{2}<\pi$, then $u(t_{2})<u(t_{1})$.
\end{lemma}

\begin{proof}
We only show the case that $u(t_{1})$ and $u(t_{2})$ are successive maxima of
$u$ and the other case can be proved similarly. Assume that $u$ reaches its
maximum between $t_{1}$ and $t_{2}$ at $t=t_{m}$. Then multiply (\ref{1.1}) by
$u^{\prime}$ and integrate over $[t_{-},t_{+}]$, where $t_{1}<t_{-}%
<t_{m}<t_{+}<t_{2}$ and $u(t_{-})=u(t_{+})$, to yield
\begin{align*}
\frac{1}{2}\varepsilon^{2}[(u^{\prime})^{2}(t_{+})-(u^{\prime})^{2}(t_{-})]
&  =\int_{t_{-}}^{t_{+}}u^{\prime}(t)\cos t\,dt=(\int_{t_{-}}^{t_{m}}%
+\int_{t_{m}}^{t_{+}})u^{\prime}(t)\cos t\,dt\\
&  =\int_{u(t_{-})}^{u(t_{m})}(\cos t_{-}(u)-\cos t_{+}(u))\,du>0.
\end{align*}
Here we use that $\cos t$ is decreasing in $(0,\pi),$ and $t_{-}(u)$ and
$t_{+}(u)$ are the inverse function of $u=u(t)$ for $t\in\lbrack t_{1},t_{m}]$
and $t\in\lbrack t_{m},t_{2}]$ respectively. It then follows that $u^{\prime
}(t_{-})=0$ before $u^{\prime}(t_{+})=0$, which implies that $u(t_{1}%
)>u(t_{2})$ as required.

\ 
\end{proof}

For $\lambda>\lambda_{0}$, recall that $\underbar{U}(t)<U_{0}(t)<\bar{U}(t)$
are the three solutions of $u^{3}(t)-\lambda u(t)+\cos t=0$. \ For given
$\lambda>\Lambda$, let $\delta=U_{0}(\pi)-V_{1}(0).$ \ We also need the fact,
easily proved, that if $V_{\kappa}$ is the homoclinic solution of $\left(
\ref{ns1}\right)  $ defined earlier, for $0<\kappa\leq1,$ \ then the minimum
values $V_{\kappa}\left(  0\right)  $ are decreasing with respect to $\kappa.$

\bigskip

Also, for each $\alpha,$ let $v_{\alpha}$ be the solution of $\left(
\ref{slow1}\right)  $ such that $v\left(  0\right)  =\alpha,$ $v^{\prime
}\left(  0\right)  =0.$ \ Let $\hat{\alpha}=U_{0}\left(  \pi\right)
-\frac{\delta}{2}$. \ 

\bigskip

Choose $N>0.$ When $\varepsilon=0,$ $v_{\hat{\alpha}}$ is periodic, and so has
$N+2$ local maxima in some interval $\left[  0,T_{N}\right]  .$ For
sufficiently small $\varepsilon,$ $v_{\hat{\alpha}}$ still has at least $N+1$
local maxima in $\left(  0,\frac{\pi}{\varepsilon}\right)  .$ \ Hence
$u_{\hat{\alpha}}$ has at least $N+1$ maxima in $\left(  0,\pi\right)  .$
\ Suppose that these are at $0<t_{1}<t_{2}<...<t_{N+1}<\pi.$ \ \ Further, from
Lemma \ref{alem7.2} \ we see that $u_{\hat{\alpha}}\left(  t_{i+1}\right)
<u_{\hat{\alpha}}\left(  t_{i}\right)  $ for $0<i\leq N.$

\bigskip

If any of the maxima $t_{i}$ for $i<N$ are such that $u_{\hat{\alpha}}\left(
t_{i}\right)  <\underline{U}\left(  t_{i}\right)  ,$ then because
\underline{$U$} \ is increasing in $\left(  0,\pi\right)  ,$ we must have
$u<\underline{U}$ in $\left(  t_{i},\pi\right)  $ (if the solution exists out
to $\pi),$ and so $u_{\hat{\alpha}}$ could not have any more maxima. \ This
shows that for $1\leq i\leq N,$ the maxima of $u_{\hat{\alpha}}$ must lie in
the interval $\left[  U_{0}\left(  t\right)  ,\bar{U}\left(  t\right)
\right]  .$ \ However, if $u^{\prime}=u^{\prime\prime}=0$ at some $t$\ in
$\left(  0,\pi\right)  ,$ then $u^{\prime\prime\prime}=-\sin\left(  t\right)
<0,$ so this could not be a maximum. \ Hence the first $N$ maxima of
$u_{\hat{\alpha}}$ must lie in $\left(  U_{0}\left(  t\right)  ,\bar{U}\left(
t\right)  \right)  .$

\bigskip

\ The solution $u_{\hat{\alpha}}$ must have a minimum at $s_{0}=0$\ and
further minima $s_{1},....s_{N}\in\left(  0,\pi\right)  $. \ Also $\hat
{\alpha}>u\left(  s_{1}\right)  >u\left(  s_{2}\right)  >...>u\left(
s_{N}\right)  $, \ by Lemma \ref{alem7.2}. \ Therefore, $u_{\hat{\alpha}%
}\left(  t_{i}\right)  \geq u_{\hat{\alpha}}\left(  s_{i-1}\right)
+\frac{\delta}{2}>u_{\hat{\alpha}}\left(  s_{i}\right)  +\frac{\delta}{2}$ for
$i=1,...,N.$ \ 

\bigskip

Now decrease $\alpha$ from $\hat{\alpha}.$ Since $u^{\prime\prime}<0$ at a
maximum, the maxima are continuous in $\alpha,$ and so remain above $U_{0}$ as
long as they exist. \ If $\alpha=-b,$ $u_{\alpha}$ has no local minima. \ But
as long as there are, say, $M\leq N$ local maxima in $\left(  0,\pi\right)  $
with $u>U_{0}$ and with $\alpha<\hat{\alpha}$ these maxima and their
intervening minima are separated, in that $u_{\alpha}\left(  t_{i}\left(
\alpha\right)  \right)  -u_{\alpha}\left(  s_{i}\left(  \alpha\right)
\right)  \geq\frac{\delta}{2}$ for $1\leq i\leq M.$ This means that the number
of maxima cannot decrease until one crosses $t=\pi.$ This must happen
successively for each maximum, which proves Theorem \ref{ath7.1}. \ \ 
\end{proof}

\begin{remark}
With further estimates it can be shown that for small $\varepsilon$ there are
$2\pi$ periodic solutions with as many as $\frac{K}{\varepsilon} $ maxima in
$(0,\pi)$, where $K>0$ is independent of $\varepsilon$.
\end{remark}

\bigskip

In the next theorem we shall describe the asymptotic behavior, as
$\varepsilon\rightarrow0$, of $2\pi$-periodic solutions of (\ref{1.1}) with
$m$ maxima in $(0,\pi]$ and with $u\left(  0\right)  <0,$ $u\left(
\pi\right)  >0$. \ The result, informally, is that if $m$ is kept fixed as
$\varepsilon\rightarrow0$ then the internal maxima and minima all tend to
$\frac{\pi}{2},$ while near $0$ and $\pi$ \ the solutions tend to homoclinic
orbits of the appropriate limiting equation. \ Thus these solutions have
``spikes'' at $0$ and $\pi$, and internal layers where they are close to
appropriate heteroclinic orbits, near $\frac{\pi}{2}.$ \ In theorem
\ref{ath10} we obtain $2\pi$-periodic solutions without spikes at $0$ and
$\pi,$ using the symmetry of the cosine. \ In forthcoming work we expect to
consider other multilayer solutions, combining spike and non-spike behaviors
at multiples of $\pi,$ without reliance on symmetry. \ 

\bigskip

\begin{theorem}
\label{ath7.5} Suppose that $\lambda>\Lambda.$ \ For an integer $N\geq2$
choose $\varepsilon_{0}$ as in Theorem \ref{ath7.1}. \ Let $m$ be an integer
with $1\leq m\leq N,$ and for each $\varepsilon\in\left(  0,\varepsilon
_{0}\right)  $ let $u=u_{\varepsilon}$ be a $2\pi$ periodic solution of
(\ref{1.1}) satisfying $u^{\prime}(0)=u^{\prime}(\pi)=0$ and $-\sqrt{\lambda
}<u(0)<U_{0}(0)$, and having exactly $m$ maxima in $(0,\pi]$. Let
$t_{1},\cdots,t_{m}$ and $s_{1},\cdots,s_{m-1}$ be the maxima and the minima
of $u$ in $(0,\pi]$ respectively such that $0<t_{1}<s_{1}<t_{2}<s_{2}%
<\cdots<t_{m-1}<s_{m-1}<t_{m}=\pi$. Let $V_{1}$ be the homoclinic solution of
$\ddot{v}=v^{3}-\lambda v+1$ such that $\lim_{\tau\rightarrow\pm\infty}%
v(\tau)=\bar{U}(0)$ and $\dot{v}(0)=0$. Suppose that $\delta$ is a positive
number with $\delta<\bar{U}(0)-\sqrt{\frac{\lambda}{3}}$. Let $V_{-1}$ be the
homoclinic solution of $\ddot{v}=v^{3}-\lambda v-1$ such that $\lim
_{\tau\rightarrow\pm\infty}v(\tau)=\underline{U}(\pi)$ and $\dot{v}(0)=0$. Let
$V_{0^{+}}$ be the heteroclinic solution of $\ddot{v}=v^{3}-\lambda v$ with
$\lim_{\tau\rightarrow-\infty}v(\tau)=\bar{U}(\frac{\pi}{2})$, $\lim
_{\tau\rightarrow\infty}v(\tau)=\underline{U}(\frac{\pi}{2})$ and
$v(0)=\bar{U}(\frac{\pi}{2})-\delta$. Also, let $V_{0^{-}}\left(  t\right)
=-V_{0^{+}}\left(  -t\right)  .$ Then, as $\varepsilon\rightarrow0$,

(i) $u_{\varepsilon}(0)\rightarrow V_{1}(0)$, $u_{\varepsilon}(\pi)\rightarrow
V_{-1}(0)$.

(ii) for each $1\leq j\leq m-1$, $t_{j}\rightarrow\frac{\pi}{2}$,
$s_{j}\rightarrow\frac{\pi}{2}$, $u_{\varepsilon}(t_{j})\rightarrow\bar
{U}(\frac{\pi}{2})$, and $u_{\varepsilon}(s_{j})\rightarrow\underline{U}%
(\frac{\pi}{2})$.

(iii) for any given $T>0$, $|u_{\varepsilon}(\varepsilon\tau)-V_{1}%
(\tau)|+|\varepsilon u_{\varepsilon}^{\prime}(\varepsilon\tau)-\dot{V}%
_{1}(\tau)|)\rightarrow0$ uniformly for $\tau\in\lbrack0,T]$, and
$|u_{\varepsilon}(\pi+\varepsilon\tau)-V_{-1}(\tau)|+|\varepsilon
u_{\varepsilon}^{\prime}(\pi+\varepsilon\tau)-\dot{V}_{-1}(\tau)|)=0$
uniformly for $\tau\in\lbrack-T,0]$.

Also, for $1\leq j\leq m-1$, let
\[
T_{j}=\sup\{t\in(t_{j},\pi):\,u_{\varepsilon}>u_{5}-\delta\quad\mbox
{ in }(t_{j},t]\}
\]
and
\[
S_{j}=\sup\{t\in(s_{j},\pi):\,u_{\varepsilon}<u_{1}+\delta\quad\mbox
{ in }(s_{j},t]\}.
\]
Then $0<T_{j}-t_{j}<\frac{3\varepsilon}{K}|\ln\varepsilon|$ and $0<S_{j}%
-s_{j}<\frac{3\varepsilon}{K}|\ln\varepsilon|$, where $K$ is defined in
Lemma~\ref{lem4'}. Further, for any given $T>0$, as $\varepsilon\rightarrow0,$
$|u_{\varepsilon}(T_{j}+\varepsilon\tau)-V_{0^{+}}(\tau)|+|\varepsilon
u_{\varepsilon}^{\prime}(T_{j}+\varepsilon\tau)-\dot{V}_{0^{+}}(\tau
)|\rightarrow0$ and $|u_{\varepsilon}(S_{j}+\varepsilon\tau)+V_{0^{-}}%
(\tau)|+|\varepsilon u_{\varepsilon}^{\prime}(T_{j}+\varepsilon\tau)+\dot
{V}_{0^{-}}(\tau)|\rightarrow0,$ uniformly for $\tau\in\lbrack0,T]$.

\bigskip

Finally, for any small $\mu>0$ there is an $\epsilon_{\mu}>0$ such that if
$0<\epsilon<\epsilon_{\mu}$, then (\ref{a651}) and (\ref{a650}) hold in
$[\mu,\frac{\pi}{2}-\mu]$ and $[\frac{\pi}{2}+\mu,\pi-\mu]$ respectively.
\end{theorem}

\bigskip

\begin{proof}
\bigskip

In the proof we shall suppress the dependence of $u$ on $\varepsilon$. \ We
first show that $\lim_{\varepsilon\rightarrow0}u(0)=V_{1}(0)$. \ Assume that
this is false. Since $u(0)$ is bounded, there is a sequence $\varepsilon_{n}$
with $\lim_{n\rightarrow\infty}\varepsilon_{n}=0$ such that as $n\rightarrow
\infty$, the corresponding $u(0)$ approaches $\beta_{0}\in\lbrack
-\sqrt{\lambda},U_{0}(0)]$. For simplicity we assume that $\lim_{\varepsilon
\rightarrow0}u(0)=\beta_{0}$. Let $z(\tau)=u(\varepsilon\tau)$. Then $z$
satisfies
\[
\ddot{z}=z^{3}-\lambda z+\cos\varepsilon\tau,\quad z(0)=u(0),\quad\dot
{z}(0)=0.
\]
Let $z_{0}$ be the solution of
\[
\ddot{z}_{0}=z_{0}^{3}-\lambda z_{0}+1,\quad z_{0}(0)=\beta_{0},\quad\dot
{z}_{0}(\tau)=0.
\]
Since $\lim_{\varepsilon\rightarrow0}z(0)=\beta_{0}$, it follows from the
continuity of solutions with respect to parameters that for any given $T>0$
\ if $z_{0}$ is defined on $[0,T]$, then $\lim_{\varepsilon\rightarrow0}%
z(\tau)=z_{0}(\tau)$ uniformly for $\tau\in\lbrack0,T]$.\bigskip

Assume that $\beta_{0}<V_{1}(0)$. Then there is a $T_{0}>0$ such that $\dot
{z}_{0}(\tau)>0$ for $\tau\in(0,T_{0})$ and $\lim_{\tau\rightarrow T_{0}^{-}%
}z_{0}(\tau)=\infty$. It follows that for $\varepsilon$ sufficiently small $z$
crosses $b$, which is impossible. Hence $\beta_{1}>V_{1}(0)$. Then $z_{0}$ is
a periodic function with a period $\tilde{T}_{0}>0$ and so $z_{0}$ has $m+1$
maximum in the interval $(0,(m+2)\tilde{T}_{0})$. Hence by continuity, for
$\varepsilon$ sufficiently small, $z(\tau)$ also has $m+1$ maximum in the
interval $(0,(m+2)\tilde{T}_{0})$ , which implies that $u$ has $m+1$ maxima in
$(0,\varepsilon(m+2)\tilde{T}_{0})\subset(0,\pi]$ for $\varepsilon$
sufficiently small, contradicting the assumption on $u$. Therefore,
$\lim_{\varepsilon\rightarrow0}u(0)=V_{1}(0)$ and so (iii) and the first part
of (i) follow. The rest of (i) and (iii) can be proved similarly.

\bigskip

We next show that $\lim_{\varepsilon\rightarrow0}t_{1}=\frac{\pi}{2}$. Suppose
not. Since $t_{1}\in(0,\pi]$, there is a sequence $\varepsilon_{n}$ with
$\lim_{n\rightarrow\infty}\varepsilon_{n}=0$ such that $\lim_{n\rightarrow
\infty}t_{1}=\bar{t}\neq\frac{\pi}{2}$ for some $\bar{t}\in\lbrack0,\pi]$.
Again, for simplicity, we assume that $\lim_{\varepsilon\rightarrow0}%
t_{1}=\bar{t}$.

\bigskip

Since $u(\pi)\rightarrow V_{-1}(0)$ as $\varepsilon\rightarrow0$, it follows
from Lemma~\ref{alem7.2} that $\ \bar{t}\neq\pi$. We suppose now that
$\frac{\pi}{2}<\bar{t}<\pi$. Let $z(\tau)=u(t_{1}+\varepsilon\tau)$. Then $z$
satisfies
\[
\ddot{z}=z^{3}-\lambda z+\cos(t_{1}+\varepsilon\tau),\quad z(0)=u(t_{1}%
),\quad\dot{z}(0)=0.
\]
For $b=b\left(  \lambda\right)  =\sqrt{\lambda+\frac{1}{2\lambda}}$ (as in
Lemma~\ref{alem6.1}), let
\[
T_{b}=\int_{-b}^{u(t_{1})}\frac{dy}{\sqrt{\frac{1}{2}(y^{4}-u(t_{1}%
)^{4})-\lambda(y^{2}-u(t_{1})^{2})+2(y-u(t_{1}))\cos t_{1}}},
\]
and
\[
T=\sup\{\tau\in(0,T_{b}+1):\,z>-b,\quad\dot{z}<0\quad\mbox{ in }(0,\tau)\,\}.
\]

Observe that $\bar{U}(t_{1})>z(0)=u(t_{1})>\bar{U}(t_{1})-\delta$ if
$\varepsilon$ is sufficiently small. Hence, $T_{b}$ is defined and
$T_{b}<\infty$; $\ddot{z}(0)>0$ and so $T$ is well defined. Assume that
$\varepsilon$ is so small that $t_{1}+\varepsilon(T_{b}+1)\leq\pi$. \ \ Then
on $\left(  0,T\right)  $ we have $\ddot{z}\leq z^{3}-\lambda z+\cos t_{1}$
and so $(\dot{z})^{2}>\frac{1}{2}(z^{4}-z^{4}(0))-\lambda(z^{2}-z^{2}%
(0))+2(z-z(0))\cos t_{1}$. We see that $T<T_{b}$ and $\dot{z}(T)<0$.
Therefore, by the definition of $T$, it follows that $u(t_{1}+\varepsilon
T)=z(T)=-b$, which is impossible since $u(\tau)>-b$ for all $\tau$.

\bigskip

We assume now that $\bar{t}\in\lbrack0,\frac{\pi}{2})$. \ Note that
$u(t_{1})<\bar{U}(t_{1})$. We shall show that $\lim_{\varepsilon\rightarrow
0}u(t_{1})=\bar{U}(\bar{t})$. \ For if this is false, then the boundedness of
$u(t_{1})$ implies that there is a sequence $\varepsilon_{n},$ with
$\varepsilon_{n}\rightarrow0$ as $n\rightarrow\infty,$ such that
$\lim_{n\rightarrow\infty}u(t_{1})<\bar{U}(\bar{t})$. It then follows by
arguments similar to those above that $u$ has more than $m$ maxima in
$(0,\pi]$ for large $n$, which is a contradiction.

\bigskip

We now claim that $|T_{1}-t_{1}|<\frac{3\varepsilon}{K}|\ln\varepsilon|$ for
$\varepsilon>0$ sufficiently small, where $T_{1}$ is defined in the statement
of this Theorem $\left(  \ref{ath7.5}\right)  $. Assume that this is false.
Let $w=u_{5}-u$. Then in $(t_{1},T_{1})$, $w>0$ and $\varepsilon^{2}%
w^{\prime\prime}\geq K^{2}w$. Since $u^{\prime}<0$ and $u_{5}^{\prime}>0$ just
to the right of $t_{1,}$ we see that $w^{\prime}>0$ and $w^{\prime\prime
}w^{\prime}>(K^{2}/\varepsilon^{2})ww^{\prime}$ in $(t_{1},T_{1}]$. Then, as
in $\left(  \ref{a600}\right)  $(where $w$ was slightly different), \ we again
get
\begin{equation}
w(T)+\sqrt{w^{2}(T)-w^{2}(t)+\frac{\varepsilon^{2}}{K^{2}}(w^{\prime})^{2}%
(t)}>[w(t)+\frac{\varepsilon}{K}w^{\prime}(t)]e^{\frac{K}{\varepsilon}(T-t)}.
\label{a100a}%
\end{equation}
\bigskip

We assume now that $\bar{t}>0$. Since $\lim_{\varepsilon\rightarrow0}\bar
{U}^{\prime}(t_{1})=\bar{U}^{\prime}(\bar{t})=\frac{\sin\bar{t}}{3\bar{U}%
^{2}(\bar{t})-\lambda}$, we use $\left(  \ref{a7.21}\right)  $ to obtain%
\[
w^{\prime}(t_{1})=u_{5}^{\prime}(t_{1})\geq\bar{U}^{\prime}(t_{1}%
)-M_{0}^{\prime}\varepsilon\geq\frac{1}{2}\bar{U}^{\prime}(\bar{t}%
)-M_{0}^{\prime}\varepsilon>\frac{\sin\bar{t}}{4(3\bar{U}^{2}(\bar{t}%
)-\lambda)}>0.
\]
\bigskip

This estimate gives a lower bound on $w^{\prime}\left(  t_{1}\right)  $ as
$\varepsilon$ tends to zero. Evaluating (\ref{a100a}) at $t=t_{1}$ and
$T=t_{1}+\frac{3\varepsilon}{K}|\ln\varepsilon|$ and letting $\varepsilon
\rightarrow0$ shows that the right side of (\ref{a100a}) goes to infinity,
while the left side is bounded. \ \ This contradiction proves the claim if
$\bar{t}\in\left(  0,\frac{\pi}{2}\right)  .$

\bigskip

We now assume that $\bar{t}=0$. We observe that
\begin{align*}
w^{\prime}(t_{1}+\frac{\varepsilon}{K}|\ln\varepsilon|)  &  >u_{5}^{\prime
}(t_{1}+\frac{\varepsilon}{K}|\ln\varepsilon|)\geq\frac{1}{2}\bar{U}^{\prime
}(t_{1}+\frac{\varepsilon}{K}|\ln\varepsilon|)-M_{0}^{\prime}\varepsilon\\
&  >\frac{\sin(\frac{\varepsilon}{K}|\ln\varepsilon|)}{4(3\bar{U}%
^{2}(0)-\lambda)}>\frac{\frac{\varepsilon}{K}|\ln\varepsilon|}{2\pi(3\bar
{U}^{2}(0)-\lambda)}.
\end{align*}

\bigskip

Then evaluating (\ref{a100a}) at $t=t_{1}+\frac{\varepsilon}{K}|\ln
\varepsilon|$ and $T=t_{1}+\frac{3\varepsilon}{K}|\ln\varepsilon|$ and letting
$\varepsilon\rightarrow0$, the same contradiction will be obtained. This shows
that $|T_{1}-t_{1}|<\frac{3\varepsilon}{K}|\ln\varepsilon|,$ as claimed. In
particular, $T_{1}-t_{1} \to0$ and so $T_{1}\to\bar{t}$.

\bigskip

Continuing with our proof that $\bar{t}=\frac{\pi}{2},$ and under the
assumption that this is false and $\bar{t}\in\lbrack0,\frac{\pi}{2}),$ let
$z(\tau)=u(T_{1}+\varepsilon\tau)$. Then $z$ satisfies
\begin{equation}
\ddot{z}=z^{3}-\lambda z+\cos(T_{1}+\varepsilon\tau),\quad z(0)=u(T_{1}%
)=u(T_{1})-\delta,\quad\dot{z}(0)=\varepsilon u^{\prime}(T_{1}).
\end{equation}
We shall show that $\lim_{\varepsilon\rightarrow0}\dot{z}(0)=\dot{z}_{0}(0)$,
where $z_{0}$ is the homoclinic solution of
\begin{equation}
\ddot{z}_{0}=z_{0}^{3}-\lambda z_{0}+\cos(\bar{t})
\end{equation}
with $\ z_{0}(0)=\bar{U}(\bar{t})-\delta$ and $\dot{z}_{0}(0)<0.$

\bigskip

Suppose that $\lim_{\varepsilon\rightarrow0}\dot{z}(0)\neq\dot{z}_{0}(0)$.
Since $\dot{z}(0)<0$ and $\dot{z}(0)$ is bounded, which is easily verified,
there is sequence of values of $\varepsilon,$ and corresponding solutions $u$
with corresponding $z,$ such that $\dot{z}\left(  0\right)  $ approaches a
number $\sigma\neq\dot{z}_{0}(0)$ with $\sigma\leq0$. We first assume that
$\sigma>\dot{z}_{0}(0)$. Let $z_{1}$ be the periodic solution, with period
$\tilde{T},$ to
\begin{equation}
\ddot{z}_{1}=z_{1}^{3}-\lambda z_{1}+\cos(\bar{t}),\quad z_{1}(0)=\bar{U}%
(\bar{t})-\delta,\quad\dot{z}_{1}(0)=\sigma. \label{a405}%
\end{equation}
Since $T_{1}\to\bar{t}$, $z(\tau)$ approaches $z_{1}(\tau)$ as $\varepsilon
\rightarrow0$ uniformly in compact intervals of $\tau,$ and it follows that
for sufficient small $\varepsilon$, $z$ oscillates more than $m+1$ times in
$[0,(m+2)\tilde{T}]$, and so $u$ has more than $m$ maxima in $[0,\pi]$, a contradiction.

\bigskip

Hence we can assume that $\sigma<\dot{z}_{0}(0)$ and again assume that $z_{1}$
solves (\ref{a405}). Then there is a $\hat{T}>0$ such that $\dot{z}_{1}%
(\tau)<0$ for $\tau\in\lbrack\hat{T},0]$ and $z_{1}(\hat{T})>b$, which implies
that for $\varepsilon$ sufficiently small, $u(t_{1}+\varepsilon\hat{T})>b$,
again a contradiction.

\bigskip

Therefore, as $\varepsilon\rightarrow0$, $\dot{z}(0)\rightarrow\dot{z}_{0}(0)$
and $z$ goes to $z_{0}$ uniformly in any compact interval. Since $z_{0}$ is
homoclinic to $\bar{U}(\bar{t})$, it follows from continuity that after
$T_{1}$, $u$ will return to the any given neighborhood of $\bar{U}(\bar{t})$
before $t=T_{1}+M\varepsilon$ for some $M$ independent of $\varepsilon$. Then,
since $m=2$, $u$ increases and stays close to $\bar{U}$ till $\pi$, which
implies that $\lim_{\varepsilon\rightarrow0}u(\pi)=\bar{U}(\pi)$ and so
$u(t_{2})=u(\pi)>u(t_{1})$ for $\varepsilon$ small, contradicting
Lemma~\ref{alem7.2}. This shows that $\bar{t}=\frac{\pi}{2}$ for $m=2$.

\bigskip

Then, $z_{0}$ is the heteroclinic solution connecting $\bar{U}(\frac{\pi}{2})$
as $t\rightarrow-\infty$ to $\underline{U}(\frac{\pi}{2})$ as $t\rightarrow
\infty$. So, as $\varepsilon\rightarrow0$ the point $\left(  t,u(t)\right)  $
reaches a point as close to $(\frac{\pi}{2},\underline{U}\left(  \frac{\pi}%
{2}\right)  )$ as we like, and then Lemma \ref{alem7.2} \ implies that
$\lim_{\varepsilon\rightarrow0}s_{1}=\frac{\pi}{2}$ and $\lim_{\varepsilon
\rightarrow0}u\left(  s_{1}\right)  =\underline{U}(\frac{\pi}{2})$. For
$t\in\lbrack t_{1},\pi],$ $u\left(  t\right)  <u\left(  t_{1}\right)  $, again
by Lemma \ref{alem7.2}. \ Further, $u(t)$ remains close to $\underline{U}(t)$
until $t$ is close to $\pi,$ for otherwise there would be similar
contradictions to those obtained above. Hence, $u$ must, after its last
minimum, follow $\underline{U}$ until close to $\pi$ and then , since its last
maximum is above $U_{0}(\pi)$, follow $V_{-1}$. The bound on $S_{j}-s_{j}$
follows in the same way as the bound for $T_{j}-t_{j}$. For the final
statement in the Theorem, we use (i)-(iii) and Corollary \ref{acor40}. This
proves the theorem for $m=2.$ \ The proof for $m>2$ is similar.
\end{proof}

\bigskip

The solutions discussed in Theorem \ref{ath7.5} have upward and downward
pointing ``spikes'' at each multiple of $\pi,$ corresponding to homoclinic
orbits in the phase plane, as well as ``layers'' at odd multiples of
$\frac{\pi}{2}$ corresponding to heteroclinic orbits. \ In the next theorem we
use the symmetry of cosine around $\frac{\pi}{2}$ to give a quick proof that
there are also solutions with layers but without the spikes. \ In our future
work we expect to show that these solutions exist without reliance on
symmetry, but it is also valuable, we believe, to show how quickly a proof can
be obtained in the symmetric case. We need a preliminary lemma.

\bigskip

\begin{lemma}
\label{alem30} Assume that $\lambda>\lambda_{0}$. Let $u$ be a solution of
(\ref{1.1}) with $u(\frac{\pi}{2})=0$ and $u^{\prime}(\frac{\pi}{2})=\beta$.

(i) If $\beta<-\frac{\lambda}{\sqrt{2}\varepsilon}$, then either $u^{\prime
}<0$ in $[\frac{\pi}{2},\pi]$ or there is $t_{0}\in(\frac{\pi}{2},\pi)$ such
that $u(t_{0})=\underline{U}(t_{0})$ and $u^{\prime}(t)<0$ in $[\frac{\pi}%
{2},t_{0}]$.

(ii) Let $H(u)=\lambda u^{2}-u^{4}/2$. Assume that $\beta<-\frac
{\sqrt{H(\underline{U}(\pi))}}{\varepsilon}$. If $t_{1}\in(\frac{\pi}{2},\pi)$
is the first time that $u^{\prime}=0$ (if such a point exists), then
$u(t_{1})<\underline{U}(\pi)$.
\end{lemma}

\begin{proof}
From (\ref{1.1}) we have
\begin{equation}
\varepsilon^{2}(u^{\prime})^{2}(t)=\varepsilon^{2}(u^{\prime})^{2}(\frac{\pi
}{2})-H(u)+2\int_{\frac{\pi}{2}}^{t}u^{\prime}\cos s\,ds. \label{a200}%
\end{equation}

(i) Let $t_{0}=\sup\{t\in(\frac{\pi}{2},\pi):\,u^{\prime}<0\mbox
{ and }u>-\sqrt{\lambda}\mbox{ in }(\frac{\pi}{2},t)\,\}$. It follows from
(\ref{a200}) that for $t\in(\frac{\pi}{2},t_{0})$,
\[
\varepsilon^{2}(u^{\prime})^{2}(t)>\varepsilon^{2}\beta^{2}-H(u)>\varepsilon
^{2}\beta^{2}-\frac{\lambda^{2}}{2}>0,
\]
from which the assertion in (i) follows.

(ii) From (\ref{1.1}), we first have $-\sqrt{\lambda}<u(t_{1})<0$. Evaluate
(\ref{a200}) at $t=t_{1}$ to give%

\[
0=\varepsilon^{2}\beta^{2}-H\left(  u\left(  t_{1}\right)  \right)
+2\int_{\frac{\pi}{2}}^{t_{1}}u^{\prime}\cos\left(  s\right)  \,ds
\]
so that
\[
H\left(  u\left(  t_{1}\right)  \right)  >\varepsilon^{2}\beta^{2}.
\]
This implies that $u\left(  t_{1}\right)  <\underline{U}\left(  \pi\right)  $,
\ as desired.
\end{proof}

\begin{theorem}
\label{ath10} For any given integer $m$, there is an $\varepsilon_{m}>0$ such
that for $\varepsilon\in(0,\varepsilon_{m})$, (\ref{1.1}) has a solution $u$
such that $u^{\prime}(0)=u^{\prime}(\pi)=0$, $u(\frac{\pi}{2})=0$, $u^{\prime
}(\frac{\pi}{2})<0$, and $u$ has $m$ minima and $m$ maxima in $[\frac{\pi}%
{2},\pi]$. \ If we denote these minima and maxima by $s_{j}$ and $t_{j}$
respectively, with $1\leq j\leq m,$ then
\begin{align*}
&  \frac{\pi}{2}<s_{1}<t_{1}<s_{2}<t_{2}<\cdots<s_{m}<t_{m}=\pi,\\
&  \lim_{\varepsilon\rightarrow0}s_{1}=\cdots=\lim_{\varepsilon\rightarrow
0}s_{m}=\frac{\pi}{2},\quad\lim_{\varepsilon\rightarrow0}u(s_{1})=\cdots
=\lim_{\varepsilon\rightarrow0}u(s_{m})=\underline{U}(\frac{\pi}{2}),\\
&  \lim_{\varepsilon\rightarrow0}t_{1}=\cdots=\lim_{\varepsilon\rightarrow
0}t_{m-1}=\frac{\pi}{2},\quad\lim_{\varepsilon\rightarrow0}u(t_{1}%
)=\cdots=\lim_{\varepsilon\rightarrow0}u(t_{m-1})=\bar{U}(\frac{\pi}{2}),\\
&  u^{\prime}(t)>0,\quad\underline{U}(\frac{\pi}{2})<u(t)<\underline{U}%
(\pi)\quad\mbox{ for }s_{m}<t<\pi.
\end{align*}
\end{theorem}

\begin{proof}
Let $u_{\beta}$ denote the solution $u$ of (\ref{1.1}) satisfying $u(\frac
{\pi}{2})=0$ and $u^{\prime}(\pi)=\beta$. Choose $\beta_{1}\in(-\frac{\lambda
}{\sqrt{2}\varepsilon},-\frac{\sqrt{H(\underline{U}(\pi))}}{\varepsilon})$ and
$\beta_{2}<-\frac{\lambda}{\sqrt{2}\varepsilon}$. By a phase plane argument,
there is an $\varepsilon_{m}$ such that for $\varepsilon<\varepsilon_{m}$, the
solution $u_{\beta_{1}}$ has at least $m$ minima and $m$ maxima in $(\frac
{\pi}{2},\pi)$ with all its first $m$ minima lying between $\underline{U}$ and
$U_{0}$ and all its first $m$ maxima lying between $U_{0}$ and $\bar{U}$.
Denote the $m$-th minimum and maximum by $s_{m}$ and $t_{m}$. We consider the
change in these minima as $\beta$ decreases from $\beta_{1}$. Part (ii) of
lemma~\ref{alem30} shows that all the minima in $(\frac{\pi}{2},\pi]$ lie
below the line $u=\underline{U}(\pi)$, and so they can neither pass through
$t=\pi$ nor disappear at the middle branch $U_{0}$ of $u^{\prime\prime}=0$.
Since $u_{\beta_{2}}$ does not have any minimum at all, it follows that as we
decrease $\beta_{1}$ to $\beta_{2}$, all the minima will disappear by crossing
the lower branch of $u^{\prime\prime}=0$ before $t=\pi$.

\bigskip

Let $\beta_{3}=\inf\left\{  \bar{\beta}|s_{m}\text{ is defined continuously
for }\beta\in(\bar{\beta},\beta_{1})\text{ as the m }^{th}\text{ \ minimum
after }\frac{\pi}{2}\right\}  .$ \ \ Then $\ u_{\beta_{3}}^{\prime}(s_{m})=0$
and $u_{\beta_{3}}(s_{m})=\underline{U}(s_{m})$ . For $\beta-\beta_{3}>0$
small, the $m$-th maximum $t_{m}(\beta)$ of $u_{\beta}$ exists. Clearly
$\lim_{\beta\rightarrow\beta_{3}}t_{m}(\beta)=s_{m}(\beta_{3}).$ As before, if
$u_{\beta}(t_{m})=\underline{U}(t_{m})$, we get $u_{\beta}^{\prime}(t)<0$ for
all $t-t_{m}$ small, contradicting the definition of $t_{m}$. \ Hence
$u_{\beta}(t_{m}(\beta))<\underline{U}(t_{m}(\beta)).$ Therefore, as we raise
$\beta$ from $\beta_{3}$, $t_{m}(\beta)$ has to move toward to $\pi$ under the
lower branch $\underline{U}$ of $u^{\prime\prime}=0$. This maximum
$t_{m}\left(  \beta\right)  $ cannot disappear by merger with another minimum
in $\left[  \frac{\pi}{2},\pi\right]  $, \ since this minimum would have to
lie above $u_{\beta}\left(  s_{m}\right)  ,$ contradicting Lemma
\ref{alem7.2}. \ Since $u_{\beta_{1}}(t_{m})>U_{0}(\pi)>\underline{U}(\pi)$,
there must be a point of discontinuity of the $m$th maximum before reaching
$\beta_{1},$ $\ $and so there is a $\beta_{4}\in(\beta_{3},\beta_{1})$ such
that $t_{m}\left(  \beta\right)  $ is continuous in $\left[  \beta_{3}%
,\beta_{4}\right]  $, $t_{m}(\beta_{4})=\pi$ and $u_{\beta_{4}}^{\prime}%
(t_{m}(\beta_{4}))=0$. Then $u_{\beta_{4}}$ is the desired solution. The
asymptotic formulas stated in this theorem can be proved in a similar way to
that of Theorem~\ref{ath7.5}.
\end{proof}

\bigskip

Figure 7 shows graphs of two of the new solutions, with one or three
``layers'' near $\frac{\pi}{2},$ for $\varepsilon=0.5,\lambda=3,$ together
with $u_{2}$ \ and the function $u_{0}\left(  t\right)  .$ It is of interest
that these solutions behave in a very different way from the solution $u_{p}$
discussed earlier. \ While $u_{p}$ increases near $\frac{\pi}{2}$ \ and
decreases near $\frac{3\pi}{2},$ \ these solutions behave in the opposite way. \ 

\bigskip

\ \ The existence of these solutions, and many more, with different numbers of
oscillations between different zeros of $\cos x,$ was conjectured
independently by H. Matano (private communication). \ This was for a more
general class of equations, but the result was limited to a finite interval,
so that ``chaos'' was not involved. \ \ Professor Matano calls the solutions
with the opposite behavior from $u_{p}$ ``up wind'' solutions, \ a graphic
term which we adopted in the caption of Figure 8. \ Similar solutions were
obtained for a different equation by Nakashima \cite{nak}. \ We wish to thank
Professor Matano for helpful correspondence, and in particular for sharing his
conjecture with us. \ While the results in this paper, about periodic
solutions with multiple internal layers, were obtained independently, our
thoughts concerning non-periodic solutions with more than three internal
layers near odd multiples of $\frac{\pi}{2}$ were previously somewhat vague,
and we have been inspired to pursue this topic further by Professor Matano's
conjecture. \ We expect that a proof of this conjecture using shooting
methods, at least for some class of forcing functions, will be part of our
next paper. The methods proposed by Matano for obtaining these solutions are
very different from ours.

\bigskip

Nakashima, in \cite{nak}, \ has also studied the stability of the oscillating
solutions for her equation, including the dimensions of the unstable manifolds
(Morse Index). \ This has implications for the dimension of the global
attractor for the problem $\left(  \ref{1.2}\right)  -\left(  \ref{1.3}%
\right)  .$ Results of this type were also obtained in \cite{hm2}, for a
different equation, though they were not stated in these terms. \ We expect to
study these topics for the current problem in future work, using the methods
of \cite{hm2}, \ which are very different from the approach taken by
Nakashima. \ 

{%
\begin{center}
\includegraphics[
height=2.405in,
width=5.1785in
]%
{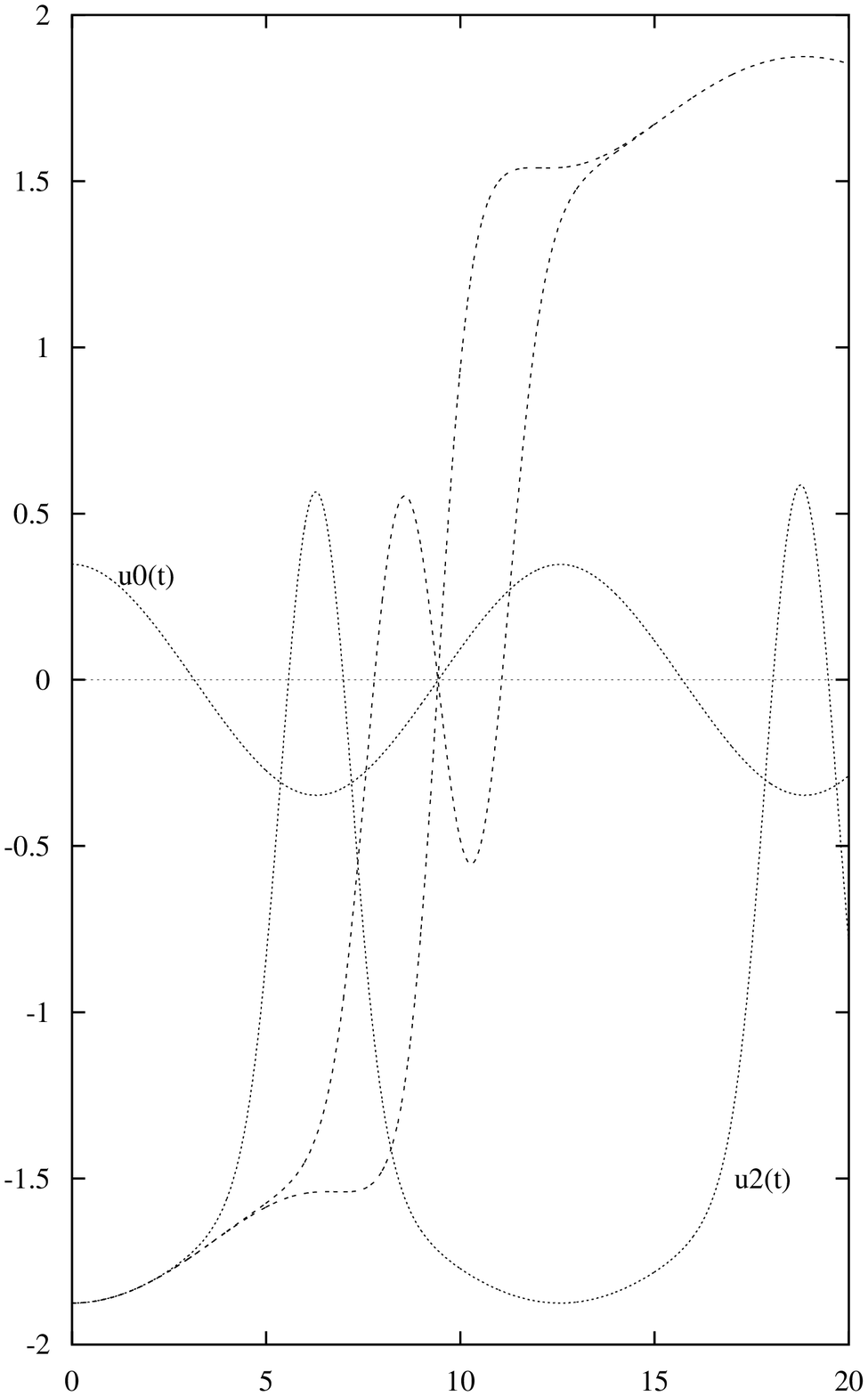}%
\\
Figure 7
\end{center}
}

\begin{center}
Graphs of two ``up wind'' solutions, together with $u_{2}$ and $U_{0}.$
\end{center}

We end this section with a numerically computed solution suggesting how the
different kinds of solutions can be combined. \ In this figure, $\varepsilon
=1,$ which explains why the oscillations of the new solution are so small.
\ Results about this kind of solution will appear in part II. \ 

{%
\begin{center}
\includegraphics[
height=2.2174in,
width=5.4933in
]%
{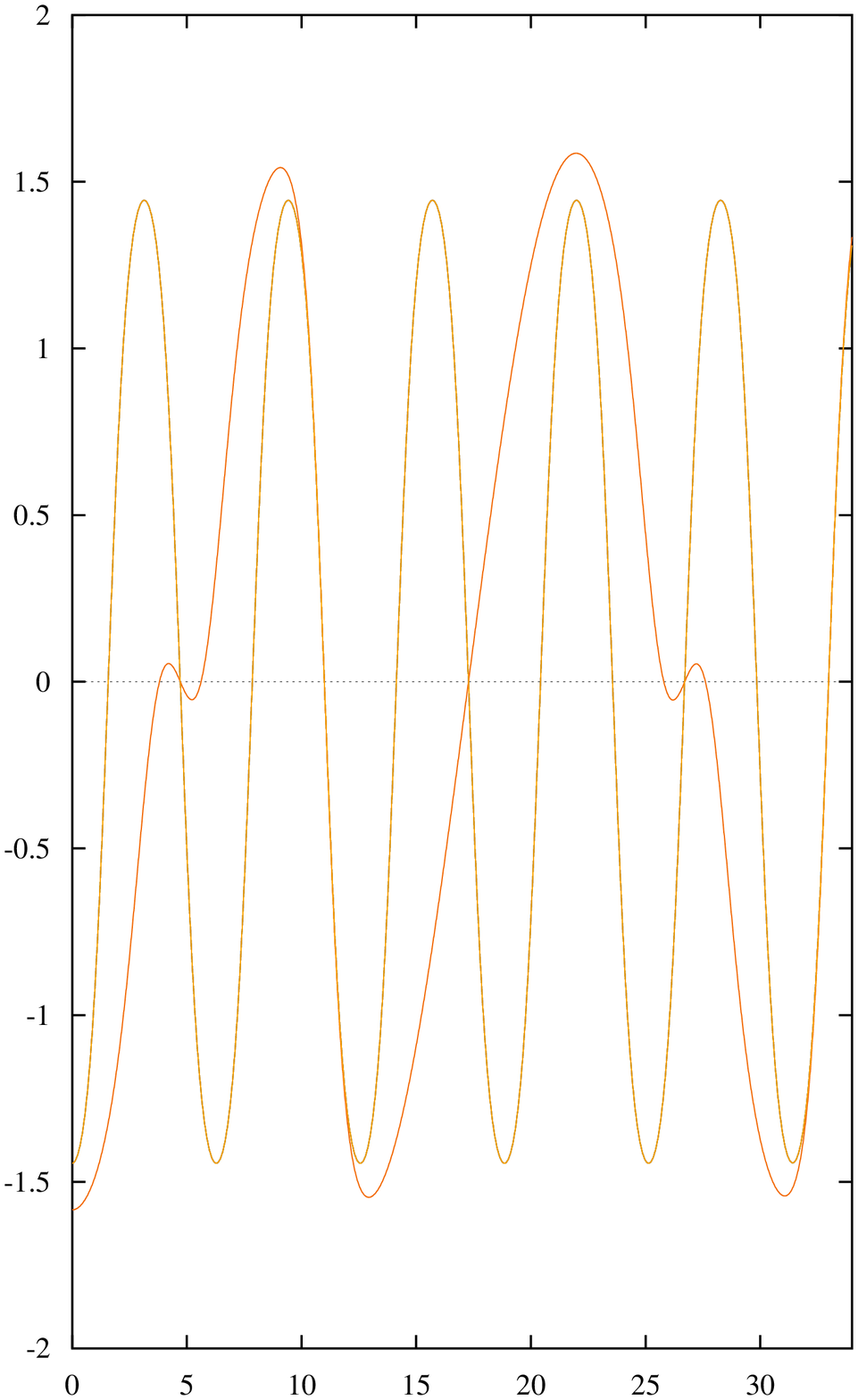}%
\\
Figure 8.
\end{center}
}

\begin{center}
$u_{p}$ and a new solution with three ``up wind'' sections
\end{center}

\subsection{\label{sec4.3}\bigskip Isolation and stability of u$_{p}$}

\bigskip

In this section we show how a shooting method can give results similar to
those in \cite{ampp} for $\left(  \ref{1.1}\right)  ,$ without the use of
infinite dimensional analysis or abstract dynamical systems. \ Using standard
ode methods we prove the linearized stability of the three solutions found in
section \ref{sec1} with respect to $\left(  \ref{1.1}\right)  -\left(
\ref{1.2}\right)  $, and the isolation of $u_{p}$ from other solutions
satisfying $u^{\prime}\left(  0\right)  =u^{\prime}\left(  k\pi\right)  =0.$
\ Full nonlinear stability follows by standard methods, laid out explicitly
for this problem \ in \cite{bf}. \ However, stability is inherently an
infinite dimensional problem, and we certainly do not claim to prove
(nonlinear) stability without the use of functional analysis. \ 

\begin{theorem}
\label{thm8}Suppose that $\lambda>\lambda_{0}.$ \ Choose $r_{1}$ with
$0<r_{1}<\frac{1}{2}.$ Let
\[
\mu=\int_{\left(  1-r_{1}\right)  \pi}^{\left(  1-\frac{r_{1}}{2}\right)  \pi
}\sin s\,ds.
\]
Choose $\delta>0$ so that
\begin{align*}
\left(  i\right)  \text{ \ \ \ }\delta &  <\underline{U}\left(  \frac{r_{1}%
}{2}\pi\right)  -\underline{U}\left(  0\right) \\
\left(  ii\right)  \text{ \ \ If }\bar{U}\left(  t\right)  -\delta &
<u\left(  t\right)  <\bar{U}\left(  t\right)  \text{ for some }t,\text{\ then
}\left|  u^{3}-\lambda u+\cos t\right|  <\mu.
\end{align*}
Let $S_{\delta}$ denote the set of points $\left(  t,u\right)  $ such that
$0\leq t\leq\pi$, $-b\leq u\leq b,$ and
\begin{align}
-b  &  \leq u\leq\underline{U}\left(  t\right)  +\delta\text{ if }0\leq t\leq
r_{1}\pi\nonumber\\
\bar{U}\left(  t\right)  -\delta &  \leq u\leq b\text{ if }\left(
1-r_{1}\right)  \pi\leq t\leq\pi.\nonumber
\end{align}
Then for sufficiently small $\varepsilon$, \ $u_{p}$ is the only $2\pi
$-periodic solution of $\left(  \ref{1.1}\right)  $ with $u^{\prime}\left(
0\right)  =u^{\prime}\left(  \pi\right)  =0$ whose graph over $\left[
0,\pi\right]  $ lies in $S_{\delta}.$ \ Further, $\ u_{p}|_{[0,\pi]}$ is a
stable attractor for the problem $\left(  \ref{1.2}\right)  -\left(
\ref{1.3}\right)  .$with $L=\pi$ \ \ 

\bigskip

Also, suppose that $S_{\delta}$ is extended to $\left[  0,2\pi\right]  $ by
the operation $\left(  t,u\right)  \rightarrow\left(  2\pi-t,u\right)  ,$ and
then periodically to $[0,\infty).$ \ Then for any positive integer $k,$
$u_{p}$ is the unique solution of $\left(  \ref{1.1}\right)  $ satisfying
$u^{\prime}\left(  0\right)  =u^{\prime}\left(  k\pi\right)  =0$ whose graph
lies in the extended set $S_{\delta},$ and its restriction to $\left[
0,k\pi\right]  $ is a stable attractor for $\left(  \ref{1.2}\right)  -\left(
\ref{1.3}\right)  $ with $L=k\pi.$ \ 
\end{theorem}

To compare this with the result in \cite{ampp} \ we notice that the width of
the vertical strip in $S_{\delta}$ where the solution can increase from near
\underline{$U$} \ to near $\bar{U}$ is $\left(  1-2r_{1}\right)  \pi,$ and
this is free to be chosen within the constraint $0<r_{1}<\frac{1}{2},$ \ By
contrast, \cite{ampp} state only that there is some strip, with width
independent of $\varepsilon$, which contains the ``internal layer'' of the
solution we have denoted by $u_{p,}$ and no other solution satisfying the
boundary conditions\ has a jump upward within this strip. \ 

\bigskip

The proof in \cite{ampp} \ is by a detailed construction of $u_{p}$ using
sub-and super-solutions. \ It also uses abstract results from dynamical
systems (to get the uniqueness.) \ Our proof of uniqueness is more direct, and
starts with the existence of $u_{p}$ as given in Theorem $\ref{thm02a}.$

\bigskip

\begin{proof}
To prove Theorem \ref{thm8} \ we consider the variational equation and initial
conditions satisfied by $v=\frac{\partial u_{\alpha}}{\partial\alpha}.$ These
are
\begin{align}
\varepsilon^{2}v^{\prime\prime}  &  =\left(  3u^{2}-\lambda\right)
v\label{st1}\\
v\left(  0\right)   &  =1,\,\,v^{\prime}\left(  0\right)  =0.\nonumber
\end{align}
We will also be concerned with $w=u^{\prime},$ \ which satisfies
\begin{align}
\varepsilon^{2}w^{\prime\prime}  &  =\left(  3u^{2}-\lambda\right)  w-\sin
t\label{st2}\\
w\left(  0\right)   &  =0,w^{\prime}\left(  0\right)  =u^{\prime\prime}\left(
0\right) \nonumber
\end{align}
We observe that $w^{\prime}\left(  0\right)  >0$ when $u\left(  0\right)
\in\left(  \underline{U}\left(  0\right)  ,0\right)  .$ Multiplying $\left(
\ref{st1}\right)  $ by $w$ \ and $\left(  \ref{st2}\right)  $ by $v,$
\ subtracting, \ integrating by parts and using the initial conditions on $v$
and $w,$ we obtain
\begin{equation}
wv^{\prime}-vw^{\prime}|_{t}=-u^{\prime\prime}\left(  0\right)  +\int_{0}%
^{t}\frac{v\left(  s\right)  }{\varepsilon^{2}}\sin s\,\,ds. \label{st3}%
\end{equation}

\begin{lemma}
\label{aa} \ If $u$ is a solution with $u^{\prime}(0)=0$ which remains in
$S_{\delta}$ \ on $\left[  0,\pi\right]  ,$ then $v>0$ on $\left[
0,\pi\right]  $ \ and $v^{\prime}\left(  \pi\right)  >0.$ \ 
\end{lemma}

\begin{proof}
Because $u<\underline{U}\left(  t\right)  +\delta$ \ over the interval
$\left[  0,r_{1}\pi\right]  $, and therefore $3u^{2}-\lambda>0$ \ in this
interval, $\left(  \ref{st1}\right)  $ implies that $v$ grows exponentially
large. \ More precisely, there are positive numbers $K_{1}\geq1$ and $\gamma,$
independent of $\varepsilon,$ \ such that $v\geq K_{1}e^{\frac{\gamma
}{\varepsilon}t}$ in $\left[  0,r_{1}\pi\right]  .$ Also, $u^{\prime\prime
}\left(  0\right)  =O\left(  \frac{1}{\varepsilon^{2}}\right)  $\ as
$\varepsilon\rightarrow0.$ Therefore, for small $\varepsilon,$ \ the right
side of $\left(  \ref{st3}\right)  $ is positive as long after $t=r_{1}\pi$
\ as $v$ is positive (up to $t=\pi).$ \ 

\bigskip

We now show that $u^{\prime}\left(  \left(  1-\frac{r_{1}}{2}\right)
\pi\right)  >0.$ \ If not, then $u^{\prime}<0$ on $\left(  \left(
1-\frac{r_{1}}{2}\right)  \pi,\pi\right)  ,$ since $u^{\prime\prime}<0$ when
$\bar{U}-\delta<u<\bar{U}$, and $\bar{U}^{\prime}>0$ in $(\frac{\pi}{2},\pi)$.
Therefore $u\left(  \pi\right)  <u\left(  \left(  1-\frac{r_{1}}{2}\right)
\pi\right)  $. From (i) it follows that $u\left(  \pi\right)  <\bar{U}\left(
\pi\right)  -\delta$, \ contradicting the assumption that $u$ \ remains in
$S_{\delta}.$ \ Hence, $u^{\prime}\left(  \left(  1-\frac{r_{1}}{2}\right)
\pi\right)  >0.$ \ This implies that $u^{\prime}>0$ on $(0,\left(
1-\frac{r_{1}}{2}\right)  \pi],$ for if not, then $u$ would have a minimum in
this interval, and this minimum would lie above $u\left(  0\right)  ,$
contradicting Lemma \ref{alem7.2}.

\bigskip

Now suppose in $\left(  \ref{st3}\right)  $ that $v=0$ somewhere in $\left(
0,\left(  1-\frac{r_{1}}{2}\right)  \pi\right]  .$ Then the right side of
$\left(  \ref{st3}\right)  $ is positive, while the left side is negative.
\ Hence, $v>0$ on $\left[  0,\left(  1-\frac{r_{1}}{2}\right)  \pi\right]  .$ \ 

\bigskip

It is possible that $v^{\prime}$ becomes negative somewhere in $\left[
0,\pi\right]  .$ Indeed, numerically this is seen to happen. However, we will
show that $v^{\prime}\left(  \left(  1-\frac{r_{1}}{2}\right)  \pi\right)
>0.$ Suppose that $v^{\prime}\left(  \left(  1-\frac{r_{1}}{2}\right)
\pi\right)  \leq0.$ Then $v^{\prime}((1-r_{1})\pi)<0$, for otherwise
$v^{\prime}$ would be positive and increasing on $[(1-r_{1})\pi,(1-\frac
{r_{1}}{2})\pi)$, because $(3u^{2}-\lambda) >0$ there. Thus $v$ is positive
but decreasing on $\left[  \left(  1-r_{1}\right)  \pi,\left(  1-\frac{r_{1}%
}{2}\right)  \pi\right]  .$ \ From $\left(  \ref{st3}\right)  $ we obtain
\begin{align*}
&  v^{\prime}\left(  \left(  1-\frac{r_{1}}{2}\right)  \pi\right)  u^{\prime
}\left(  \left(  1-\frac{r_{1}}{2}\right)  \pi\right)  -v\left(  \left(
1-\frac{r_{1}}{2}\right)  \pi\right)  u^{\prime\prime}\left(  \left(
1-\frac{r_{1}}{2}\right)  \pi\right) \\
&  =-u^{\prime\prime}\left(  0\right)  +\int_{0}^{\left(  1-r_{1}\right)  \pi
}\frac{v\left(  s\right)  }{\varepsilon^{2}}\sin s\,ds+\int_{\left(
1-r_{1}\right)  \pi}^{\left(  1-\frac{r_{1}}{2}\right)  \pi}\frac{v\left(
s\right)  }{\varepsilon^{2}}\sin s\,ds\\
&  \geq v\left(  (1-\frac{r_{1}}{2})\pi\right)  \int_{\left(  1-r_{1}\right)
\pi}^{\left(  1-\frac{r_{1}}{2}\right)  \pi}\frac{\sin s}{\varepsilon^{2}%
}\,ds.
\end{align*}
Since $u^{\prime}\left(  \left(  1-\frac{r_{1}}{2}\right)  \pi\right)  >0,$ we
get a contradiction from condition $\left(  ii\right)  $ in the statement of
the theorem, and so $v^{\prime}$ could not remain negative on $\left[  \left(
1-r_{1}\right)  \pi,\left(  1-\frac{r_{1}}{2}\right)  \pi\right]  $.
\ \ Therefore, at $\left(  1-\frac{r_{1}}{2}\right)  \pi$ we have both $v$ and
$v^{\prime}>0,$ and since $3u^{2}-\lambda>0$ on $\left[  \left(  1-\frac
{r_{1}}{2}\right)  \pi,\pi\right]  ,$ both remain positive out to $\pi.$ This
proves the lemma.
\end{proof}

\ \ 

\bigskip

In particular, this applies to the solution $u_{p}.$ It follows by standard
stability arguments \cite{bf} that the solution $u_{p}$ is a stable attractor
for the problem $\left(  \ref{1.2}\right)  -\left(  \ref{1.3}\right)  $ on
$\left[  0,\pi\right]  $. \ 

\bigskip

To prove that $u_{p}$ is the only $2\pi$-periodic solution which remains in
$S_{\delta},$ it is convenient to truncate the nonlinearity in $\left(
\ref{1.1}\right)  $ by replacing $u^{3}-\lambda u$ \ by $b^{3}-\lambda b$ for
all $u\geq b,$ and similarly, if $u\leq-b,$ \ replace $u^{3}-\lambda u$ in
$\left(  \ref{1.1}\right)  $ with $\lambda b-b^{3}.$ This means that all
solutions $u_{\alpha}$ exist on $\left[  0,\pi\right]  ,$ and we can consider
$u_{\alpha}^{\prime}\left(  \pi\right)  $ \ to be defined continuously for all
negative $\alpha$. \ \ We also note that if $u_{\alpha}$ is a solution which
does leave the region $\left[  -b,b\right]  ,$ then from the point where
$\left|  u_{\alpha}\right|  =b,$ $\left|  u_{\alpha}^{\prime}\right|  $
continuous to increase and cannot satisfy $u_{\alpha}^{\prime}=0.$ \ 

\bigskip

Recall that there is a unique $\alpha_{1}$ such that $u_{\alpha_{1}}$ is
$2\pi$-periodic and lies entirely below $-\sqrt{\frac{\lambda}{3}};$ this is
the ``minimal'' bounded solution, and its graph does not lie in $S_{\delta}.$
In considering the possibility of a second $2\pi$-periodic solution, besides
$u_{p},$ which lies in $S_{\delta},$ we need only consider $\alpha\in\left(
\alpha_{1},\underline{U}\left(  0\right)  +\delta\right)  .$

\bigskip

\bigskip The solution $u_{p}$ remains in the interior of $S_{\delta}$ on
$\left[  0,\pi\right]  ,$ and the same is true for $u_{\alpha}$ \ if $\left|
\alpha-\alpha_{p}\right|  $ is sufficiently small. \ However as we raise or
lower $\alpha$ from $\alpha_{p},$ \ we reach values where $u_{\alpha}$ leaves
$S_{\delta}$ in $\left[  0,\pi\right]  .$ Let $I=\left[  \hat{\beta}%
,\hat{\alpha}\right]  $ be the maximal interval containing $\alpha_{p}$ such
that $u_{\alpha}$ remains in $S_{\delta}$ on $\left[  0,\pi\right]  $ for all
$\alpha\in I.$ $I$ is well-defined and closed because $S_{\delta}$ is a closed
set. \ 

\bigskip

\begin{lemma}
\label{lema} If $\alpha\notin I,$ \ then $u_{\alpha}$ does not remain in
$S_{\delta}$ on $\left[  0,\pi\right]  .$
\end{lemma}

\begin{proof}
: \ We claim that $u_{\hat{\beta}}$ leaves $S_{\delta}$ at $\left(  \pi
,\bar{U}\left(  \pi\right)  -\delta\right)  $ and $u_{\hat{\alpha}}$ \ exits
$S_{\delta}$ \ at $\left(  \pi,b\right)  .$ \ \ If not, then one of these
solutions is tangent to the boundary of $S_{\delta}$ at some $t<\pi.$ For
example, a tangency could occur at $\left(  \left(  1-r_{1}\right)  \pi
,\bar{U}\left(  \left(  1-r_{1}\right)  \pi\right)  \right)  .$ But then, no
matter what the slope of $u$ is at this point, a phase plane argument shows
that for sufficiently small $\varepsilon,$ $u$ must cross $-b$ on one
direction or the other, and so nearby solutions also leave $S_{\delta}$ before
$t=\pi,$ \ contradicting the definition of $\hat{\alpha}$ or $\hat{\beta}.$
Similar considerations apply at any other possible tangent point in $[0,\pi).$

\bigskip

Consider the case of $u_{\hat{\beta}}.$ By Lemma \ref{aa}, \ $\frac{\partial
}{\partial\alpha}u_{\alpha}\left(  \pi\right)  |_{\alpha=\hat{\beta}}>0.$
Hence, for $\alpha<\hat{\beta}$ and close to $\hat{\beta},$ $u_{\alpha}\left(
\pi\right)  <\bar{U}\left(  \pi\right)  -\delta,$ \ and we cannot, as we lower
$\alpha,$ find a lower $\alpha$ where $u_{\alpha}$ remains in $S_{\delta}$ and
$u_{\alpha}\left(  \pi\right)  =\bar{U}\left(  \pi\right)  -\delta,$ for at
the first such point we would have $\frac{\partial}{\partial\alpha}u_{\alpha
}\left(  \pi\right)  \leq0,$ a contradiction to Lemma \ref{aa}. \ Similar
remarks apply to $u_{\hat{\alpha}},$ completing the proof of Lemma \ref{lema}. \ \ 
\end{proof}

\bigskip

The proof that if $u=u_{p},$ then $v>0$ on $\left[  0,\pi\right]  $ and
$v^{\prime}\left(  \pi\right)  >0,$ shows that in some neighborhood of
$\alpha_{p},$ $\left(  \alpha-\alpha_{p}\right)  u_{\alpha}^{\prime}\left(
\pi\right)  >0.$ \ \bigskip Now, if there is an $\check{\alpha}\in(\alpha
_{p},\hat{\alpha})$ with $u_{\check{\alpha}}^{\prime}(\pi)=0$, then choose the
smallest such $\check{\alpha}$. By Lemma \ref{aa} we get $\frac{\partial
}{\partial{ \alpha}}u_{\alpha}^{\prime}(\pi)>0 $ at $\alpha_{p}$ and
$\check{\alpha}$, a contradiction because these are adjacent zeros of
$u_{\alpha}^{\prime}(\pi)$. This contradiction can be reached similarly in the
interval $\left(  \hat{\beta},\alpha_{p}\right)  $, completing the proof of
the uniqueness of $u_{p}$ among solutions with period $2\pi$ which remain in
$S_{\delta}.$ \ \ 

\bigskip

\ To extend the uniqueness and stability statements to larger intervals we
note that starting with $v\left(  \pi\right)  >0,$ $\ v^{\prime}\left(
\pi\right)  >0$, \ the same analysis allows us to show inductively that $v$
remains positive, and $v^{\prime}$ \ is positive at any multiple of $\pi.$
\ The arguments about uniqueness and stability can then be extended to
$\left[  0,k\pi\right]  $ completing the proof of Theorem \ref{thm8}. \ 
\end{proof}

\bigskip

\ \ 

\subsection{\label{sec4.4}Sensitivity with respect to initial conditions}

As this paper is already quite long, we will content ourselves with a few
remarks. \ Up until now there has been no mention of ``horseshoes'' in this
paper, or of Poincar\'{e} maps, because our technique is to follow complete
solutions of the ode, rather than to take snapshots at regular intervals.
\ The results, however, are related to standard dynamical systems concepts
such as horseshoes and sensitivity to initial conditions.

\bigskip

In Theorems \ref{thm4} \ and \ref{thm5} we obtain a weak kind of sensitivity
to initial conditions. \ For each sequence there is a corresponding solution,
but the relation established is not 1:1. \ There could be an interval of
initial values $\alpha$ \ in which $u_{\alpha}$ intersects the same sequence
of $w_{k}$. \ This corresponds to a so-called ``topological'' horseshoe,
without the hyperbolicity that was a key feature of Smale's original
derivation. \ (See $\cite{gh}$ \ for general discussion and references.) It is
noted, however, that the solutions in Theorem \ref{thm4} corresponding to the
sequences of all odd integers and of all even integers, are our periodic
solutions $u_{1}$ and $u_{5},$ \ and for these, uniqueness of the
correspondence is established. \ The difficulty is the stability, or at least
hyperbolicity, of solutions $u_{2},u_{3},$ and $u_{4.}$

\bigskip

In Section \ref{sec4.3}, \ hyperbolicity was established for $u_{3}$ when
$\varepsilon$ is sufficiently small.\ Theorem \ref{thm8} implies that any
solution which corresponds to a sequence with no $2$ \ or $4$ \ is isolated
from any other such solution. Therefore, in Theorem \ref{thm8}, if we consider
a sequence with no $2$ or $4$, then the infinite intersection of closed
intervals used in the construction contains exactly one point. \ This means
that the set of solutions found in Theorem \ref{thm4} corresponding to
sequences chosen from the set $\{1,3,5\}$\ is, for sufficiently small
$\varepsilon,$ \ in 1:1 correspondence with the set of allowed symbol
sequences. (The rule that $1$ and $5$ cannot follow each other must still be
obeyed.) \ Hence the desired degree of sensitivity to initial conditions, in
which any small perturbation of the initial condition leads to a deviation
from the given sequence, is achieved. \ \ 

\bigskip

In the paper \cite{hm2}, \ a similar result was obtained. \ There, however,
the analysis was not only for stable solutions. \ Indeed, only one of all the
solutions found in that paper is stable in the linearized sense. \ However
they are all hyperbolic. \ We believe that a similar analysis will allow us to
study linearizations around solutions $u_{2}$ \ and $u_{4}$ and prove their
isolation as well. \ However we will not attempt this analysis here. \ \ 

\subsection{\label{sec4.5}\bigskip Bifurcation in $\lambda$}

We saw in the Introduction that at some $\lambda_{b}\in(0,\lambda_{0}]$ new
periodic solutions appear. \ While we know that for $\lambda>\lambda_{0}$
there are at least five solutions with period $2\pi,$ we have not discussed
the nature of the bifurcation. \ It is possible that the number of solutions
goes from one to three at $\lambda_{b}$ and later increases to five.
\ Stability considerations do not rule this out.

\bigskip

\ One way of studying this numerically is to consider the graph of $G\left(
\alpha\right)  =u_{\alpha}^{\prime}\left(  \pi\right)  $, since we are only
considering periodic solutions with the properties $u^{\prime}\left(
0\right)  =u^{\prime}\left(  \pi\right)  =0.$ If the bifurcation is
``pitchfork'', then the graph of $G$ will qualitatively resemble that of the
function $\alpha^{3}-\mu\alpha$ \ as the parameter $\mu$ changes from negative
to positive. \ (The zeros of $G$ will not be at $\alpha=0.$) \ However a
numerical study of the function $G$ quickly suggests that this is not the
case. \ In Figure 9 we show the resulting graph of $G$ for $\varepsilon=1$
close to the bifurcation point, $\lambda=1.023.$ \ This indicates that two
pairs of solutions bifurcate at the same value of $\lambda.$ This is partly a
trivial observation, however, for the symmetry in the problem shows that
solutions other than those for which $u\left(  \frac{\pi}{2}\right)  =0$ occur
in symmetric pairs. \ The essential nature of the bifurcation can be seen by
looking only at the left branch of the bifurcation curve shown in Figure 9.

\bigskip

For comparison we also show a graph, Figure 10, of $G$ for the forcing
function $g\left(  t\right)  =\cos t-\sin2t.$ \ This is computed with
$\lambda$ a bit above the bifurcation point. \ In this problem symmetry is
lost, and solutions with $u^{\prime}=0$ at $0$ and $\pi$ are not necessarily
periodic. \ However the problem of finding solutions satisfying these two
boundary conditions is still of interest, as an example of steady-states for
the corresponding pde boundary value problem. \ The figure suggests that for
some $g$ the bifurcation may be of pitchfork type. \ We used the program xpp
of G. B. Ermentrout $\cite{erm}.$

\bigskip{%
\begin{center}
\includegraphics[
height=2.5183in,
width=4.0759in
]%
{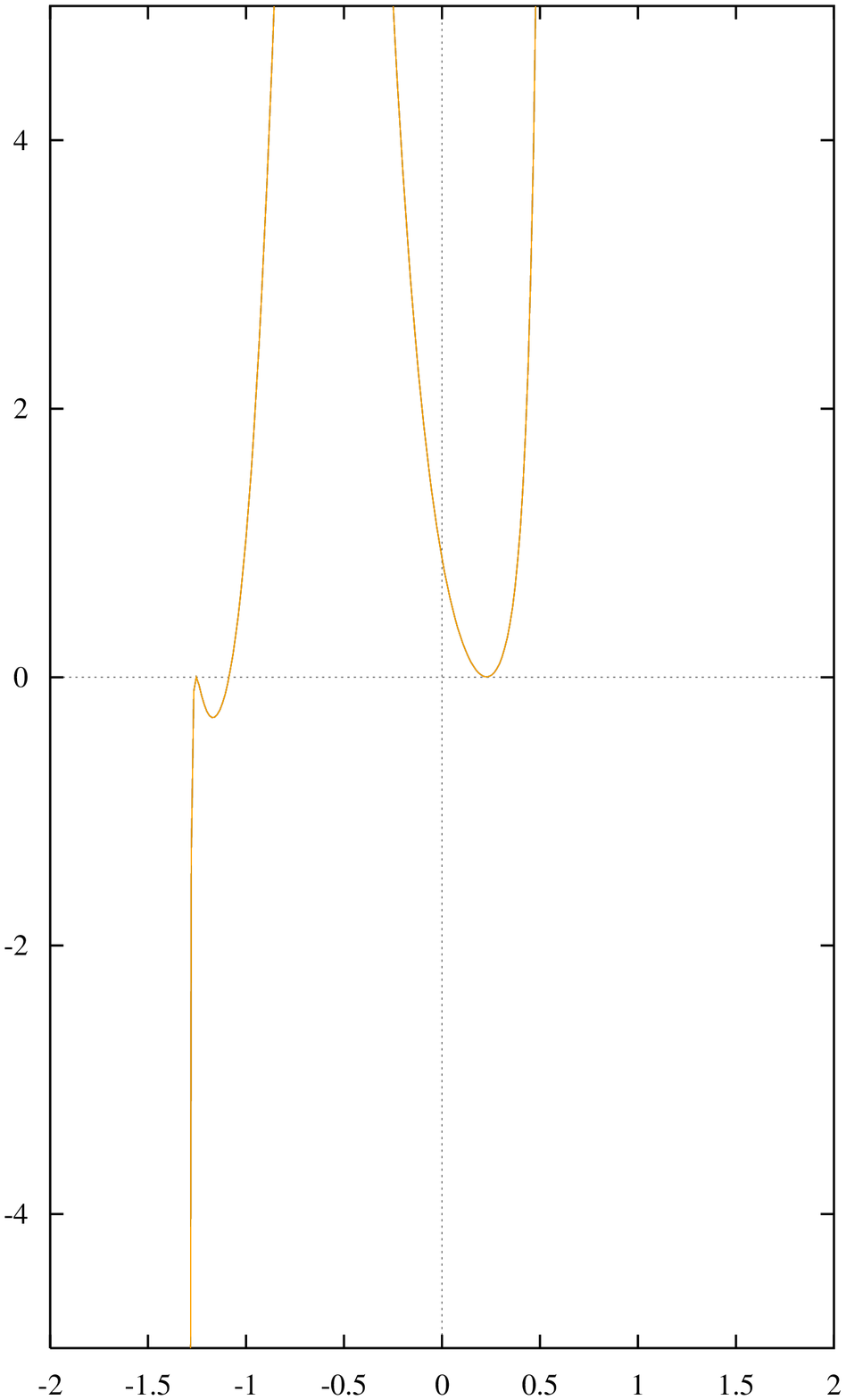}%
\\
Figure 9
\end{center}
}

{%
\begin{center}
\includegraphics[
height=2.4915in,
width=4.1675in
]%
{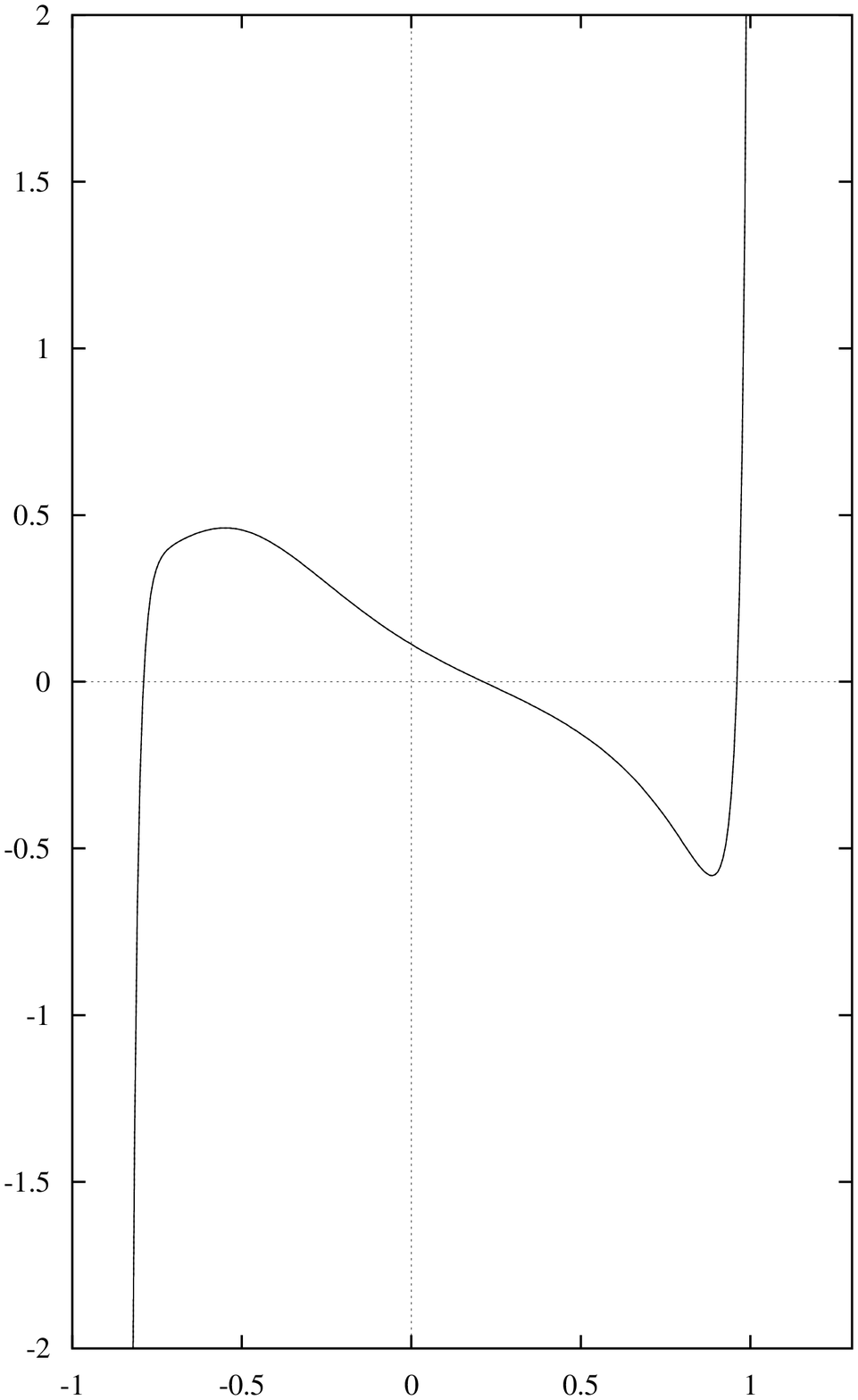}%
\\
Figure 10
\end{center}
}

\bigskip

To study the bifurcation analytically we again use $v=\frac{\partial
u}{\partial\alpha}.$ \ We are considering only the specific equation $\left(
\ref{1.1}\right)  .$ \ In order to prove that the bifurcation is not
pitchfork, we study the linearization around the antisymmetric solution
$u_{p}$. \ We determine the slope of the function $G$ \ at the antisymmetric
solution. \ Since $G\left(  \alpha\right)  =u_{\alpha}^{\prime}\left(
\pi\right)  ,$ $G^{\prime}\left(  \alpha\right)  =v_{\alpha}^{\prime}\left(
\pi\right)  $. \ But in the proof of Theorem \ref{thm8} \ we showed that when
$u=u_{p},$ $v>0$ \ on $\left[  0,\pi\right]  $ and\ $v^{\prime}\left(
\pi\right)  >0.$ This proves that the bifurcation is not of pitchfork type,
\ because the new solutions must appear at a positive distance away from
$u_{p}$.

\bigskip

However, this analysis does not eliminate the possibility of several
bifurcation points, in which solutions appear and then disappear. \ To get a
complete picture more work is required. We have seen that for $\lambda
>\lambda_{0}$ there are at least five solutions. \ We also showed that
$\lambda_{b},$ the bifurcation point, is less than $\lambda_{0}.$ But in the
following result, we hold $\lambda$ fixed and less than $\lambda_{0},$ and let
$\varepsilon$ tend to zero. \ 

\begin{theorem}
\label{thmbif1}\bigskip Suppose that $\lambda<\lambda_{0}$. \ \ Then for
sufficiently small $\varepsilon,$ $u_{p}$ is the only solution with
$u^{\prime}\left(  0\right)  =0,$ $u^{\prime}\left(  \pi\right)  =0.$
\end{theorem}

\begin{corollary}
\label{cor1}$\lim_{\varepsilon\rightarrow0}\lambda_{b}=\lambda_{0}.$
\end{corollary}

\bigskip

\begin{proof}
\ The corollary follows immediately from the Theorem, the proof of which is
more easily understood by reference to figure 1. We use the original scaling
$\left(  \ref{1.1}\right)  .$ \ For any $\delta>0$ there are $\mu>0$ and
$\nu>0$\ such that $f\left(  t,u\right)  \leq-\nu$ on the set%
\[
\Omega_{\delta,\mu}=\left\{  \left(  t,u\right)  |\left(  1-\mu\right)
\pi\leq t\leq\left(  1+\mu\right)  \pi,\ u\leq\bar{U}\left(  t\right)
-\delta\right\}  .
\]
\ \ \ Hence, for sufficiently small $\varepsilon,$ \ no periodic solution can
intersect the region $\Omega_{\delta,\mu}$ at a point where $\left(
1-\frac{\mu}{2}\right)  \pi\leq t\leq\left(  1+\frac{\mu}{2}\right)  \pi$,
\ \ since then $u$ \ would be forced below $-b$ \ before $\left|
t-\pi\right|  =\mu.$ \ \ \ So any periodic solution must have lie above
$\bar{U}\left(  t\right)  -\delta$ \ on $\left(  1-\frac{\mu}{2}\right)
\pi\leq t\leq\left(  1+\frac{\mu}{2}\right)  \pi.$ The argument used to prove
Theorem \ref{thm8} can then be used to show that for sufficiently small
$\varepsilon,$ $u_{p}$ is the only periodic solution $u_{\alpha}$ with a
maximum at $\pi.$ \ This completes the proof of the theorem and Corollary.
\end{proof}

\bigskip

We now show that for small $\varepsilon$ there is no ``reverse bifurcation''
as $\lambda$ increases from $\lambda_{b}$. In other words, as we increase
$\lambda$ from $\lambda_{b}$, there is no return to the case where $u_{p}$ is
the unique periodic solution.

\bigskip

\begin{theorem}
\label{thm3} \ For sufficiently small $\varepsilon>0,$ if $\lambda>\lambda
_{b}$ then there are at least five periodic solutions, \ $u_{1},...,u_{5}$
with $\alpha_{i}=u_{i}\left(  0\right)  <0.$ Further, $u_{i}<u_{i+1}$ for
$i=1,\cdots,4.$ \ Here $u_{3}$ is the antisymmetric solution, also denoted by
$u_{p}.$ \ The solutions $u_{1}$ and $u_{5}$ are symmetric reflections, with
$u_{5}\left(  t\right)  =-u_{1}\left(  t+\pi\right)  .$ Similarly,
$u_{4}\left(  t\right)  =-u_{2}\left(  t+\pi\right)  $. \ 
\end{theorem}

\begin{proof}
\ We have seen that $G\left(  \alpha\right)  <0$ for large negative $\alpha$
and also $G<0$ just below $\alpha_{p}.$ The existence of at least five
solutions for $\lambda>\lambda_{b}$ follows from the remarks above by showing
that if
\[
{\alpha_{1}}=\inf\left\{  \alpha|\,\,u_{\alpha}^{\prime}\left(  \pi\right)
=0\right\}  ,
\]
and if $G\left(  {\alpha_{1}}\right)  =G^{\prime}\left(  {\alpha_{1}}\right)
=0$, \ then $\frac{\partial G\left(  {\alpha_{1}}\right)  }{\partial\lambda
}>0.$ Let $u=u_{{\alpha_{1}}},$ $h\left(  t\right)  =\frac{\partial u\left(
t\right)  }{\partial\lambda}|_{\alpha={\alpha_{1}}}$, \ and as above,
$v=\frac{\partial u}{\partial\alpha}.$ \ \ We see that%
\begin{align}
\varepsilon^{2}h^{\prime\prime}  &  =\left(  3u^{2}-\lambda\right)
h-u\label{a1}\\
h\left(  0\right)   &  =h^{\prime}\left(  0\right)  =0\nonumber
\end{align}
\bigskip while $v$ \ satisfies (\ref{st1}). \ \ Therefore we obtain, for any
$t\in\left[  0,\pi\right]  ,$
\begin{equation}
hv^{\prime}-vh^{\prime}|_{t}=\frac{1}{\varepsilon^{2}}\int_{0}^{t}v\left(
s\right)  \,u\left(  s\right)  \,ds. \label{eq1}%
\end{equation}

\bigskip

\begin{lemma}
\label{lem3}For sufficiently small $\varepsilon,$ $u_{{\alpha_{1}}}<0$ on
$\left[  0,\pi\right]  .$
\end{lemma}

\begin{proof}
\ \ This follows by a slight modification of the proof of Theorem
\ref{thmbif1}. Let
\[
\hat{\Omega}_{\delta,\mu}=\left\{  \left(  t,u\right)  |\left(  1-\mu\right)
\pi\leq t\leq\left(  1+\mu\right)  \pi,\ 0\leq u\leq\bar{U}\left(  t\right)
-\delta\right\}  .
\]
\ Suppose $\lambda=\lambda_{0}.$ Then for small enough $\delta$ and $\mu,$
there is an $\varepsilon_{0}$ such that if $0<\varepsilon\leq\varepsilon_{0},$
\ then any solution which intersects $\hat{\Omega}_{\delta,\mu}$ must decrease
monotonically, in at least one direction, to below $-b,$ within a time which
is bounded over $\hat{\Omega}_{\delta,\mu}.$ \ By continuity, this is also
true for $\lambda$ sufficiently close to $\lambda_{0.}$ \ As in the proof of
Theorem \ref{thm8}, no solution other than $u_{p}$ can have a positive maximum
at $\pi$ and not intersect $\hat{\Omega}_{\delta,\mu},$ proving the result. \ 
\end{proof}

\bigskip

\ \ Also, \ if $u=u_{{\alpha_{1}},}$ then $v^{\prime}\left(  \pi\right)
=G^{\prime}\left(  {\alpha_{1}}\right)  =0.$ \ Therefore $v\left(  \pi\right)
\neq0.$ \ Suppose that $v\left(  \pi\right)  <0.$ \ Then for $\alpha$ slightly
lower than ${\alpha_{1}},$ $u_{\alpha}\left(  \pi\right)  >u_{{\alpha_{1}}%
}\left(  \pi\right)  .$ Let
\[
\beta=\inf\left\{  \alpha\,|\,\,u_{\alpha}\left(  t\right)  >u_{{\alpha_{1}}%
}\left(  t\right)  \text{ for some }t\in(0,\pi]\right\}  .
\]
Then $\beta$ is well-defined \ and $-\infty<\beta<{\alpha_{1}}$. \ Suppose
that $u_{\beta}\left(  t_{0}\right)  =u_{{\alpha_{1}}}\left(  t_{0}\right)  $
for some $t_{0}\in\left(  0,\pi\right)  .$ Since different solutions cannot be
tangent, we must have $u_{\beta}\left(  t\right)  >u_{{\alpha_{1}}}\left(
t\right)  $ for some $t\in\left(  0,\pi\right)  ,$ but this contradicts the
definition of $\beta.$ Therefore $u_{\beta}\left(  \pi\right)  =u_{{\alpha
_{1}}}\left(  \pi\right)  $ and $u_{\beta}^{\prime}\left(  \pi\right)
>u_{{\alpha_{1}}}\left(  \pi\right)  =0.$ \ But in this case, we can lower
$\alpha$ further, until we find a $\gamma<\beta$ with $u_{\gamma}^{\prime
}\left(  \pi\right)  =0.$ \ This contradicts the definition of ${\alpha_{1}}.$
A similar argument shows that $v>0$ on $[0,\pi]$. \bigskip\ 

Therefore, when $u=u_{{\alpha_{1}}},$ \ $v\left(  \pi\right)  >0.$ Then
$\left(  \ref{eq1}\right)  $ and Lemma \ref{lem3} show that $\frac{\partial
u^{\prime}\left(  \pi\right)  }{\partial\lambda}=h^{\prime}\left(  \pi\right)
>0.$ Thus, \ for $\lambda$ just above $\lambda_{b}$ there are at least five
solutions (using symmetry). \ Further, there cannot be a decrease to fewer
than five as $\lambda$ increases further, for at any point where $G=G^{\prime
}=0$ \ we would again get $h^{\prime}\left(  \pi\right)  >0.$ \ 

\bigskip

\ 

To complete the proof of Theorem \ref{thm3} it is convenient to let
$\alpha_{2}=\sup\left\{  \alpha<\alpha_{p}|u_{\alpha}^{\prime}\left(
\pi\right)  =0\right\}  .$ (Numerically, it appears there is only one $2\pi
$-periodic solution between $u_{1}$ and $u_{p}.)$ Then reflection and
translation of $u_{\alpha_{1}}$and $u_{\alpha_{2}}$ by the transformation
$u\left(  t\right)  \rightarrow-u\left(  t+\pi\right)  $ give the additional
two asymmetric solutions. \ 

\bigskip

Our construction implies that $u_{i}\left(  0\right)  <u_{i+1}\left(
0\right)  .$ \ From Proposition \ref{thm02} and the way we define $\alpha_{2}$
we have $u_{\alpha}^{\prime}\left(  \pi\right)  <0$ for $\alpha_{2}%
<\alpha<\alpha_{p},$ and in this range and close enough to $\alpha_{p,}$
$u_{\alpha}<u_{p}$ \ (This is because when $u=u_{p,}$ $v>0$ on $\left[
0,\pi\right]  .)$ \ \ Suppose, however, that for some $\alpha\in\lbrack
\alpha_{2},\alpha_{p}),$ $u_{\alpha}\left(  t\right)  \geq u_{p}\left(
t\right)  $ for some $t\in\left[  0,\pi\right]  .$ Let $\check{\alpha}%
=\sup\left\{  \alpha\in\lbrack\alpha_{2},\alpha_{p})|u_{\alpha}\left(
t\right)  =u_{p}\left(  t\right)  \text{ for some }t\in\lbrack0,\pi].\right\}
$. \ By the same argument as above we show that $u_{\check{\alpha}}\left(
\pi\right)  =u_{p}\left(  \pi\right)  $ \ and $u_{\check{\alpha}}^{\prime
}\left(  \pi\right)  >0.$ \ \ Hence there is an $\alpha\in\left(
\check{\alpha},\alpha_{p}\right)  $ with $u_{\alpha}^{\prime}\left(
\pi\right)  =0.$ \ But this contradicts the definition of $\alpha_{2}.$

\bigskip

This proves that $u_{2}<u_{p}.$ \ \ The proof that $u_{1}<u_{2}$ is similar,
and the construction of $u_{4}$ and $u_{5}$ by reflection and translation
implies the remaining order relations, namely, $u_{p}=u_{3}<u_{4}<u_{5}.$
\ This completes the proof of Theorem \ref{thm3}.\bigskip
\end{proof}

\bigskip

\section{\label{sec5}Conclusion}

\bigskip

\subsection{\label{sec5.1}How special is the cosine?}

\bigskip\ \ \ 

General non-symmetric forcing functions are a topic for further study, but a
few things are easy to see. \ First, the procedure in Theorems \ref{thm4} and
\ref{thm5}\ for obtaining some sort of chaotic behavior will carry over to a
large variety of forcing functions. \ Also, on a finite interval the technique
will give many steady states for the problem $\left(  \ref{1.2}\right)
-\left(  \ref{1.3}\right)  .$ \ The only requirement is the existence of some
set of functions $w_{k}$ whose graphs form ``fingers'' \ pointing up and down
alternately in a way similar to that in figure 2. \ We will not try to
formulate a precise result here.

\bigskip

We did use symmetry essentially to obtain the existence of the solution
$u_{p}.$ \ It is here that we have found significant differences between the
problem of finding steady-states for $\left(  \ref{1.2}\right)  -\left(
\ref{1.3}\right)  $ on a finite interval and the problem of finding bounded
solutions, and chaos, on an infinite interval. \ It is easy to find a third
solution to the boundary value problem on a finite interval using shooting.
\ From there one can go on to find many other solutions, both stable and
unstable. \ But dealing with chaos on an infinite interval seems different,
and it is for that reason that we concentrated on the particular equation
$\left(  \ref{1.1}\right)  $. \ Having done so, it is natural to make use of
symmetry to obtain simpler proofs in some cases where a shooting method may
apply even without symmetry. \ 

\bigskip

We note, however, that the proof of Theorem \ref{thm4}, which does not use
symmetry, includes the existence of a solution which does not intersect any of
the $w_{k}.$ \ \ This solution is the stable solution of \cite{ampp}, \ and
could play the role of $u_{p}$ in a study which uses shooting but does not use symmetry.

\bigskip

The proof of stability of the solution $u_{p}$ (Theorem \ref{thm8}) does not
use symmetry, but does use the fact that the derivative of cosine is negative
in $(0,\pi)$. We do not believe that this property should be crucial to
obtaining an elementary proof of linearized stability, and we expect that our
proof will extend to more general forcing functions, but this is the subject
of future investigation. The assumptions in \cite{ampp} imply that the forcing
function $g$ is monotone in a neighborhood of any of its zeros, and this
should allow our proof to go through. However the result will probably be
weaker in that $S_{\delta}$ will be restricted more.

\bigskip We make use of Lemma \ref{alem7.2}, about the monotonicity of the
maxima and minima in $(0,\pi)$, a number of times, and this also relies on the
monotonicity of cosine in this interval. Our future work will include a study
of how results will change without this condition, but we expect that
extensions will be possible.

\bigskip

The results about bifurcation as $\lambda$ increases, however, may change
radically if $g$ has less symmetry. \ For example, if $g$ has period $2\pi,$
\ and is symmetric around $\pi,$ but not anti-symmetric around $\frac{\pi}%
{2},$ then the solutions will not appear in symmetric pairs. \ The technique
we have used may still be able to prove that saddle-node bifurcations occur,
at least in a nearly-symmetric situation, but there may be two of them,
producing first two new solutions and then two more. \ This is easily seen in
numerical simulations. \ As the deviation from full symmetry increases, it
appears from some brief numerical experiments that other possibilities exist,
and we hope to explore this further.

\bigskip

\subsection{\label{sec5.2}Summary of main points.}

\bigskip

In Section 1 we introduce the problem and relate it to some previous work. \ 

\bigskip

In Section 2 we give preliminary results valid for all $\varepsilon>0.$ Some
of these are for a more general forcing function $g\left(  t\right)  .$ The
main points are that for $\lambda\leq0$ \ there is a unique bounded solution,
which is periodic, while for $\lambda>\lambda_{0}$ there are at least three
solutions. \ Therefore a bifurcation takes place, and some computations are
given in section 4.5 indicating that this can be of different types depending
on $g\left(  t\right)  .$

\bigskip

\ \ In section 3 \ we give results valid for a specific range of $\varepsilon$
\ and $\lambda.$ \ \ The main hypothesis is Condition \ref{A}, \ the existence
of ``spikes'', ( the $w_{k}),$ which are solutions tending to $\pm\infty$ in
both directions. \ Solutions are characterized by which of the $w_{k}$ they
intersect. \ Section 3.1 gives the proof of the 1:1 correspondence with
certain sequences if Condition \ref{A} holds. \ It would be possible to
rephrase this result to give a natural correspondence with sequences of three
symbols, corresponding to the three solutions which we later labeled
$u_{1},u_{3}\left(  =u_{p}\right)  ,$ and $u_{5}.$ \ In 3.2 we give a brief
discussion of ``kneading theory'' in our context. \ In Section 3.3 \ the
symbolic dynamics is extended to sequences of five symbols. In section 3.4
\ Condition \ref{A} is verified, first for ``sufficiently small''
$\varepsilon,$ where no analysis is required, and then for larger $\varepsilon.$

\bigskip

The results in Section 3 do not include uniqueness or stability, and there is
only a limited sensitivity to initial conditions demonstrated. \ \ They do not
depend at all on symmetry, and indeed, the techniques will yield a weak form
of chaos for a wide variety of forcing functions, including non-periodic forcing.

\bigskip

Section 4 contains a variety of results, all proved for sufficiently small
$\varepsilon$ with no estimate on the range of $\varepsilon$ for which they
hold. \ In 4.1 some results are given about asymptotic behavior of solutions
as $\varepsilon\rightarrow0.$ \ In 4.2 further periodic solutions are found,
including solutions with multiple internal layers and with a ``down-jump''
near $\frac{\pi}{2}$, \ in contrast to the solution $u_{p}$ which jumps upward
at $\frac{\pi}{2}$. \ \ \ In 4.3 we give a proof of the stability of $u_{p}$
using classical ode methods, and extend the result on uniqueness of
\cite{ampp} \ a bit by obtaining a larger region in which it is unique. \ The
uniqueness proof is also more direct than that in \cite{ampp}.

\bigskip

In 4.4 we discuss sensitivity with respect to initial conditions. In 4.5 \ we
consider the bifurcation problem in $\lambda$ for $2\pi$-periodic solutions of
$\left(  \ref{1.1}\right)  $.

\bigskip

Finally, in Section 5.1, \ we discuss the role of the specific cosine forcing
function in our results, and conjecture that its symmetry and monotonicity
properties may not be essential except in the bifurcation analysis of section
4.6 .

\bigskip We wish to express our appreciation to Professor H. Matano for
sharing his thoughts on this problem. \ His remarks are cited in more detail
at the end of section \ref{sec4.2}, \ together with a citation of recent work
of his coworker, Dr. K. Nakashima. \

\bigskip

\bigskip

\bigskip
\end{document}